\documentclass[11pt]{article}

\usepackage{amsmath,amssymb,amsthm}
\usepackage{multirow}
\usepackage{enumitem}
\usepackage{hyperref}
\usepackage{xcolor}
\usepackage{tikz}
\newcommand{\R}{\mathbb{R}}
\newcommand{\Om}{\Omega}
\newcommand{\B}{\mathcal{B}}
\newcommand{\Nc}{\mathcal{N}_c}
\newcommand{\vx}{\vec{x}}
\newcommand{\F}{\mathcal{F}}
\DeclareMathOperator*{\Minimize}{Minimize}
\DeclareMathOperator*{\Maximize}{Maximize}

\newtheorem{teo}{Theorem}[section]
\newtheorem{lem}[teo]{Lemma}
\newtheorem{prop}[teo]{Proposition}
\newtheorem{coro}[teo]{Corollary}
\theoremstyle{definition}
\newtheorem{defi}[teo]{Definition}
\newtheorem{assump}{Assumption}
\newtheorem{rem}[teo]{Remark}

\begin{document}

\title{Augmented Lagrangian methods with telescopic subproblems for MPCCs and cardinality
constraints with general constraints: constraint-preserving swaps and B-stationarity\thanks{This
  work was supported by FAPESP (grants 2022/05803-3 and 2023/08706-1) and CNPq (grants
  302073/2022-1 and 305058/2024-0).}}

\author{
E. G. Birgin\thanks{Department of Computer Science, Institute of Mathematics, Statistics, and Computer Science (IME),
  University of S\~ao Paulo, Rua do Mat\~ao, 1010, Cidade Universit\'aria, 05508-090,
  S\~ao Paulo, SP, Brazil. e-mail: egbirgin@ime.usp.br}
\and
J. M. Mart\'{\i}nez\thanks{Department of Applied Mathematics, Institute of Mathematics,
  Statistics, and Scientific Computing (IMECC), State University of Campinas, 13083-859
  Campinas SP, Brazil. e-mail: martinez@ime.unicamp.br}}

\date{November 5, 2026}

\maketitle

\begin{abstract}
We consider optimization problems with general equality and inequality constraints and
partitioned cardinality constraints, which include complementarity constraints. In the
augmented Lagrangian method that we consider, the bounds and the cardinality constraints are
kept in the subproblems, which are solved by the Telescopic Method. The limit points of this
method are B-stationary, without constraint qualifications, if the penalty parameters remain
bounded, but not in general. The reason is that the swaps of the Telescopic Method pass through
points at which the constraints have different values, and the corresponding variation of the
penalty term grows with the penalty parameter. With unbounded penalty parameters and a
derivative test outside the current face, the feasible limit points at which a linear
independence constraint qualification holds are M-stationary for MPCCs. We modify the swaps, by
means of a scaled Newton restoration, so that they preserve the values of the constraints, and
we prove that every feasible limit point at which this constraint qualification holds is then
B-stationary (S-stationary for MPCCs), whatever the behavior of the penalty parameters.
Counterexamples show that the modification is necessary, even with exact solutions of the
subproblems, and that the constraint qualification cannot be weakened to one of Mangasarian--
Fromovitz type. In the MacMPEC collection, the method found certified S- or B-stationary points
in 185 of the 193 instances, more than Algencan and the regularization method of Scholtes. In
experiments with cardinality-constrained portfolios and sparse principal component analysis,
the constraint-preserving swaps improved the original swaps and the methods applied to
complementarity reformulations.
\end{abstract}



\section{Introduction} \label{motiv}

In \cite{bmtel} we introduced the Telescopic Method for the problem
\begin{equation} \label{pcc}
\Minimize f(x) \mbox{ subject to } x \in \Om,
\end{equation}
where $\Om$ is defined by bounds and partitioned cardinality constraints (PCC). This means that the variables
are divided into blocks and that, in each block, at most a given number of variables may be
non-null. Bound-constrained mathematical programs with complementarity constraints (MPCCs) are
the particular case in which each block is a complementarity pair, with at most one non-null
variable. Every limit point generated by the method is B-stationary,
without constraint qualifications. In this paper we consider additional constraints
$h(x) = 0$ and $c(x) \leq 0$.

Without constraint qualifications, we do not know any method that guarantees B-stationarity
for problems with general nonlinear constraints, not even in smooth nonlinear programming. In that generality, augmented Lagrangian,
sequential quadratic programming, and interior point methods guarantee sequential optimality
conditions such as AKKT \cite{ahm11} (see \cite{krs21b} for cardinality constraints), which may hold at points that are not B-stationary.
For example, for $\Minimize -x_1 - x_2$ subject to $x_1 x_2 = 0$ and $x \geq 0$, AKKT holds
at $x^* = 0$ (take $x^k = (1/k, 1/k)$ and the multiplier $k$ for $x_1 x_2 = 0$), which is
not B-stationary.

For MPCCs with general constraints, regularization methods, such as
Schol\-tes' relaxation \cite{sch01}, and penalty methods converge, under MPEC constraint qualifications, to C- or M-stationary points, which need not be B-stationary. With
inexact solutions of the subproblems, even these guarantees may be weakened \cite{ks15}.
Relaxation methods are compared, theoretically and on MacMPEC, in \cite{hks13}. SQP methods converge locally superlinearly near strongly stationary points
\cite{flrs06}, and interior-penalty methods converge globally to strongly stationary points
under standard assumptions \cite{llcn06}. Augmented Lagrangian methods applied to the
nonlinear programming reformulation of MPCCs were analyzed in \cite{isu12,ass18} and, for
cardinality-constrained problems (see also \cite{be13,bks}), in \cite{krs21}, where
convergence to strongly stationary points is proved under a problem-tailored quasinormality
condition. Guo and Deng \cite{gd22} keep the complementarity constraints, written
with slack variables, in the subproblems of an augmented Lagrangian method, solve the
subproblems by a nonmonotone projected gradient method, and obtain M-stationary limit points
under MPCC-RCPLD. For this method, Andreani, da Rosa and Secchin \cite{ars25} proved the
boundedness of the multipliers under quasi-normality-type conditions, and observed
that keeping the complementarity constraints in the subproblems improves the numerical
stability. Closest to our approach, Jia, Kanzow, Mehlitz and Wachsmuth \cite{jkmw23} keep
nonconvex sets with inexpensive projections, such as complementarity or cardinality
constraints, in the subproblems of an augmented Lagrangian method, and obtain M-stationary
limit points under an asymptotic regularity condition. The straightforward augmented Lagrangian method with telescopic subproblems also keeps the nonconvex constraints in the subproblems and, if their stopping criterion includes a derivative test outside the current face, its feasible limit points are M-stationary under MPEC-LICQ (Section~\ref{alsecunb}). In the method proposed in this paper,
the subproblems are also solved by the Telescopic Method, but its swaps are modified so that
they preserve the values of the constraints, and the limit points are B-stationary
(S-stationary for MPCCs) under PCC-LICQ (MPEC-LICQ for MPCCs). Active-set and hybrid methods, such as the one of Lin and Fukushima
\cite{lf06}, assume MPEC-LICQ and use the signs of the multipliers to change the active set
or to certify S-stationarity. The method of Nurkanovi\'c and Leyffer \cite{nl26} solves a
sequence of nonlinear programs associated with branches of the feasible set (BNLPs) and of
linear programs with equilibrium constraints (LPECs). The LPECs either certify B-stationarity
or provide a new branch. Global convergence to B-stationary points is proved assuming
MPEC-MFCQ at every feasible point and that the BNLP solver returns stationary points with a uniform
sufficient decrease. For bound-constrained MPCCs, the method of Kirches, Larson, Leyffer and
Manns \cite{kllm} converges to B-stationary points without constraint qualifications. Other
constraints may be handled by an augmented Lagrangian method, but, as in
Section~\ref{alsec}, the B-stationarity of the subproblems is not inherited by the limit
points.
To our knowledge, no method guarantees B-stationarity for PCC problems with general
constraints. For MPCCs, the closest result is the one of \cite{nl26}, which requires
MPEC-MFCQ at every feasible point and solves LPECs, whereas our result requires PCC-LICQ only at the
limit point, uses only swaps, and applies to general partitioned cardinality constraints. The
two assumptions are not comparable.

The contributions of this paper are the following.
\begin{itemize}[itemsep=2pt,parsep=0pt,topsep=4pt]
\item The straightforward safeguarded augmented Lagrangian method \cite{abms,bmbook} with telescopic subproblems
  (Algorithm~\ref{alsec}.1). Its limit points are B-stationary, without constraint
  qualifications, if the penalty parameters remain bounded (Theorem~\ref{teoal}), but not in
  general (counterexample of Proposition~\ref{propcex}). If the penalty parameters tend to infinity, its feasible limit points at which PCC-LICQ holds are M-stationary for MPCCs (Theorem~\ref{teounb} and Corollary~\ref{corounb}), thanks to condition (c) of the stopping criterion of the subproblems, which cannot be dropped (counterexample of Proposition~\ref{propcexm}).
\item A characterization of B-stationarity for PCC problems with general constraints
  (Proposition~\ref{propbstat}).
\item The constraint-preserving swap (Algorithm~\ref{cps}.2), based on a scaled Newton
  restoration with a computable tolerance (Lemma~\ref{lemrest}).
\item The modified augmented Lagrangian method (Algorithm~\ref{cps}.3), whose feasible limit
  points at which PCC-LICQ holds are B-stationary, also when the penalty parameters tend to
  infinity and for the quadratic penalty method (Theorem~\ref{teoconv}).
\item The modification that gives rise to Algorithm~\ref{cps}.3 is necessary, even under MPEC-LICQ and with exact solutions of the
  subproblems (counterexamples of Propositions~\ref{propcex2}, \ref{propcex3}
  and~\ref{propcex4}), and PCC-LICQ cannot be replaced by a condition of
  Mangasarian--Fromovitz type (counterexample of Proposition~\ref{propmfcq}).
\item Numerical confirmation of the counterexamples (Section~\ref{numer}), and experiments with
  the 193 instances of the MacMPEC collection (Section~\ref{macmpec}), in which Algorithm~\ref{cps}.3
  certified more solutions than Algencan \cite{abms,bmbook} applied to the nonlinear
  programming reformulation and the regularization method of Scholtes \cite{sch01}.
  In the 191 problems used by Nurkanovi\'c and Leyffer \cite{nl26}, the success rate of Algorithm~\ref{cps}.3 was
  97.4\%, while the highest success rate reported in \cite{nl26}, obtained by MPECopt, was
  94.24\%.
\item Experiments with cardinality-constrained portfolios and sparse principal component
  analysis (Section~\ref{cardexp}). In these problems, the constraint-preserving swaps
  improved the original swaps. In most problems, Algorithm~\ref{cps}.3 also found better points
  than Algencan and the regularization method of Scholtes, applied to the
  complementarity-type reformulation of \cite{bks}. In sparse principal component analysis,
  it also found better points than the truncated power method of \cite{yz13}.
\end{itemize}

Section~\ref{overview} gives an overview of the paper, and Section~\ref{probl} describes the
problem and the stationarity concepts. Section~\ref{alsec} presents the augmented Lagrangian
method with telescopic subproblems, its convergence when the penalty parameters remain
bounded, the counterexample, and the properties of the limit points when the penalty parameters are unbounded. Constraint-preserving swaps are introduced in
Section~\ref{cps}, and the augmented Lagrangian method that uses them is defined in
Section~\ref{algal}. Section~\ref{conv} contains its convergence analysis and shows that the
modification is necessary and that MPEC-LICQ cannot be replaced by MPEC-MFCQ. Numerical
experiments are reported in Section~\ref{numexp}, and conclusions are given in
Section~\ref{concl}.

\medskip
\noindent\textbf{Notation.} Throughout the paper, $\|\cdot\|$ is the Euclidean norm, $\|\cdot\|_\infty$ is the sup norm and, for a vector $v$,
$\|v\|_0$ is the number of non-null entries of $v$. For $t \in \R$, $(t)_+ = \max\{0, t\}$,
applied componentwise to vectors. $P_C$ denotes the Euclidean projection onto a closed convex
set $C$, $T_C(x)$ is the (Bouligand) tangent cone of a set $C$ at $x \in C$, and $e^{i,\nu}$
(or $e^v$) is the canonical vector associated with the variable $x_{i,\nu}$ (or $x_v$). For a differentiable function $\phi$, $\partial_v \phi(x) = \partial \phi(x)/\partial x_v$ (or $\partial_{i,\nu} \phi(x) = \partial \phi(x)/\partial x_{i,\nu}$).

\section{Overview} \label{overview}

We consider the problem
\begin{equation} \label{pccnl}
\Minimize f(x) \mbox{ subject to } h(x) = 0, \; c(x) \leq 0, \; x \in \Om,
\end{equation}
where $f$, $h$ and $c$ are smooth and $\Om$ is defined by bounds on the variables and by
partitioned cardinality constraints. This means that the variables are divided into blocks and, in each
block, at most a given number of variables are allowed to be non-null. Complementarity constraints
$0 \leq y \perp z \geq 0$ are a particular case of partitioned cardinality constraints, in which each pair $(y_i, z_i)$ is a block
with at most one non-null variable. Binary variables are also a particular case of partitioned cardinality constraints (Section~\ref{probl}).

In \cite{bmtel} we considered problem~(\ref{pccnl}) without the constraints $h(x) = 0$ and
$c(x) \leq 0$, and we introduced the Telescopic
Method. For this purpose, the set $\Om$ is divided into faces. Each face is defined by the variables that are allowed to
be non-null on it and is a box. The Telescopic Method alternates iterations of a
bound-constrained algorithm on the current face with \emph{swaps}, in which a variable of the
face leaves it (it is set to zero) and a variable outside the face enters it, moving
away from zero along a one-dimensional search. A swap is accepted if it decreases $f$
sufficiently, and then the method continues on the new face. Every limit point of the
sequence generated by the method is B-stationary, without constraint qualifications.

To handle the additional constraints, it is natural to use an augmented Lagrangian method,
as Algencan \cite{abms,bmbook}, in which the constraints $h(x) = 0$ and $c(x) \leq 0$ are
penalized and the set $\Om$ is kept in the subproblems. At each outer iteration, the subproblem (minimization
of the augmented Lagrangian on $\Om$) is solved by the Telescopic Method. We will show that, if
the penalty parameters remain bounded, every limit point is B-stationary, without
constraint qualifications (Section~\ref{alsec}). However, if the penalty parameters tend to
infinity, the limit points may fail to be B-stationary, even if every subproblem is solved
exactly (counterexample of Proposition~\ref{propcex}). The reason is that a swap could move the point to places
where the constraints have different values than the ones at the current point, in such a way
that, due to a large penalty parameter, a promising swap could be rejected, in spite of having
decreased the objective function. When the penalty parameters tend to infinity and the stopping criterion of the subproblems also looks outside the current face, the feasible limit points are M-stationary for MPCCs under MPEC-LICQ (Section~\ref{alsecunb}), but not necessarily B-stationary. Moreover, in this case, every
limit point is stationary for the minimization of the infeasibility measure
$\|h(x)\|^2 + \|c(x)_+\|^2$ on a face of $\Om$ (Proposition~\ref{propinf}). In particular,
this necessarily occurs if the feasible set is empty.

In order to remedy the inconvenience described above, we introduce \emph{constraint-preserving
swaps} (Section~\ref{cps}). In a constraint-preserving swap, the leaving variable
is set to zero and the entering variable is moved away from zero by a backtracking search on its
value. For each trial value, the other free variables of the face are corrected by a
Newton-type method, so that the resulting point $w$ satisfies $h(w) \approx h(x^k)$. If $L_k$
is the augmented Lagrangian of the $k$-th subproblem and $h(w) = h(x^k)$, then
\[
L_k(w) - L_k(x^k) = f(w) - f(x^k)
\]
for all values of the penalty parameter and of the multipliers. Therefore, decreases of $L_k$
along such swaps correspond to decreases of $f$, and the acceptance of promising swaps is not
impaired by large penalty parameters.
We will prove that, under a linear independence constraint qualification at the
limit point, the correction is well defined and small. Then the analysis of the Telescopic
Method can be repeated with the Lagrangian of problem~(\ref{pccnl}) in the place of $L_k$. With this
modification, every feasible limit point at which this constraint
qualification holds is B-stationary, even if the penalty parameters tend to infinity
(Sections~\ref{algal} and~\ref{conv}). The stationarity of the infeasible limit points for the infeasibility measure (Proposition~\ref{propinf}) is also preserved. Finally, we show that the constraint-preserving swap is necessary and
that the constraint qualification cannot be weakened to a Mangasarian--Fromovitz type
condition (Section~\ref{nec}).

\section{Problem and stationarity} \label{probl}

We consider problem~(\ref{pccnl}), where $f: \R^n \to \R$, $h: \R^n \to \R^m$, $c: \R^n \to \R^q$, and $\Om$ is defined by
bounds and partitioned cardinality constraints, as follows. The vector of variables is
partitioned into disjoint blocks $\vx_\nu \in \R^{n_\nu}$, $\nu = 1, \ldots, N$, with $n_1 + \cdots + n_N = n$, and we write $x_{i,\nu}$
for the $i$-th variable of block $\nu$ and $x_v$, $v = 1, \ldots, n$, when the block
structure is not needed. Given a box $\B = [a, b]$ and integers $p_\nu \in \{1, \ldots, n_\nu\}$,
\[
\Om = \{ x \in \B \mid \|\vx_\nu\|_0 \leq p_\nu, \; \nu = 1, \ldots, N \},
\]
that is, at most $p_\nu$ variables of block $\nu$ may be non-null. We denote by
$\Nc = \{\nu \mid p_\nu < n_\nu\}$ the set of blocks in which the cardinality constraint is
effective. Mathematical programs with complementarity constraints (MPCCs) are an important
particular case, which we use in the examples and in the experiments, in which every $\nu \in \Nc$ is
a complementarity pair, with $n_\nu = 2$, $p_\nu = 1$ and $\vec{a}_\nu = 0$.
Problem~(\ref{pccnl}) also includes binary variables, because $x_j \in \{0, 1\}$ is equivalent to
$x_j + \bar x_j = 1$, where $\bar x_j$ is an additional variable and $(x_j, \bar x_j)$ is a
complementarity pair, that is, a block with $n_\nu = 2$ and $p_\nu = 1$, with bounds $[0, 2]$ for
both variables. (With the bounds $[0, 1]$, the variable that is equal to one would be at its
upper bound, and PCC-LICQ, defined below, would fail at every feasible point.) All the results
of this paper apply to this reformulation.

Writing $c(x) + s = 0$ with slack variables
$s \in [0, s_{\max}]$, where $s_{\max} \geq \max_{x \in \B} \|(-c(x))_+\|_\infty$,
problem~(\ref{pccnl}) becomes a problem of the same form with equality constraints only, in
which the slacks form an additional block $\nu$ with $p_\nu = n_\nu$, so that $\nu \notin \Nc$. Its feasible set is
the graph of $x \mapsto -c(x)$ on the feasible set of~(\ref{pccnl}), so that B-stationarity
is preserved, and so are the constraint qualifications (see the comment after
Definition~\ref{defpccl}). Therefore, except in Section~\ref{alsec}, where the inequality
constraints are handled directly by the augmented Lagrangian method, as in Algencan, we
consider
\begin{equation} \label{mpcc}
\Minimize f(x) \mbox{ subject to } h(x) = 0, \; x \in \Om,
\end{equation}
where, with some abuse of notation, $x$ includes the slack variables, $h: \R^n \to \R^m$
includes the constraints $c(x) + s = 0$, and we denote by $X = \{x \in \Om \mid h(x) = 0\}$
its feasible set. We denote by
$J(x) \in \R^{m \times n}$ the Jacobian of $h$, by $J_v(x)$ its $v$-th column and, for
$S \subseteq \{1, \ldots, n\}$, by $J_S(x)$ the submatrix formed by the columns with indices in $S$.

We now define the faces of $\Om$ and the projected derivatives. A \emph{face index} is an $N$-tuple $I = (I_1, \ldots, I_N)$ such that
$I_\nu \subseteq \{1, \ldots, n_\nu\}$ and $\# I_\nu = p_\nu$ for all $\nu$. The
corresponding \emph{face} is the box
\[
F_I = \{ x \in \B \mid x_{i,\nu} = 0 \mbox{ whenever } i \notin I_\nu, \; \nu = 1, \ldots, N \},
\]
and $\Om$ is the union of the finitely many faces. For a continuously differentiable function
$\phi$, $x \in \B$, $\nu = 1, \ldots, N$ and $i = 1, \ldots, n_\nu$, the \emph{projected
derivative} is
\begin{equation} \label{gdef}
g^\phi_{i,\nu}(x) = P_{[a_{i,\nu}, b_{i,\nu}]}\left( x_{i,\nu} - \frac{\partial \phi}{\partial x_{i,\nu}}(x) \right) - x_{i,\nu},
\end{equation}
which vanishes if and only if $x_{i,\nu}$ satisfies the first-order condition for minimizing
$\phi$ with respect to $x_{i,\nu} \in [a_{i,\nu}, b_{i,\nu}]$, the other variables being
fixed. For $x \in F_I$, the projected gradient $P_{F_I}(x - \nabla \phi(x)) - x$ is the vector $g^\phi_I(x) \in \R^n$ with components $g^\phi_{i,\nu}(x)$ if $i \in I_\nu$ and $0$ otherwise.

We use the stationarity concept of \cite{bmtel}, given in Definition~\ref{defstat}. By
Proposition~\ref{propstatb}, it is equivalent to B-stationarity on $\Om$. Assumption~\ref{assum},
stated after that proposition, is used in the rest of the paper.

\begin{defi} \label{defstat}
A point $x^* \in \Om$ is \emph{stationary for the minimization of $\phi$ on $\Om$} if, for
all $\nu = 1, \ldots, N$,
\begin{equation} \label{s1}
x^*_{i,\nu} \neq 0 \;\Rightarrow\; g^\phi_{i,\nu}(x^*) = 0, \quad i = 1, \ldots, n_\nu,
\end{equation}
and
\begin{equation} \label{s2}
\|\vx^*_\nu\|_0 < p_\nu \;\Rightarrow\; g^\phi_{i,\nu}(x^*) = 0 \mbox{ for all } i = 1, \ldots, n_\nu.
\end{equation}
\end{defi}

Condition~(\ref{s1}) concerns the non-null variables, whereas condition~(\ref{s2}) says that, if a
block has room for one more non-null variable, then no variable of the block can decrease
$\phi$ by moving away from zero. In the complementarity case, if both components of a pair
vanish, (\ref{s2}) says that both partial derivatives are non-negative.

\begin{prop} \label{propstatb}
A point $x^* \in \Om$ is stationary for the minimization of $\phi$ on $\Om$ if and only if
$\nabla \phi(x^*)^T d \geq 0$ for all $d \in T_\Om(x^*)$.
\end{prop}

\begin{proof}
See Proposition~2.3 of \cite{bmtel}, with $\phi$ in the place of $f$.
\end{proof}

\begin{assump} \label{assum}
For all $v$, $a_v \leq b_v$ and, if $v$ is in a block $\nu \in \Nc$, also $a_v \leq 0 \leq b_v$ and $a_v < b_v$. The box $\B$ is bounded. The functions $f$ and $h$ are
continuously differentiable on an open set containing $\B$, and $\nabla f$ and $J$ are
Lipschitz continuous on $\B$.
\end{assump}

The constraint qualification of our main result involves the set of \emph{free variables} at
$x^* \in \B$,
\[
\F(x^*) = \{ v \in \{1, \ldots, n\} \mid a_v < x^*_v < b_v, \mbox{ and } x^*_v \neq 0 \mbox{ if } v
\mbox{ belongs to a block } \nu \in \Nc \}.
\]

\begin{defi} \label{defpccl}
We say that \emph{PCC-LICQ} holds at $x^* \in X$ if the matrix $J_{\F(x^*)}(x^*)$ has
rank~$m$.
\end{defi}

For example, let $N = 1$, $n = 5$, $p_1 = 3$, $m = 2$, and $x^* = (4, 0, 0, 2, 0)^T$, with $a_v < x^*_v < b_v$ for $v = 1, 4$. Although the block admits one more non-null variable, the faces that contain $x^*$ differ in the variable that they allow to move, and only $x_1$ and $x_4$ are free in all of them. Thus $\F(x^*) = \{1, 4\}$, and PCC-LICQ holds if and only if $\nabla h_1(x^*)$, $\nabla h_2(x^*)$, $e^2$, $e^3$, and $e^5$ are linearly independent. This choice of $\F(x^*)$ is what the proof of Proposition~\ref{propbstat} requires, and it agrees with the convention of MPEC-LICQ (see below). In a block $\nu \notin \Nc$, such as the block of slacks, a variable is free if it is strictly between its bounds. If $\F(x^*) = \emptyset$, then PCC-LICQ holds only if $m = 0$.

For problem~(\ref{pccnl}), written with slack variables, the slack of an inactive inequality
constraint is a free variable, and the slack of an active one is at its lower bound. The column
of the slack $s_i$ in the Jacobian of the constraints is the canonical vector of the row
of~$c_i$, and it belongs to $J_{\F(x^*)}$ exactly when $c_i$ is inactive. Therefore, PCC-LICQ is
equivalent to the linear independence of the vectors of partial derivatives, with respect to
the free components of~$x$, of $h_1, \ldots, h_m$ and of the active inequality constraints
$c_i$ of~(\ref{pccnl}). Equivalently, the gradients of $h_1, \ldots, h_m$ and of the active
$c_i$, together with the canonical vectors $e^v$ of the components $v$ of $x$ that are not free,
are linearly independent. Each variable that is not free contributes a single canonical vector,
even if it is null and at a bound at the same time. In the MPCC case, a null component of a
complementarity pair is at its lower bound, so that $\F(x^*)$ is the set of variables that are
strictly between their bounds, and PCC-LICQ says that the gradients of $h_1, \ldots, h_m$ and
the canonical vectors of the variables at a bound are linearly independent. This is
the usual MPEC-LICQ, that is, the LICQ of the nonlinear program in which the null components
of the pairs are fixed at zero (the tightened nonlinear program of \cite{ss00}), when the
bounds are considered as constraints. In this case we keep the name MPEC-LICQ.

Proposition~\ref{propbstat} below relates the two stationarity concepts used in this paper.
A point $x^* \in X$ is B-stationary for~(\ref{mpcc}) if $\nabla f(x^*)^T d \geq 0$ for all
$d \in T_X(x^*)$. For $\lambda \in \R^m$, let $\ell_\lambda(x) = f(x) + \lambda^T h(x)$. By Proposition~\ref{propstatb}, $x^*$ is stationary for the minimization of $\ell_\lambda$ on
$\Om$ (Definition~\ref{defstat}) if and only if $\nabla \ell_\lambda(x^*)^T d \geq 0$ for all
$d \in T_\Om(x^*)$. The second concept is the one that the methods of this paper can reach,
because the constraints $h(x) = 0$ are penalized and only $\Om$ is kept in the subproblems.

\begin{prop} \label{propbstat}
Let Assumption~\ref{assum} hold and $x^* \in X$.
\begin{enumerate}
\item If there exists $\lambda \in \R^m$ such that $x^*$ is stationary for the minimization
  of $\ell_\lambda$ on $\Om$, then $x^*$ is B-stationary for~(\ref{mpcc}).
\item If PCC-LICQ holds at $x^*$ and $x^*$ is B-stationary for~(\ref{mpcc}), then there exists
  a unique $\lambda^* \in \R^m$ such that $x^*$ is stationary for the minimization of
  $\ell_{\lambda^*}$ on $\Om$. Moreover, $\lambda^*$ is the unique solution of
  $\nabla_{\F} f(x^*) + J_{\F}(x^*)^T \lambda = 0$, where $\F = \F(x^*)$.
\end{enumerate}
\end{prop}

\begin{proof}
\emph{Item 1.} Let $d \in T_X(x^*)$, and let $t_j \downarrow 0$ and $d^j \to d$ be such that
$x^* + t_j d^j \in X$. Since $X \subseteq \Om$, $d \in T_\Om(x^*)$ and, since
$h(x^* + t_j d^j) = h(x^*) = 0$, $J(x^*) d = 0$. Therefore,
$\nabla f(x^*)^T d = \nabla \ell_\lambda(x^*)^T d \geq 0$.

\emph{Item 2.} Let $\mathcal{I}(x^*)$ be the set of face indices $I$ such that $x^* \in F_I$
and, for $I \in \mathcal{I}(x^*)$, let $X_I = F_I \cap \{x \mid h(x) = 0\} \subseteq X$.
Since $T_{X_I}(x^*) \subseteq T_X(x^*)$, $x^*$ is B-stationary for the minimization of $f$
on $X_I$. The constraints of this problem are $h(x) = 0$ and the bounds of the box $F_I$, and the variables that are strictly inside the bounds of $F_I$ at $x^*$ include those of $\F$, because a variable $x_v$ with $v \in \F$ that belongs to a block of $\Nc$ is non-null at $x^*$ and, therefore, is not fixed at zero in $F_I$. Thus, PCC-LICQ implies the LICQ for this problem at $x^*$ and, therefore,
B-stationarity implies the KKT conditions, that is, there exists $\lambda^I$ such that
$g^{\ell_{\lambda^I}}_I(x^*) = 0$. Since, for every $v \in \F$, the variable $x_v$ is not fixed at zero in any face $F_I$
with $I \in \mathcal{I}(x^*)$ and $a_v < x^*_v < b_v$, we have that
$\partial \ell_{\lambda^I}(x^*) / \partial x_v = 0$. Hence,
$\nabla_\F f(x^*) + J_\F(x^*)^T \lambda^I = 0$ and, since $J_\F(x^*)$ has rank $m$,
$\lambda^I = \lambda^*$ for all $I \in \mathcal{I}(x^*)$. Thus,
$g^{\ell_{\lambda^*}}_I(x^*) = 0$ for all $I \in \mathcal{I}(x^*)$ which, as in the proof of
Proposition~2.3 of \cite{bmtel}, is equivalent to the stationarity of $x^*$ for the
minimization of $\ell_{\lambda^*}$ on $\Om$. Uniqueness follows from the same argument.
\end{proof}

In the MPCC case, stationarity of $x^*$ for the minimization of $\ell_{\lambda^*}$ on $\Om$
means the following. For the variables that are not null components of complementarity pairs,
the partial derivatives of $\ell_{\lambda^*}$ at $x^*$ satisfy the usual sign conditions at
the bounds. There is no condition on the derivative with respect to a null component whose
partner is non-null, and both derivatives are non-negative for a biactive pair. This is
S-stationarity \cite{ss00}.

\section{Augmented Lagrangian method with telescopic subproblems} \label{alsec}

In this section we consider problem~(\ref{pccnl}), in which, besides $x \in \Om$, there are
general equality and inequality constraints $h(x) = 0$ and $c(x) \leq 0$. A natural way
of handling this problem is the augmented Lagrangian method of Algencan \cite{abms,bmbook},
in which the constraints $h(x) = 0$ and $c(x) \leq 0$ are penalized, while $x \in \Om$ (the bounds and
the partitioned cardinality constraints) is kept in the subproblems, which are solved by the
Telescopic Method. We call it the \emph{straightforward} method. We prove that, if the
penalty parameters remain bounded, every limit point of the sequence generated by the
straightforward method is B-stationary, without constraint qualifications
(Theorem~\ref{teoal}). However, if the penalty parameters tend to infinity, the limit points
may fail to be B-stationary, even if every subproblem is solved exactly
(counterexample of Proposition~\ref{propcex}, illustrated numerically at the end of Section~\ref{alseccex}). In this case, we show in Section~\ref{alsecunb} that the infeasible limit points are stationary for the infeasibility measure on a face and that, in the MPCC case, the feasible limit points at which MPEC-LICQ holds are M-stationary. This is why, in Section~\ref{algal}, we introduce a less
obvious method, in which the swaps of the Telescopic Method are replaced by
constraint-preserving swaps (Section~\ref{cps}).

\begingroup
\subsection{The Telescopic Method} \label{telsec}

In this subsection we summarize the Telescopic Method of \cite{bmtel} and the results of
\cite{bmtel} that are used in this paper. The method minimizes a function $\phi$ on $\Om$,
where the bounds satisfy the first part of Assumption~\ref{assum}, $\B$ is bounded, and $\phi$
is continuously differentiable on an open set containing $\B$, with $\nabla \phi$ Lipschitz
continuous on $\B$ with constant $L_\phi$. (In \cite{bmtel}, these are Assumptions~A1--A3,
the last one being the boundedness of a level set, which holds because $\B$ is bounded.) The
method alternates swaps, which change the face, with iterations of a bound-constrained
algorithm on the current face.

The swaps use a one-dimensional \emph{escaping procedure} (Algorithm~3.1 of \cite{bmtel}).
Let $\bar y \in \B$ with $\bar y_{j,\nu} = 0$ and $g^\phi_{j,\nu}(\bar y) \neq 0$, and let
$\psi(t) = \phi(\bar y + t e^{j,\nu})$ for $t \in [a_{j,\nu}, b_{j,\nu}]$. Given
$\alpha_{esc} > 0$ and $\sigma_{\min} > 0$, the procedure computes, for
$\sigma = \sigma_{\min}, 2\sigma_{\min}, 4\sigma_{\min}, \ldots$, the minimizer
$t = P_{[a_{j,\nu}, b_{j,\nu}]}(-\psi'(0)/\sigma)$ of the model
$\psi(0) + \psi'(0) t + \sigma t^2/2$ on $[a_{j,\nu}, b_{j,\nu}]$, and stops with
$t_{out} = t$ at the first $\sigma$ for which $\psi(t) \leq \psi(0) - \alpha_{esc} t^2$. By
Lemma~3.1 of \cite{bmtel}, it stops after finitely many trials, and
\begin{equation} \label{escdec}
\phi(\bar y + t_{out} e^{j,\nu}) \leq \phi(\bar y) - \theta \, g^\phi_{j,\nu}(\bar y)^2,
\quad \mbox{where} \quad
\theta = \frac{\alpha_{esc}}{\max\{1, \bar\sigma\}^2}
\quad \mbox{and} \quad
\bar\sigma = \max\{\sigma_{\min}, L_\phi + 2\alpha_{esc}\}.
\end{equation}
Thus, the decrease is bounded below in terms of the projected derivative of the entering
variable, with a constant $\theta$ that depends only on the parameters and on $L_\phi$.

On a face $F_I$, the minimization of $\phi$ is a bound-constrained problem, and $x \in F_I$
is a stationary point of this problem if and only if $g^\phi_I(x) = 0$. An iteration of the
method used on the faces, called Algorithm MIA in \cite{bmtel} (Gencan in our
implementation), computes a point $\hat w \in F_I$ from a point $w \in F_I$. Since the
Telescopic Method leaves a face after one iteration of Algorithm MIA and may return to it
later at a different point, the assumptions on Algorithm MIA concern the following sequences.
A sequence $\{w^j\} \subseteq F_I$ is an \emph{MIA sequence with replacements} if, for all
$j$, $\phi(w^{j+1}) \leq \phi(\hat w^{j+1})$, where $\hat w^{j+1}$ is computed by an
iteration of Algorithm MIA from $w^j$.

\begin{assump} \label{amia}
(Assumptions A4--A6 of \cite{bmtel}.) For every face $F_I$, (i) an iteration of
Algorithm MIA from $w \in F_I$ computes $\hat w \in F_I$ with $\phi(\hat w) \leq \phi(w)$,
(ii) every limit point $w^*$ of an MIA sequence with replacements satisfies
$g^\phi_I(w^*) = 0$, and (iii) if $\{w^j\}$ is an MIA sequence with replacements and
$\lim_{j \in J} w^j = w^*$ for an infinite set $J$, then
$\lim_{j \in J} \|\hat w^{j+1} - w^j\| = 0$.
\end{assump}

These are natural properties of bound-constrained methods, and (iii) holds, for example, if
$\|\hat w - w\| \leq C \|g^\phi_I(w)\|$ for some constant $C$.

The Telescopic Method (Algorithm~5.1 of \cite{bmtel}) is given $\delta \in (0, +\infty]$,
$\eta_{\rm br} \in (0, +\infty]$, the parameters of the escaping procedure, and $x^0 \in F_{I^0}$. At
iteration $k$, it performs a \emph{swap search} at $x^k \in F_{I^k}$. For every
$\nu \in \Nc$, every $i \in I^k_\nu$ with $|x^k_{i,\nu}| \leq \delta$ and every
$j \notin I^k_\nu$ such that $g^\phi_{j,\nu}(\bar y) \neq 0$, where
$\bar y = x^k - x^k_{i,\nu} e^{i,\nu}$ is $x^k$ with the leaving variable set to zero, the
escaping procedure gives the point $w = \bar y + t_{out} e^{j,\nu}$, which belongs to the
face obtained from $I^k$ by replacing $i$ with $j$ in $I^k_\nu$, and satisfies
$\phi(w) < \phi(\bar y)$. The best of these points and its value $f_{best}$ are kept, and
the search is interrupted as soon as $f_{best} \leq \phi(x^k) - \eta_{\rm br}$ (brute acceptance). If
$f_{best} < \phi(x^k)$, the best point becomes $z^k$ and its face index becomes $J^k$.
Otherwise, $z^k = x^k$ and $J^k = I^k$. Then $x^{k+1} \in F_{J^k}$ is computed by one
iteration of Algorithm MIA from $z^k$, and $I^{k+1} = J^k$. The analysis uses only that
$z^k \in F_{J^k}$ and $\phi(z^k) \leq f_{best}$. Therefore, as observed in Remark~5.1 of
\cite{bmtel}, $z^k$ may be the best of several candidates improved by Algorithm MIA in their
new faces (the \emph{look-ahead} strategy of our implementation), and any finite number of
iterations of Algorithm MIA may be performed on each visited face. The following lemma
collects the results of \cite{bmtel} used in this paper.

\begin{lem} \label{lemtel}
Let the assumptions above on $\B$ and $\phi$ and Assumption~\ref{amia} hold, and let
$\{x^k\}$ and $\{z^k\}$ be generated by the Telescopic Method applied to the minimization of
$\phi$ on $\Om$.
\begin{enumerate}
\item For all $k$, $\phi(x^{k+1}) \leq \phi(z^k) \leq \phi(x^k)$, and
  $\lim_{k \to \infty} [\phi(x^k) - \phi(z^k)] = 0$.
\item If $J^k = I$ for all $k$ in an infinite set $K$, then every limit point $z^*$ of
  $\{z^k\}_{k \in K}$ belongs to $F_I$ and satisfies $g^\phi_I(z^*) = 0$ and, if
  $\lim_{k \in K} z^k = z^*$, then $\lim_{k \in K} x^{k+1} = z^*$.
\item Every limit point of $\{x^k\}$ is stationary for the minimization of $\phi$ on $\Om$
  (Definition~\ref{defstat}) or, equivalently, B-stationary.
\item For every $\varepsilon > 0$, there exists $k$ such that the swap search at iteration $k$
  is not interrupted by the brute acceptance test, $\|g^\phi_{I^k}(x^k)\|_\infty \leq
  \varepsilon$, and $\phi(w) \geq \phi(x^k) - \varepsilon$ for every point $w$ computed in
  that swap search.
\end{enumerate}
\end{lem}

\begin{proof}
Items 1 and 2 are Lemmas~6.1, 6.2 and~6.3 of \cite{bmtel}, item~3 is Theorem~6.4 of
\cite{bmtel}, and item~4 is Lemma~6.5 of \cite{bmtel}. The proof of item~3 combines items~1
and~2 with~(\ref{escdec}). If $J^k = I$ and $z^k \to x^*$ for $k \in K$, and a block
$\nu \in \Nc$ with $\|\vx^*_\nu\|_0 < p_\nu$ had a variable $x_{j,\nu}$, $j \notin I_\nu$,
with $g^\phi_{j,\nu}(x^*) \neq 0$, then, since $I_\nu$ contains an index $i$ with
$x^*_{i,\nu} = 0$, by~(\ref{escdec}), the swap $x_{i,\nu} \to x_{j,\nu}$ at $x^{k+1}$ (or the brute acceptance
test, if it interrupts the search first) would give $\phi(z^{k+1}) \leq \phi(x^{k+1}) - c$ for
a constant $c > 0$ and all $k \in K$ large enough, contradicting item~1.
\end{proof}
\endgroup

\subsection{Algorithm and convergence with bounded penalty parameters} \label{alsecconv}

We denote by $X = \{x \in \Om \mid h(x) = 0, \; c(x) \leq 0\}$ the feasible set of~(\ref{pccnl}), and say that $x^* \in X$ is B-stationary for~(\ref{pccnl}) if
$\nabla f(x^*)^T d \geq 0$ for all $d \in T_X(x^*)$. We denote by $\nabla h(x)$ and
$\nabla c(x)$ the $n \times m$ and $n \times q$ matrices whose columns are the gradients of
the components of $h$ and $c$, respectively. For $\rho > 0$, $\lambda \in \R^m$ and $\mu \in \R^q_+$, the
(Powell--Hestenes--Rockafellar) augmented Lagrangian is
\[
L_\rho(x, \lambda, \mu) = f(x) + \frac{\rho}{2} \left[ \left\| h(x) + \frac{\lambda}{\rho} \right\|^2
+ \left\| \left( c(x) + \frac{\mu}{\rho} \right)_+ \right\|^2 \right].
\]

In this section we use the following assumption, under which the augmented Lagrangian has the
regularity required by the Telescopic Method.

\begin{assump} \label{aal}
The bounds satisfy the first part of Assumption~\ref{assum}, $\B$ is bounded, and $f$, $h$ and $c$ are continuously
differentiable on an open set containing $\B$, with Lipschitz continuous gradients on $\B$.
\end{assump}

Under Assumption~\ref{aal}, for all $\rho$, $\lambda$ and $\mu$, the function
$L_\rho(\cdot, \lambda, \mu)$ satisfies the assumptions on $\phi$ of Section~\ref{telsec}, since
$t \mapsto (t)_+^2$ is continuously differentiable with Lipschitz continuous derivative.
For such functions $\phi$, the projected derivatives $g^\phi_{i,\nu}$ and $g^\phi_I$ are
defined in Section~\ref{probl}. We denote by $\bar\varepsilon$ half of the smallest absolute value of the non-null bounds $a_v$, $b_v$ of the variables of the blocks of $\Nc$.

The following algorithm uses the safeguarded multiplier updates and the penalty parameter
update of Algencan \cite{abms,bmbook}.

\medskip
\noindent
{\bf Algorithm \ref{alsec}.1 (Augmented Lagrangian with telescopic subproblems)}

\noindent
Let $\tau \in (0,1)$, $\gamma > 1$, $\rho_1 > 0$, $\lambda_{\min} < \lambda_{\max}$,
$\mu_{\max} > 0$, $\bar\lambda^1 \in [\lambda_{\min}, \lambda_{\max}]^m$,
$\bar\mu^1 \in [0, \mu_{\max}]^q$, a sequence $\{\varepsilon_k\}$ with $\varepsilon_k > 0$ and
$\lim_{k\to\infty} \varepsilon_k = 0$, the parameters $\delta$ and $\eta_{\rm br}$ of
the Telescopic Method, the parameters of the escaping procedure, and $x^0 \in \Om$ be
given. Set $k \leftarrow 1$.

\smallskip
\noindent
{\bf Step 1 (Subproblem).}
Let $L_k(x) = L_{\rho_k}(x, \bar\lambda^k, \bar\mu^k)$. Starting from $x^{k-1}$, apply
the Telescopic Method to the minimization of $L_k$ on $\Om$ until a point $x^k \in \Om$ and a
face index $I^k$ with $x^k \in F_{I^k}$ are found such that
\begin{itemize}
\item[(a)] $\|g^{L_k}_{I^k}(x^k)\|_\infty \leq \varepsilon_k$;
\item[(b)] no original swap decreases $L_k$ significantly from $x^k$: for every
  $\nu \in \Nc$, $i \in I^k_\nu$ with $|x^k_{i,\nu}| \leq \delta$ and $j \notin I^k_\nu$, let
  $\bar y = x^k - x^k_{i,\nu} e^{i,\nu}$ ($x^k$ with the leaving variable set to zero) and, if
  $g^{L_k}_{j,\nu}(\bar y) \neq 0$, the point $w = \bar y + t_{out} e^{j,\nu}$, where $t_{out}$
  is computed by the one-dimensional escaping procedure of the Telescopic Method
  (Section~\ref{telsec}) applied to $t \mapsto L_k(\bar y + t e^{j,\nu})$,
  $t \in [a_{j,\nu}, b_{j,\nu}]$, satisfies $L_k(w) \geq L_k(x^k) - \varepsilon_k$;
\item[(c)] for every $\nu \in \Nc$, $i \in I^k_\nu$ with
  $|x^k_{i,\nu}| \leq \min\{\varepsilon_k, \bar\varepsilon\}$ and $j \notin I^k_\nu$,
  $|\partial_{i,\nu} L_k(x^k)| \leq \varepsilon_k$ or $|g^{L_k}_{j,\nu}(x^k)| \leq \varepsilon_k$.
\end{itemize}

\smallskip
\noindent
{\bf Step 2 (Penalty parameter).}
Let $V^k = \max\{c(x^k), -\bar\mu^k/\rho_k\}$ (componentwise) and
$M_k = \max\{\|h(x^k)\|_\infty, \|V^k\|_\infty\}$. If $k > 1$ and $M_k > \tau M_{k-1}$, set $\rho_{k+1} = \gamma \rho_k$. Otherwise, set $\rho_{k+1} = \rho_k$.

\smallskip
\noindent
{\bf Step 3 (Multipliers).}
Set
\[
\begin{array}{ll}
\lambda^{k+1} = \bar\lambda^k + \rho_k h(x^k), &
\mu^{k+1} = (\bar\mu^k + \rho_k c(x^k))_+,
\end{array}
\]
choose $\bar\lambda^{k+1} \in [\lambda_{\min}, \lambda_{\max}]^m$ and
$\bar\mu^{k+1} \in [0, \mu_{\max}]^q$ (for example, the projections of $\lambda^{k+1}$ and $\mu^{k+1}$ onto these boxes, which are used in our implementation), set $k \leftarrow k + 1$, and go to Step~1.

\bigskip

Conditions (a) and (b) say that $x^k$ is approximately stationary in its face (the stopping criterion of the bound-constrained solver) and that no swap
examined by the Telescopic Method decreases $L_k$ significantly, which are the properties used in the proof of item~3 of Lemma~\ref{lemtel}. The practical stopping criterion of
Section~5 of \cite{bmtel} (no swap found and~(14) of \cite{bmtel} with tolerance $\varepsilon_k$)
implies (a) and (b). Condition (c) says that, if a variable of the current face is almost null, then either it is approximately stationary also as an interior variable, or the null variables of its block cannot decrease $L_k$ by moving away from zero. Conditions (a) and (c) together are an approximate version, in terms of derivatives, of the stationarity of $x^k$ for the minimization of $L_k$ on $\Om$ (Definition~\ref{defstat}), and condition (c) is the part of~(\ref{s2}) that looks outside the current face. Unlike the decrease guaranteed by~(\ref{escdec}), which is used in (b), its tolerance does not deteriorate when $\rho_k$ grows. Condition (c) is not used in Theorem~\ref{teoal}, and it is what gives M-stationarity when $\rho_k \to \infty$ (Section~\ref{alsecunb}). The following lemma shows that Step~1 is well defined, and
Theorem~\ref{teoal} shows that, if the penalty parameters remain bounded, the limit points
are B-stationary.

\begin{lem} \label{lemalwd}
Let Assumption~\ref{aal} hold and let $\phi = L_\rho(\cdot, \lambda, \mu)$ for some
$\rho > 0$, $\lambda$ and $\mu \geq 0$. Let the Telescopic Method be applied to the
minimization of $\phi$ on $\Om$ and let $\varepsilon > 0$. Then, after a finite number of
iterations, the point $x^j$ at which the swap search of the Telescopic Method is performed
satisfies (a), (b) and (c) with $\phi$, $x^j$, $I^j$ and $\varepsilon$ in the places of $L_k$,
$x^k$, $I^k$ and $\varepsilon_k$, respectively.
\end{lem}

\begin{proof}
Since $\phi$ satisfies the assumptions on $\phi$ of Section~\ref{telsec}, the proofs of items~3 and~4 of Lemma~\ref{lemtel} (Theorem~6.4 and Lemma~6.5 of \cite{bmtel}) apply to $\phi$. By the proof of item~4, (b) holds for all $j$ large enough. By the proof of item~3, there exist an infinite set $K$, a face index $I$ and a point $x^*$, stationary for the minimization of $\phi$ on $\Om$, such that $\lim_{j \in K} x^j = x^*$ and $I^j = I$ for all $j \in K$. By the continuity of $g^\phi_I$, (a) holds for $j \in K$ large enough. Let $\bar\varepsilon_1 = \min\{\varepsilon, \bar\varepsilon\}$, $\nu \in \Nc$, $i \in I_\nu$ and $\ell \notin I_\nu$. If $x^*_{i,\nu} = 0$, then $\|\vx^*_\nu\|_0 < p_\nu$, because $\# I_\nu = p_\nu$ and $x^*_{\ell',\nu} = 0$ for all $\ell' \notin I_\nu$. Then, by~(\ref{s2}), $g^\phi_{\ell,\nu}(x^*) = 0$ and, by continuity, $|g^\phi_{\ell,\nu}(x^j)| \leq \varepsilon$ for $j \in K$ large enough. If $0 < |x^*_{i,\nu}| \leq \bar\varepsilon_1$, then $x^*_{i,\nu}$ is not a bound, because the non-null bounds have absolute values at least $2\bar\varepsilon$. Thus, by~(\ref{s1}), $\partial_{i,\nu} \phi(x^*) = 0$ and, by continuity, $|\partial_{i,\nu} \phi(x^j)| \leq \varepsilon$ for $j \in K$ large enough. Finally, if $|x^*_{i,\nu}| > \bar\varepsilon_1$, then $|x^j_{i,\nu}| > \bar\varepsilon_1$ for $j \in K$ large enough, and (c) imposes no condition on $i$. Therefore, (a), (b) and (c) hold for $j \in K$ large enough.
\end{proof}

\begin{teo} \label{teoal}
Let Assumption~\ref{aal} hold, and let $\{x^k\}$ be generated by Algorithm~\ref{alsec}.1.
If $\{\rho_k\}$ is bounded, then every limit point $x^*$ of $\{x^k\}$ is feasible
for~(\ref{pccnl}) and B-stationary. Moreover, there exist $\hat\lambda \in \R^m$ and
$\hat\mu \in \R^q_+$, with $\hat\mu_i = 0$ whenever $c_i(x^*) < 0$, such that $x^*$ is
stationary for the minimization of $f + \hat\lambda^T h + \hat\mu^T c$ on $\Om$.
\end{teo}

\begin{proof}
Since $\{\rho_k\}$ is bounded and non-decreasing, there exist $\bar\rho$ and $k_1$ such that
$\rho_k = \bar\rho$ for all $k \geq k_1$. By Step~2, $M_k \leq \tau^{k - k_1} M_{k_1}$ for
$k \geq k_1$, so that $h(x^k) \to 0$ and $V^k \to 0$. Since $c(x^k) \leq V^k$, every limit
point of $\{x^k\}$ is feasible.

Let $x^*$ be a limit point. Since the safeguarded multipliers are bounded and there are
finitely many face indices, there exist an infinite set $K$, $\tilde\lambda$, $\tilde\mu$
and a face index $I$ such that $\lim_{k \in K} x^k = x^*$,
$\lim_{k \in K} \bar\lambda^k = \tilde\lambda$, $\lim_{k \in K} \bar\mu^k = \tilde\mu$, and
$I^k = I$ for all $k \in K$. If $c_i(x^*) < 0$, then $c_i(x^k) < c_i(x^*)/2$ for $k \in K$
large enough and, since $V^k_i \to 0$, $\bar\mu^k_i / \bar\rho \to 0$, so that $\tilde\mu_i = 0$. Define $\Phi(x) = L_{\bar\rho}(x, \tilde\lambda, \tilde\mu)$.

Since $(x, \lambda, \mu) \mapsto L_{\bar\rho}(x, \lambda, \mu)$ and its gradient with respect
to $x$ are continuous, and $\B$ and the safeguarding boxes are compact, $L_k \to \Phi$,
$\nabla L_k \to \nabla \Phi$ and, consequently, $g^{L_k}_{i,\nu} \to g^\Phi_{i,\nu}$ uniformly
on $\B$ for $k \in K$. Moreover, there exist $G$ and $\hat L$ such that, for all
$k \geq k_1$, $\|\nabla L_k(x)\| \leq G$ for all $x \in \B$ and $\nabla L_k$ is Lipschitz
continuous on $\B$ with constant $\hat L$. Hence (\ref{escdec}) holds for all the
functions $L_k$, $k \geq k_1$, with the same constant
$\hat\theta = \alpha_{esc} / \max\{1, \hat\sigma\}^2$, where
$\hat\sigma = \max\{\sigma_{\min}, \hat L + 2\alpha_{esc}\}$.

We prove that $x^*$ is stationary for the minimization of $\Phi$ on $\Om$ in the sense of
Definition~\ref{defstat}. By (a) and the uniform convergence, $g^\Phi_{i,\nu}(x^*) = 0$ for
all $i \in I_\nu$, $\nu = 1, \ldots, N$, which, since $x^* \in F_I$, implies condition~(\ref{s1}). Let
$\nu \in \Nc$ be such that $\|\vx^*_\nu\|_0 < p_\nu$ and let $\ell \notin I_\nu$ (for
$\nu \notin \Nc$, condition~(\ref{s2}) follows from the previous observation). As in the proof of item~3 of Lemma~\ref{lemtel}, there exists $i \in I_\nu$ with $x^*_{i,\nu} = 0$. Assume, by contradiction, that
$g^\Phi_{\ell,\nu}(x^*) = \kappa \neq 0$, and let
$\bar y^k = x^k - x^k_{i,\nu} e^{i,\nu}$, so that $\lim_{k \in K} \bar y^k = x^*$. For
$k \in K$ large enough, $|x^k_{i,\nu}| \leq \delta$,
$|g^{L_k}_{\ell,\nu}(\bar y^k)| \geq |\kappa|/2$,
$L_k(\bar y^k) \leq L_k(x^k) + G |x^k_{i,\nu}| \leq L_k(x^k) + \hat\theta \kappa^2/8$, and
$\varepsilon_k < \hat\theta \kappa^2/8$. Then the pair $(i, \ell)$ is examined in the swap
search from $x^k$ and, by (\ref{escdec}), the corresponding point $w$ satisfies
\[
L_k(w) \leq L_k(\bar y^k) - \hat\theta \kappa^2/4 \leq L_k(x^k) - \hat\theta \kappa^2/8
< L_k(x^k) - \varepsilon_k,
\]
which contradicts (b). Therefore, condition~(\ref{s2}) holds and, by Proposition~\ref{propstatb} (applied
to $\Phi$), $\nabla \Phi(x^*)^T d \geq 0$ for all $d \in T_\Om(x^*)$.

Since $h(x^*) = 0$,
$\nabla \Phi(x^*) = \nabla f(x^*) + \nabla h(x^*) \hat\lambda + \nabla c(x^*) \hat\mu$, where
$\hat\lambda = \tilde\lambda$ and $\hat\mu = (\tilde\mu + \bar\rho c(x^*))_+ \geq 0$, with $\tilde\mu_i = \hat\mu_i = 0$ if $c_i(x^*) < 0$. This proves the second statement
of the theorem. Let $d \in T_X(x^*)$, and let $t_j \downarrow 0$ and $d^j \to d$ be such that
$x^* + t_j d^j \in X$ for all $j$. Since $X \subseteq \Om$, $d \in T_\Om(x^*)$. Since
$h(x^* + t_j d^j) = h(x^*) = 0$, $\nabla h(x^*)^T d = 0$ and, if $c_i(x^*) = 0$, since
$c_i(x^* + t_j d^j) \leq 0$, $\nabla c_i(x^*)^T d \leq 0$. Therefore,
\[
\nabla f(x^*)^T d = \nabla \Phi(x^*)^T d - \hat\lambda^T \nabla h(x^*)^T d
- \sum_{i \,:\, c_i(x^*) = 0} \hat\mu_i \nabla c_i(x^*)^T d \geq 0,
\]
so that $x^*$ is B-stationary for~(\ref{pccnl}).
\end{proof}

No constraint qualification is used in Theorem~\ref{teoal}. The only property of the
constraints used in the last part of the proof is that the tangent cone of $X$ is contained in
the intersection of the tangent cone of $\Om$ with the linearized cone of the constraints
$h(x) = 0$ and $c(x) \leq 0$, which is always true. The hypothesis of the theorem is the
boundedness of the penalty parameters. For nonlinear programming problems with general
lower-level constraints, the boundedness of the penalty parameters of Algencan was proved in
\cite{abms} under the linear independence of the gradients of the active constraints, a
second-order condition, and the assumption that the Lagrange multipliers lie in the interior
of the safeguarding boxes. Whether similar conditions guarantee the boundedness of
$\{\rho_k\}$ in Algorithm~\ref{alsec}.1, in which the lower-level set $\Om$ is not convex, is
an open question.

\subsection{Unbounded penalty parameters: a counterexample} \label{alseccex}

We first describe the mechanism of the failure. B-stationarity with respect to $\Om$
is not a closed condition. At an iterate $x^k$ at which a complementarity pair has one
non-null component, the stationarity condition says nothing about the derivative with respect
to the null component, but at the limit point both components may vanish, and B-stationarity then
requires the corresponding partial derivative (of the Lagrangian) to be non-negative. In the
Telescopic Method, this is taken care of by the swaps, in which the non-null component is set to zero
and the null one is moved away from zero by a one-dimensional procedure whose decrease is
bounded below by $\theta g^2$, where $g$ is the projected derivative of the entering variable.
For equality constraints only, the augmented Lagrangian is
\[
L_k(x) = f(x) + \frac{\rho_k}{2} \left\| h(x) + \frac{\bar\lambda^k}{\rho_k} \right\|^2,
\]
and both parts of a swap are penalized. Setting a variable $x_i$ to zero changes $h$ by
$O(|x_i|)$, which costs $O(\rho_k x_i^2)$ in the penalty term. Moreover, the constant $\theta$
deteriorates as $\rho_k$ grows, because the curvature of $L_k$ along the entering variable
contains the term $\rho_k (\partial h / \partial x_j)^2$. In the counterexample below, the cost
of the first part of the swap tends to infinity.

When $\rho_k \to \infty$, the argument of the proof of Theorem~\ref{teoal} breaks down. The
Lipschitz constants of $\nabla L_k$ grow with $\rho_k$, so that the constant $\theta$ of~(\ref{escdec}) tends to zero, and the difference $L_k(\bar y^k) - L_k(x^k)$, which
involves the penalty of the (possibly infeasible) point $\bar y^k$, need not tend to zero. The
following example shows that, in this case, the limit points of Algorithm~\ref{alsec}.1 may
fail to be B-stationary, even if every subproblem is solved exactly.

\paragraph{Example 1.}
Let $n = 3$, $N = 2$, $\vx_1 = (x_1, x_2)$ with $p_1 = 1$, $\vx_2 = x_3$ with $p_2 = 1$ (so that
$\Nc = \{1\}$), $\B = [0,2] \times [0,1] \times [0,2]$, $c \equiv 0$ (no inequality
constraints), and
\begin{equation} \label{cex}
f(x) = -\frac{x_2}{2} - \frac{x_1^2}{2} + \frac{x_3^2}{2} + \frac{x_3^4}{4},
\qquad
h(x) = \frac{(x_1 - x_3)^2}{2} + x_2 x_3^2.
\end{equation}
Since $h(x) = 0$ with $x \geq 0$ implies $x_1 = x_3$ and $x_2 x_3 = 0$, and $x_1 x_2 = 0$ in
$\Om$, the feasible set is $X = \{(t, 0, t) \mid t \in [0, 2]\} \cup \{(0, t, 0) \mid t \in
[0, 1]\}$, and the global minimizer is $(0, 1, 0)$, with $f = -1/2$. In the
iterations of Algorithm~\ref{alsec}.1, $L_k(x) = f(x) + \bar\lambda^k h(x) + \rho_k h(x)^2/2$
up to an additive constant.

\begin{prop} \label{propcex}
Consider problem~(\ref{pccnl}) with the data~(\ref{cex}) and Algorithm~\ref{alsec}.1 with
$\lambda_{\min} \leq 0 < \lambda_{\max} < \infty$, $\bar\lambda^1 \in [0, \lambda_{\max}]$,
$\rho_1 > 4$, $\bar\lambda^{k+1} = P_{[\lambda_{\min}, \lambda_{\max}]}(\lambda^{k+1})$, and arbitrary $\tau$, $\gamma$, $\delta$, $\eta_{\rm br}$ and $\{\varepsilon_k\}$. For
each $k$, let $s_k$ be the unique solution in $(0,1)$ of
\begin{equation} \label{cexs}
\bar\lambda^k + \frac{\rho_k s^6}{2} = \frac{1 + s^2}{s^2},
\end{equation}
and let $x^k = (s_k + s_k^3, 0, s_k)$ and $I^k = (\{1\}, \{1\})$. Then:
\begin{enumerate}
\item $x^k$ is a strict local minimizer of $L_k$ on $\Om$ and satisfies (a), (b) and (c) with
  $\varepsilon_k = 0$. Moreover, condition (b) holds trivially, because no swap is a
  candidate ($g^{L_k}_{j,\nu}(\bar y) = 0$ for the only possible swap).
\item $\lim_{k \to \infty} \rho_k = \infty$ and $\lim_{k \to \infty} x^k = x^* = 0$, which is
  feasible.
\item $x^*$ is not B-stationary for~(\ref{pccnl}).
\end{enumerate}
\end{prop}

\begin{proof}
For all $k$, $\bar\lambda^k \in [0, \lambda_{\max}]$ and $\rho_k \geq \rho_1 > 4$. Indeed,
$\{\rho_k\}$ is non-decreasing and, if $\bar\lambda^k \geq 0$, then, as shown below,
$\lambda^{k+1} = \bar\lambda^k + \rho_k h(x^k) = (1 + s_k^2)/s_k^2 > 0$, so that
$\bar\lambda^{k+1} = \min\{\lambda^{k+1}, \lambda_{\max}\} \in (0, \lambda_{\max}]$. The
left-hand side of~(\ref{cexs}) is increasing and the right-hand side is decreasing in $s > 0$, the right-hand side tends to $+\infty$ when $s \to 0$, and at $s = 1$ the left-hand side,
$\bar\lambda^k + \rho_k/2 > 2$, is larger than the right-hand side. Thus, $s_k \in (0,1)$ is
well defined and $x^k \in F_{I^k}$, with $0 < x^k_1, x^k_3 < 2$.

\emph{Item 1.} Let $s = s_k$. Then $x^k_1 - x^k_3 = s^3$, $h(x^k) = s^6/2$ and, by~(\ref{cexs}),
$\bar\mu := \bar\lambda^k + \rho_k h(x^k) = (1 + s^2)/s^2$. Since
$\nabla L_k(x) = \nabla f(x) + (\bar\lambda^k + \rho_k h(x)) \nabla h(x)$, with
$\nabla f(x) = (-x_1, -1/2, x_3 + x_3^3)^T$ and
$\nabla h(x) = (x_1 - x_3, x_3^2, -(x_1 - x_3) + 2 x_2 x_3)^T$, we have that
$\partial L_k(x^k) / \partial x_1 = -(s + s^3) + \bar\mu s^3 = 0$ and
$\partial L_k(x^k) / \partial x_3 = s + s^3 - \bar\mu s^3 = 0$. Since $x^k_1$ and $x^k_3$ are
strictly between their bounds, $g^{L_k}_{I^k}(x^k) = 0$, so (a) holds. In a neighborhood of $x^k$,
$x_1 > 0$ and, therefore, $\Om$ coincides with $F_{I^k} = \{x \in \B \mid x_2 = 0\}$. The
Hessian of $L_k$ with respect to $(x_1, x_3)$ at $x^k$ is
\[
\begin{pmatrix} P - 1 & -P \\ -P & P + 1 + 3 s^2 \end{pmatrix},
\quad \mbox{where} \quad
P = \bar\lambda^k + \frac{3}{2} \rho_k s^6 \geq \bar\mu = \frac{1 + s^2}{s^2}.
\]
Its determinant is $3 s^2 P - 1 - 3 s^2 \geq 3 s^2 \bar\mu - 1 - 3 s^2 = 2$ and
$P - 1 \geq 1/s^2 > 0$, so it is positive definite and $x^k$ is a strict local minimizer of
$L_k$ on $\Om$. In the swap search from $x^k$, the only possibility is $i = 1$, $j = 2$ in
block $\nu = 1$, with $\bar y = (0, 0, s)$ and $h(\bar y) = s^2/2$. By~(\ref{cexs}),
$\rho_k s^4/2 = (1 + s^2)/s^4 - \bar\lambda^k / s^2$ and, since
$\bar\lambda^k \leq \bar\mu = (1 + s^2)/s^2$,
\[
\frac{\partial L_k}{\partial x_2}(\bar y) = -\frac{1}{2} + \left( \bar\lambda^k + \frac{\rho_k s^2}{2}
\right) s^2 = -\frac{1}{2} + \frac{1 + s^2}{s^4} - \bar\lambda^k \left( \frac{1}{s^2} - s^2 \right)
\geq -\frac{1}{2} + 1 + s^2 > 0.
\]
Since $\bar y_2 = 0$ is the lower bound of $x_2$, $g^{L_k}_{2,1}(\bar y) = 0$. Therefore, no
point $w$ is computed in (b), and (b) holds trivially. Condition (c) also holds, because $\partial L_k(x^k)/\partial x_1 = 0$.

\emph{Item 2.} Assume that $\{\rho_k\}$ is bounded. Then there exist $\bar\rho$ and $k_1$ such
that $\rho_k = \bar\rho$ and $h(x^k) \leq \tau h(x^{k-1})$ for all $k > k_1$. Thus,
$s_k^6 / 2 = h(x^k) \to 0$ and, by~(\ref{cexs}),
$\bar\lambda^k + \bar\rho\, h(x^k) = (1 + s_k^2)/s_k^2 \to \infty$, which contradicts
$\bar\lambda^k \leq \lambda_{\max}$. Therefore, $\rho_k \to \infty$. Since $\bar\lambda^k \geq 0$
and $s_k < 1$, (\ref{cexs}) implies that $\rho_k s_k^6 / 2 \leq 2 / s_k^2$, i.e.,
$s_k \leq (4/\rho_k)^{1/8} \to 0$. Hence, $x^k \to 0$, and $h(0) = 0$.

\emph{Item 3.} For all $t \in [0,1]$, $(0, t, 0) \in X$. Thus, $e^2 = (0,1,0)^T \in T_X(0)$ and
$\nabla f(0)^T e^2 = -1/2 < 0$.
\end{proof}

At the limit point of this example, $\nabla h(0) = 0$, so that no constraint qualification
holds. In Section~\ref{nec} we show that the method may also fail under MPEC-LICQ, with the
quadratic penalty method (counterexample of Proposition~\ref{propcex2}) or with inexact
solutions of the subproblems (counterexample of Proposition~\ref{propcex3}). It may fail even with
exact solutions of the subproblems and multipliers in the interior of the safeguarding box
(counterexample of Proposition~\ref{propcex4}).

Thus, the iterates of Algorithm~\ref{alsec}.1 may be strict local minimizers of the
subproblems, which are B-stationary for the subproblems and are accepted by the
Telescopic Method with any tolerance, and yet converge to a point that is not B-stationary
for~(\ref{pccnl}). (Whether the Telescopic Method started at $x^{k-1}$ actually returns the point $x^k$ of
Proposition~\ref{propcex} depends on the algorithm that it uses for the iterations within a
face, called Algorithm MIA in \cite{bmtel}. In \cite{bmtel}, Algorithm MIA is the
bound-constrained solver Gencan \cite{gencan,bmbook}, and the numerical results below show
that, with Gencan, it does.) The same conclusion holds without safeguards ($\lambda_{\max} = +\infty$,
i.e., $\bar\lambda^{k+1} = \lambda^{k+1}$). Indeed, if $\rho_k = \bar\rho$ for $k > k_1$, then
$\bar\lambda^{k+1} = \bar\lambda^k + \bar\rho h(x^k)$, so that $\{\bar\lambda^k\}$ is bounded,
because $\{h(x^k)\}$ converges to zero linearly, and the contradiction of Item~2 follows as
before. Since no swap is a candidate in (b), the look-ahead strategy of Section~\ref{telsec} does not change the iterates either.

In the example, the multiplier estimates $\lambda^{k+1} = (1 + s_k^2)/s_k^2$ tend to infinity. In nonlinear programming, in the absence of constraint qualifications,
augmented Lagrangian methods are known to generate sequences that satisfy sequential
(approximate) optimality conditions, which hold at every local minimizer but are weaker than
the KKT conditions \cite{bmbook}. Here, each $x^k$ is exactly B-stationary for the
minimization of $f + \lambda^{k+1} h$ on $\Om$ (by Item~1, since
$\nabla L_k(x^k) = \nabla f(x^k) + \lambda^{k+1} \nabla h(x^k)$) and $h(x^k) \to 0$, but the
limit point is not B-stationary. As described at the beginning of this subsection, the swap
that would detect the descent along $x_2$ passes through the point $\bar y^k$, which is far from
feasible ($h(\bar y^k) = s_k^2/2$, whereas $h(x^k) = s_k^6/2$), and its penalty grows with
$\rho_k$, since
$L_k(\bar y^k) - L_k(x^k) \geq \rho_k [h(\bar y^k)^2 - h(x^k)^2]/2 = \rho_k s_k^4 (1 - s_k^8)/8$.

To illustrate this behavior numerically, we applied Algorithm~\ref{alsec}.1 to problem~(\ref{pccnl}) with the data~(\ref{cex}), using
the parameters of Algencan ($\tau = 0.5$, $\gamma = 10$, $\lambda_{\min} = -10^{20}$,
$\lambda_{\max} = 10^{20}$), $\bar\lambda^1 = 0$, $\rho_1 = 10$, and $x^0 = (1, 0, 0.9)$.
The subproblems were solved by our implementation of the Telescopic Method with Gencan,
exact second derivatives and the default parameters of Section~7 of \cite{bmtel} (in particular,
$\varepsilon = 10^{-8}$, $\delta = \eta_{\rm br} = +\infty$, and look-ahead), started at $x^{k-1}$.
Table~\ref{tabal} shows the results of the first 40 iterations. Every subproblem ended at an
approximately B-stationary point of the subproblem, with no swaps, and $x^k_2 = 0$ for all
$k$. The iterates coincide, to about six significant digits, with the points of
Proposition~\ref{propcex}, computed independently by solving~(\ref{cexs}) with the same
updates of $\bar\lambda^k$ and $\rho_k$. The last two columns show that the swap that would
lead to the global minimizer becomes increasingly unprofitable. Its first step (setting
$x_1$ to zero) increases $L_k$ by an amount that grows with $\rho_k$, and at $\bar y^k$ the
derivative of $L_k$ with respect to $x_2$ is positive, so that $x_2$ is not a candidate to
enter the face. With a feasibility tolerance of $10^{-8}$, an implementation of the augmented
Lagrangian method would stop at iteration 21, at $x^{21} \approx (0.049, 0, 0.049)$, with
$f(x^{21}) \approx -4.3 \times 10^{-6}$, whereas $f(0, 1, 0) = -1/2$. Computing the points of
Proposition~\ref{propcex} with the same updates in 120-digit arithmetic, we verified that the
same behavior persists up to iteration 80, with $\rho_{80} = 10^{40}$,
$x^{80} \approx (9.5 \times 10^{-6}, 0, 9.5 \times 10^{-6})$,
$L_{80}(\bar y^{80}) - L_{80}(x^{80}) \approx 1.0 \times 10^{19}$ and
$\partial L_{80}(\bar y^{80}) / \partial x_2 \approx 4.1 \times 10^{19}$, and that the
conclusions are the same with $\lambda_{\max} = 10$.

\begin{table}[ht!]
\centering
\small
\setlength{\tabcolsep}{4pt}
\begin{tabular}{rrrrrrrrr}
\hline
$k$ & $\rho_k$ & $\bar\lambda^k$ & $x^k_1$ & $x^k_3$ & $h(x^k)$ & $f(x^k)$ &
$L_k(\bar y^k) - L_k(x^k)$ & $\partial L_k(\bar y^k)/\partial x_2$ \\
\hline
1 & $10^{1}$ & $0$ & $1.556\mathrm{e}{0}$ & $8.784\mathrm{e}{\mbox{-}1}$ & $2.3\mathrm{e}{\mbox{-}1}$ & $-6.8\mathrm{e}{\mbox{-}1}$ & $1.7\mathrm{e}{0}$ & $2.5\mathrm{e}{0}$ \\
2 & $10^{1}$ & $2.30\mathrm{e}{0}$ & $1.079\mathrm{e}{0}$ & $7.143\mathrm{e}{\mbox{-}1}$ & $6.6\mathrm{e}{\mbox{-}2}$ & $-2.6\mathrm{e}{\mbox{-}1}$ & $1.3\mathrm{e}{0}$ & $2.0\mathrm{e}{0}$ \\
4 & $10^{2}$ & $3.35\mathrm{e}{0}$ & $6.863\mathrm{e}{\mbox{-}1}$ & $5.340\mathrm{e}{\mbox{-}1}$ & $1.2\mathrm{e}{\mbox{-}2}$ & $-7.3\mathrm{e}{\mbox{-}2}$ & $1.7\mathrm{e}{0}$ & $4.5\mathrm{e}{0}$ \\
8 & $10^{4}$ & $8.45\mathrm{e}{0}$ & $3.271\mathrm{e}{\mbox{-}1}$ & $3.001\mathrm{e}{\mbox{-}1}$ & $3.7\mathrm{e}{\mbox{-}4}$ & $-6.4\mathrm{e}{\mbox{-}3}$ & $1.1\mathrm{e}{1}$ & $4.1\mathrm{e}{1}$ \\
12 & $10^{6}$ & $2.46\mathrm{e}{1}$ & $1.736\mathrm{e}{\mbox{-}1}$ & $1.688\mathrm{e}{\mbox{-}1}$ & $1.2\mathrm{e}{\mbox{-}5}$ & $-6.2\mathrm{e}{\mbox{-}4}$ & $1.0\mathrm{e}{2}$ & $4.1\mathrm{e}{2}$ \\
16 & $10^{8}$ & $7.55\mathrm{e}{1}$ & $9.575\mathrm{e}{\mbox{-}2}$ & $9.490\mathrm{e}{\mbox{-}2}$ & $3.7\mathrm{e}{\mbox{-}7}$ & $-6.1\mathrm{e}{\mbox{-}5}$ & $1.0\mathrm{e}{3}$ & $4.1\mathrm{e}{3}$ \\
20 & $10^{10}$ & $2.37\mathrm{e}{2}$ & $5.352\mathrm{e}{\mbox{-}2}$ & $5.336\mathrm{e}{\mbox{-}2}$ & $1.2\mathrm{e}{\mbox{-}8}$ & $-6.1\mathrm{e}{\mbox{-}6}$ & $1.0\mathrm{e}{4}$ & $4.1\mathrm{e}{4}$ \\
21 & $10^{10}$ & $3.52\mathrm{e}{2}$ & $4.897\mathrm{e}{\mbox{-}2}$ & $4.885\mathrm{e}{\mbox{-}2}$ & $6.8\mathrm{e}{\mbox{-}9}$ & $-4.3\mathrm{e}{\mbox{-}6}$ & $7.1\mathrm{e}{3}$ & $2.8\mathrm{e}{4}$ \\
24 & $10^{12}$ & $7.46\mathrm{e}{2}$ & $3.004\mathrm{e}{\mbox{-}2}$ & $3.001\mathrm{e}{\mbox{-}2}$ & $3.7\mathrm{e}{\mbox{-}10}$ & $-6.1\mathrm{e}{\mbox{-}7}$ & $1.0\mathrm{e}{5}$ & $4.1\mathrm{e}{5}$ \\
28 & $10^{14}$ & $2.36\mathrm{e}{3}$ & $1.688\mathrm{e}{\mbox{-}2}$ & $1.688\mathrm{e}{\mbox{-}2}$ & $1.2\mathrm{e}{\mbox{-}11}$ & $-6.1\mathrm{e}{\mbox{-}8}$ & $1.0\mathrm{e}{6}$ & $4.1\mathrm{e}{6}$ \\
32 & $10^{16}$ & $7.45\mathrm{e}{3}$ & $9.491\mathrm{e}{\mbox{-}3}$ & $9.490\mathrm{e}{\mbox{-}3}$ & $3.7\mathrm{e}{\mbox{-}13}$ & $-6.1\mathrm{e}{\mbox{-}9}$ & $1.0\mathrm{e}{7}$ & $4.1\mathrm{e}{7}$ \\
36 & $10^{18}$ & $2.36\mathrm{e}{4}$ & $5.337\mathrm{e}{\mbox{-}3}$ & $5.336\mathrm{e}{\mbox{-}3}$ & $1.2\mathrm{e}{\mbox{-}14}$ & $-6.1\mathrm{e}{\mbox{-}10}$ & $1.0\mathrm{e}{8}$ & $4.1\mathrm{e}{8}$ \\
40 & $10^{20}$ & $7.45\mathrm{e}{4}$ & $3.001\mathrm{e}{\mbox{-}3}$ & $3.001\mathrm{e}{\mbox{-}3}$ & $3.7\mathrm{e}{\mbox{-}16}$ & $-6.1\mathrm{e}{\mbox{-}11}$ & $1.0\mathrm{e}{9}$ & $4.1\mathrm{e}{9}$ \\
\hline
\end{tabular}
\caption{Augmented Lagrangian method with telescopic subproblems applied to
problem~(\ref{pccnl}) with the data~(\ref{cex}). In all the iterations, $x^k_2 = 0$, the
Telescopic Method performed no swaps, and $\bar y^k = (0, 0, x^k_3)$.}
\label{tabal}
\end{table}

\subsection{Unbounded penalty parameters: properties of the limit points} \label{alsecunb}

\begingroup
When the penalty parameters tend to infinity, Theorem~\ref{teoal} does not apply, and the
counterexample of Section~\ref{alseccex} shows that the feasible limit points of
Algorithm~\ref{alsec}.1 need not be B-stationary. In this subsection we show what can be said
in this case. Infeasible limit points are stationary for the infeasibility measure on a face
(Proposition~\ref{propinf}). Feasible limit points at which PCC-LICQ holds are KKT points of the
branch problem defined by a face and satisfy an additional condition, given by condition (c) of
Algorithm~\ref{alsec}.1, which, in the MPCC case, means M-stationarity
(Theorem~\ref{teounb} and Corollary~\ref{corounb}). Proposition~\ref{propcexm} shows that,
without condition (c), the limit points may be A-stationary and not M-stationary.

Augmented Lagrangian methods applied to MPCCs written as nonlinear programs converge, under
MPCC-tailored constraint qualifications, to C-stationary points \cite{isu12,ahss19}.
M-stationarity has been obtained with second-order information in the subproblems
\cite{ass18} and, with first-order methods, by keeping the complementarity constraints in the
subproblems \cite{gd22,jkmw23,ars25}. In \cite{ars25}, the subproblems of the method of
\cite{gd22} are solved up to approximate M-stationarity, measured with the limiting normal cone
of the complementarity set at the current iterate, and the feasible limit points are
M-stationary under constraint qualifications weaker than MPEC-LICQ. For cardinality
constraints, an augmented Lagrangian method that keeps the cardinality constraints in the
subproblems was analyzed in \cite{krs21}. Algorithm~\ref{alsec}.1 also keeps the
complementarity (and, in general, the cardinality) constraints in the subproblems. Conditions
(a) and (b) alone do not imply approximate M-stationarity of the subproblems, because (a)
concerns only the current face and the decrease guaranteed by the escaping procedure in (b),
given by~(\ref{escdec}), is of order $1/\rho_k^2$, since the Lipschitz constant of
$\nabla L_k$ grows like $\rho_k$. Condition (c) fills this gap.

Let $\Psi(x) = \frac{1}{2} \left( \|h(x)\|^2 + \|c(x)_+\|^2 \right)$ be the infeasibility
measure.
\endgroup

\begin{prop} \label{propinf}
Let Assumption~\ref{aal} hold and let $\{x^k\}$ be generated by Algorithm~\ref{alsec}.1. If $\{\rho_k\}$ is bounded, every
limit point is feasible. If $\rho_k \to \infty$, $\lim_{k \in K} x^k = x^*$ and $I^k = I$ for
all $k \in K$, then $g^\Psi_I(x^*) = 0$, that is, $x^*$ is stationary for the minimization of
the infeasibility measure $\Psi$ on the face $F_I$.
\end{prop}

\begin{proof}
If $\{\rho_k\}$ is bounded, the feasibility of the limit points follows as in the first
paragraph of the proof of Theorem~\ref{teoal}, which only uses Step~2. Assume that
$\rho_k \to \infty$, and let $\Phi_k = L_k/\rho_k$, so that
\[
\nabla \Phi_k(x) = \frac{\nabla f(x)}{\rho_k} + \nabla h(x) \left( h(x) + \frac{\bar\lambda^k}{\rho_k} \right)
+ \nabla c(x) \left( c(x) + \frac{\bar\mu^k}{\rho_k} \right)_+ .
\]
Since the safeguarded multipliers are bounded and $f$, $h$ and $c$ are continuously
differentiable on the compact set $\B$, $\nabla \Phi_k \to \nabla \Psi$ uniformly on $\B$.
Since, for an interval $[a, b]$, $x \in [a, b]$ and $d \in \R$, $|P_{[a,b]}(x - t d) - x|$ is non-decreasing in $t \geq 0$, taking $t = 1/\rho_k \leq 1$ and $d = \partial L_k(x^k)/\partial x_{i,\nu}$,
this gives $|g^{\Phi_k}_{i,\nu}(x^k)| \leq |g^{L_k}_{i,\nu}(x^k)|$, so that, by (a),
$\|g^{\Phi_k}_I(x^k)\|_\infty \leq \varepsilon_k$ for $k \in K$ large. By the uniform
convergence and the continuity of the projection, $g^\Psi_I(x^*) = 0$.
\end{proof}

Proposition~\ref{propinf} is the analogue, on a face, of the usual result for augmented
Lagrangian methods at infeasible limit points. We do not know whether $x^*$ is B-stationary
for the minimization of $\Psi$ on $\Om$. The argument of Theorem~\ref{teoal} does not apply, because,
on the scale of $\Phi_k$, the escaping procedure only guarantees a decrease of order
$1/\rho_k$, while setting the leaving variable $x_{i,\nu}$ to zero may cost a quantity of
order $x_{i,\nu}^2$, as in Proposition~\ref{propcex}.

\begingroup
In the rest of this subsection, the inequality constraints are written with slack variables,
as in Section~\ref{probl}, so that $q = 0$ and
$\nabla L_k(x^k) = \nabla f(x^k) + J(x^k)^T \lambda^{k+1} = \nabla \ell_{\lambda^{k+1}}(x^k)$.

\begin{teo} \label{teounb}
Let Assumption~\ref{aal} hold, with $q = 0$, and let $\{x^k\}$ be generated by
Algorithm~\ref{alsec}.1. Let $x^*$ be a feasible limit point at which PCC-LICQ holds, and let
$K$ be an infinite set and $I$ a face index such that $\lim_{k \in K} x^k = x^*$ and $I^k = I$
for all $k \in K$. Then $\lim_{k \in K} \lambda^{k+1} = \lambda^*$, where $\lambda^*$ is the
unique solution of $\nabla_\F f(x^*) + J_\F(x^*)^T \lambda = 0$, with $\F = \F(x^*)$, and
$g^{\ell^*}_I(x^*) = 0$, where $\ell^* = \ell_{\lambda^*}$, that is, $x^*$ is a KKT point of the
branch problem of minimizing $f$ subject to $h(x) = 0$ and $x \in F_I$. Moreover, for every
$\nu \in \Nc$, $i \in I_\nu$ with $x^*_{i,\nu} = 0$ and $j \notin I_\nu$,
\begin{equation} \label{mprop}
\partial_{i,\nu} \ell^*(x^*) = 0 \quad \mbox{or} \quad g^{\ell^*}_{j,\nu}(x^*) = 0.
\end{equation}
\end{teo}

\begin{proof}
Let $v \in \F$. Since $x^* \in F_I$ and $x^*_v \neq 0$ if $v$ belongs to a block of $\Nc$, the
variable $x_v$ is not fixed at zero in $F_I$. Since $a_v < x^*_v < b_v$, for $k \in K$ large
enough $x^k_v$ is at distance larger than $\varepsilon_k$ from its bounds, so that, by (a)
and~(\ref{gdef}), $|\partial_v L_k(x^k)| \leq \varepsilon_k$. Therefore,
$\lim_{k \in K} [\nabla_\F f(x^k) + J_\F(x^k)^T \lambda^{k+1}] = 0$. By PCC-LICQ,
$J_\F(x^k)$ has rank $m$ for $k \in K$ large enough and $J_\F(x^k) \to J_\F(x^*)$, so that
$\{\lambda^{k+1}\}_{k \in K}$ converges to the unique solution $\lambda^*$ of
$\nabla_\F f(x^*) + J_\F(x^*)^T \lambda = 0$. Hence,
$\lim_{k \in K} \nabla L_k(x^k) = \nabla \ell^*(x^*)$ and, by (a) and the continuity of the
projection, $g^{\ell^*}_I(x^*) = 0$.

Let $\nu \in \Nc$, $i \in I_\nu$ with $x^*_{i,\nu} = 0$, and $j \notin I_\nu$. If
$|x^k_{i,\nu}| > \min\{\varepsilon_k, \bar\varepsilon\}$ for infinitely many $k \in K$, then,
for these $k$ large enough, $x^k_{i,\nu}$ is at distance larger than $\varepsilon_k$ from its
bounds, because the null bound, if any, is at distance $|x^k_{i,\nu}| > \varepsilon_k$ and
the non-null bounds are far from $x^k_{i,\nu}$, which tends to zero. As in the first part of
the proof, $|\partial_{i,\nu} L_k(x^k)| \leq \varepsilon_k$ for these $k$, and
$\partial_{i,\nu} \ell^*(x^*) = 0$. Otherwise,
$|x^k_{i,\nu}| \leq \min\{\varepsilon_k, \bar\varepsilon\}$ for all $k \in K$ large enough
and, by (c), one of the two inequalities of (c) holds for infinitely many $k \in K$. Taking
limits, by the continuity of the projection, one of the two equalities of~(\ref{mprop})
holds.
\end{proof}

The penalty parameters play no role in Theorem~\ref{teounb}, which also holds when
$\{\rho_k\}$ is bounded, but then Theorem~\ref{teoal} gives B-stationarity. Condition (b) is
not used in the proof. In the MPCC case, we use the usual terminology for a feasible point
$x^*$ with a multiplier $\lambda^*$ such that the partial derivatives of $\ell^*$ satisfy the
sign conditions at the bounds for all the variables that are not null components of
complementarity pairs. There is no condition on a null component whose partner is
non-null. For a biactive pair $(x_a, x_b)$, $x^*$ is \emph{A-stationary} if
$\partial_a \ell^*(x^*) \geq 0$ or $\partial_b \ell^*(x^*) \geq 0$, and \emph{M-stationary} if
both derivatives are non-negative or one of them is null. For partitioned cardinality
constraints, (\ref{mprop}) is the natural extension of M-stationarity. It is stronger than
the M-stationarity of the complementarity reformulation of cardinality constraints
\cite{bks,krs21}, which only requires the derivatives with respect to the non-null variables to
vanish.

\begin{coro} \label{corounb}
In the MPCC case, under the assumptions of Theorem~\ref{teounb}, $x^*$ is M-stationary with
multiplier $\lambda^*$.
\end{coro}

\begin{proof}
For the variables that are not null components of pairs, the sign conditions follow
from $g^{\ell^*}_I(x^*) = 0$, because these variables are not fixed at zero in $F_I$. Let
$(x_a, x_b)$ be a biactive pair with $a \in I_\nu$ and $b \notin I_\nu$. Since $x^*_a = 0$ is
the lower bound of $x_a$, $g^{\ell^*}_a(x^*) = 0$ gives $\partial_a \ell^*(x^*) \geq 0$.
By~(\ref{mprop}), either $\partial_a \ell^*(x^*) = 0$ or $g^{\ell^*}_b(x^*) = 0$, which
means $\partial_b \ell^*(x^*) \geq 0$. In both cases, the pair satisfies the condition of
M-stationarity.
\end{proof}

In Example~1, $\nabla h(0) = 0$, so that PCC-LICQ fails at $x^* = 0$ and
$\lambda^{k+1} \to \infty$. At the limit points of the counterexamples of
Propositions~\ref{propcex2}, \ref{propcex3} and~\ref{propcex4} (Section~\ref{nec}), which
satisfy MPEC-LICQ, $\lambda^* = 0$ and the derivatives of $\ell^*$ with respect to the biactive
pair are $0$ and $-1/4$, so that these points are M-stationary, as stated by
Corollary~\ref{corounb}, but not S-stationary. The following proposition shows that, without
condition (c), the feasible limit points of Algorithm~\ref{alsec}.1 may be A-stationary and
not M-stationary, even under MPEC-LICQ and with the multipliers in the interior of the
safeguarding box.

\begin{prop} \label{propcexm}
Let $n = 3$, let $(x_1, x_2)$ be a complementarity pair, let $x_3$ form a block without
cardinality constraints, $\B = [0,1]^2 \times [-1,1]$, $h(x) = x_3 - x_2$, and
$f(x) = x_1 - x_2 + x_3^2/2$. Consider Algorithm~\ref{alsec}.1 without condition (c), with
$\tau = 1/2$, $\gamma = 10$, $\rho_1 \geq 10$, $\lambda_{\min} = -1$, $\lambda_{\max} = 1$,
$\bar\lambda^1 = 1$, $\bar\lambda^{k+1} = P_{[\lambda_{\min}, \lambda_{\max}]}(\lambda^{k+1})$,
arbitrary $\delta$ and $\eta_{\rm br}$, and $\varepsilon_k = 1/k$. For $k \geq 1$, let
$x^k = (0, 0, s_k)$, where $s_k = (1/(2k) - \bar\lambda^k)/\rho_k$, and let $I^k$ be the face
index in which $x_1$ is not fixed at zero. Then $x^k$ satisfies conditions (a) and (b) of
Algorithm~\ref{alsec}.1, $\bar\lambda^{k+1} = \lambda^{k+1} = 1/(2k) \to 0$, which belongs to
the interior of the safeguarding box, $\rho_k \to \infty$, and $x^k \to 0$, which satisfies
MPEC-LICQ and is A-stationary, but not M-stationary.
\end{prop}

\begin{proof}
For $k = 1$, $s_1 = -1/(2\rho_1)$ and $\lambda^2 = \bar\lambda^1 + \rho_1 s_1 = 1/2$. If
$\bar\lambda^k = 1/(2(k-1))$ for some $k \geq 2$, then
$h(x^k) = s_k = -1/(2k(k-1)\rho_k)$ and $\lambda^{k+1} = \bar\lambda^k + \rho_k s_k = 1/(2k)$,
so that, by induction, $\bar\lambda^{k+1} = \lambda^{k+1} = 1/(2k)$ for all $k$. Moreover,
$|s_k| \leq 1/(2\rho_1) < 1$.

\emph{Condition (a).} Since $x^k_1 = 0$ and $\partial_1 L_k(x^k) = 1$,
$g^{L_k}_1(x^k) = 0$. Also,
$\partial_3 L_k(x^k) = s_k + \bar\lambda^k + \rho_k s_k = s_k + \lambda^{k+1}$, whose
absolute value is at most $1/2 \leq \varepsilon_1$ for $k = 1$ and at most
$1/(2k) + 1/(2k(k-1)\rho_1) \leq 1/k = \varepsilon_k$ for $k \geq 2$. Since
$|g^{L_k}_3(x^k)| \leq |\partial_3 L_k(x^k)|$, (a) holds.

\emph{Penalty parameters.} We have $M_k = |s_k|$, $|s_2|/|s_1| = 1/2$ and
$|s_3|/|s_2| = 1/3$, so that $\rho_4 = \rho_3 = \rho_2 = \rho_1$. For $k \geq 4$, let
$r_k = |s_k|/|s_{k-1}| = (k-2)\rho_{k-1}/(k\rho_k)$. Then $r_4 = 1/2$, so that
$\rho_5 = \rho_1$ and, for $k \geq 5$, if $\rho_k = \rho_{k-1}$, then $r_k = (k-2)/k > 1/2$
and $\rho_{k+1} = 10\rho_k$, whereas, if $\rho_k = 10\rho_{k-1}$, then $r_k < 1/2$ and
$\rho_{k+1} = \rho_k$. Hence, $\rho_k = \rho_1 10^{\lfloor (k-4)/2 \rfloor}$ for $k \geq 4$,
and $\rho_k \to \infty$. In particular, $\rho_k \geq 9k/8$ for all $k$.

\emph{Condition (b).} The only swap is $x_1 \to x_2$ and, since $x^k_1 = 0$,
$\bar y = x^k$. The function $\psi(t) = L_k(0, t, s_k)$ is a quadratic with
$\psi''(t) = \rho_k$ and $\psi'(0) = -1 - \bar\lambda^k - \rho_k s_k = -1 - \lambda^{k+1}$.
Therefore, for all $t$,
$\psi(t) \geq \psi(0) - (1 + \lambda^{k+1})^2/(2\rho_k) \geq L_k(x^k) - 9/(8\rho_k)
\geq L_k(x^k) - \varepsilon_k$, and (b) holds.

\emph{The limit point.} Clearly, $x^k \to 0$, which is feasible. Since $\F(0) = \{3\}$ and
$\partial_3 h = 1$, MPEC-LICQ holds at $0$, and $\lambda^* = -\partial_3 f(0) = 0$.
Then $\partial_1 \ell^*(0) = 1$ and $\partial_2 \ell^*(0) = -1$, so that
$0$ is A-stationary, but not M-stationary. It is not B-stationary either, because $(0, t, t)$
is feasible for $t \in [0,1]$ and $f(0, t, t) = -t + t^2/2$.
\end{proof}

In Proposition~\ref{propcexm}, condition (c) fails at $x^k$ for $i = 1$ and $j = 2$, because
$|\partial_1 L_k(x^k)| = 1$ and $|g^{L_k}_2(x^k)| = 1$. The swap $x_1 \to x_2$ decreases
$L_k$ by at most $9/(8\rho_k)$, which is smaller than $\varepsilon_k$ because
$\rho_k \varepsilon_k \to \infty$. Thus, a decrease test with a tolerance that does not depend
on $\rho_k$, such as (b), cannot replace (c).
\endgroup

\section{Augmented Lagrangian with constraint-preserving swaps} \label{cps}

In the Telescopic Method, a swap sets a nonzero variable of a constrained block to zero (the
leaving variable) and lets a null variable of the same block move (the entering variable). In
the augmented Lagrangian method, the swaps are applied to the minimization of $L_k$, and a
swap usually changes the value of $h$. When $\rho_k$ is large, the corresponding change of the
penalty term of $L_k$ may hide a decrease of $f$, and the swap is rejected although it leads to a
better branch of the feasible set. This is what happens in the counterexample of
Section~\ref{alseccex} and in Proposition~\ref{propcex2} (see Figure~\ref{figex2}). The aim of
this section is to define a swap that does not change $h$. A \emph{constraint-preserving swap}
(CPS) performs the same exchange of variables and then corrects the other free variables of the
face, so that the resulting point $w$ satisfies $h(w) = h(y)$, where $y$ is the current point.
Then, $L_k(w) - L_k(y) = f(w) - f(y)$ for all values of the penalty parameter and of the multipliers, and
the swap is accepted or rejected according to the variation of $f$ alone. The correction is
computed by a Newton-type method (Algorithm~\ref{cps}.1), which, under PCC-LICQ, converges
and moves the variables by an amount proportional to the change of $h$ produced by the
exchange (Lemma~\ref{lemrest}). The CPS itself is Algorithm~\ref{cps}.2. In practice, $h(w)$
is computed with a tolerance, and $L_k(w) - L_k(y)$ differs from $f(w) - f(y)$ by a
controlled amount (see~(\ref{key}) below). The augmented Lagrangian method that uses the
constraint-preserving swaps (Algorithm~\ref{cps}.3) is defined in Section~\ref{algal}.

\subsection{Scaled Newton restoration and constraint-preserving swaps} \label{cpssub}

\begin{figure}[ht!]
\centering
\includegraphics[width=0.6\textwidth]{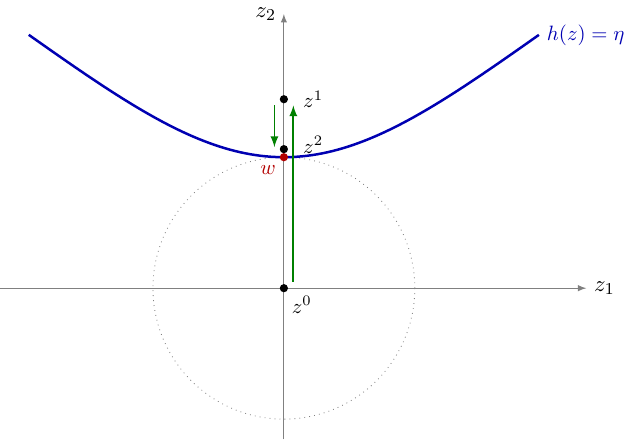}
\caption{Algorithm~\ref{cps}.1 with $S = \{1, 2\}$, $D = I$, $h(z) = e^{z_2} - z_1^2$, $\eta = 2$ and $z^0 = (0, 0)$. The iterates are $z^1 = (0, 1)$, $z^2 \approx (0, 0.736)$, $z^3 \approx (0, 0.694)$, \ldots, and converge to $w = (0, \ln 2)$, which, in this example, is the solution of $h(z) = \eta$ closest to $z^0$ (dotted circle).}
\label{fignewton}
\end{figure}

Algorithm~\ref{cps}.1 below solves the system $h(z) = \eta$, which, in general, has more unknowns (the variables of $S$) than equations, starting from a point $z^0$. Among its solutions, we want one close to $z^0$, so that the swap is changed as little as possible. Each iteration is a Newton step, in which $h$ is linearized at $z^j$ and, among the solutions $\Delta_S$ of $J_S(z^j) \Delta_S = -r^j$, the one that minimizes the scaled norm $\sum_{v \in S} \Delta_v^2 / d_v^2$ (with $\Delta_v = 0$ if $d_v = 0$) is taken. The weights are the (truncated) distances of the variables to their bounds, that is, for $x \in \B$ and $v = 1, \ldots, n$, 
\[
d_v(x) = \min\{1, x_v - a_v, b_v - x_v\} \geq 0,
\]
and $|d_v(x) - d_v(y)| \leq |x_v - y_v|$. The Newton step, $\Delta_S = -D^2 J_S(z^j)^T M_j^{-1} r^j$, is a linear combination of the scaled gradients $D^2 \nabla_S h_i(z^j)$, $i = 1, \ldots, m$, so that the variables at their bounds are not moved and those close to them move little. In Algorithm~\ref{cps}.2, $z^0$ is the point obtained from the current point $y$ by the exchange of variables and $\eta = h(y)$, so that $h(w) = \eta$ means that the swap preserves the values of the constraints. Figure~\ref{fignewton} illustrates the iteration.

\medskip
\noindent
{\bf Algorithm \ref{cps}.1 (Scaled Newton restoration)}

\noindent
Let $z^0 \in \B$, $S \subseteq \{1, \ldots, n\}$, $D = \mathrm{diag}(d_v)_{v \in S}$ with $d_v \in [0,1]$,
$\eta \in \R^m$, a tolerance $\zeta > 0$ and a maximum number of iterations
$j_{\max} \in \{1, 2, \ldots\}$ be given. For $j = 0, 1, 2, \ldots$:

\smallskip
\noindent
{\bf Step 1.} Let $r^j = h(z^j) - \eta$. If $\|r^j\| \leq \zeta$, stop and return
$w = z^j$. If $j = j_{\max}$, stop without returning a point.

\smallskip
\noindent
{\bf Step 2.} Let $M_j = J_S(z^j) D^2 J_S(z^j)^T$. If $M_j$ is singular, stop without returning a point.

\smallskip
\noindent
{\bf Step 3.} Define $z^{j+1}$ by $z^{j+1}_S = z^j_S - D^2 J_S(z^j)^T M_j^{-1} r^j$ and
$z^{j+1}_v = z^j_v$ for $v \notin S$. If $z^{j+1} \notin \B$, stop without returning a point.

\bigskip

The next lemma shows that, near a point $x^*$ at which PCC-LICQ holds, Algorithm~\ref{cps}.1 behaves as Newton's method near a regular solution. If $z^0$ is close to $x^*$ and the initial residual $\|h(z^0) - \eta\|$ is small, the algorithm is well defined, the iterates remain in the box, and the residual decreases quadratically, and at least by half at each iteration. Therefore, a point $w$ with $h(w) = \eta$ (with tolerance $\zeta$) is found within a number of iterations that can be computed in advance. Moreover, $w$ is close to $z^0$, and the variables close to their bounds move very little. The lemma does not assume that $h(z) = \eta$ has a solution. As in Kantorovich's theorem for Newton's method (see \cite[Section~5.3, Theorem~5.3.1]{ds96} and \cite{or70}, and, for Newton's method with generalized inverses for systems with more unknowns than equations, \cite{bi66}), the existence of a solution follows from the nonsingularity of the matrix near $x^*$ and from the smallness of the initial residual.

\begin{lem} \label{lemrest}
Let Assumption~\ref{assum} hold and let $x^* \in \B$ be such that $J_\F(x^*)$ has rank $m$,
where $\F = \F(x^*)$. Then there exist $r > 0$, $\bar r \in (0, 1]$ and $C > 0$ with the
following property. Let $\hat x \in \B$ with $\|\hat x - x^*\| \leq r$, $S \subseteq \{1, \ldots, n\}$ with
$\F \subseteq S$, $d_v = d_v(\hat x)$ for $v \in S$, $z^0 \in \B$ with
$\|z^0 - x^*\| \leq r$ and $z^0_S = \hat x_S$, and $\eta \in \R^m$ with
$\|h(z^0) - \eta\| \leq \bar r$. Then Algorithm~\ref{cps}.1 is well defined and, while it does not stop, $\|r^{j+1}\| \leq C \|r^j\|^2$ and
$\|r^j\| \leq 2^{-j} \|r^0\|$. If $j_{\max} \geq \log_2(1/\zeta)$, it returns a point
$w$ that satisfies $w \in \B$, $\|h(w) - \eta\| \leq \zeta$, $w_v = z^0_v$ for $v \notin S$,
\[
|w_v - z^0_v| \leq C\, d_v(\hat x)^2 \, \|h(z^0) - \eta\| \mbox{ for } v \in S,
\quad \mbox{and} \quad
\|w - z^0\| \leq C \|h(z^0) - \eta\|.
\]
The constants $r$, $\bar r$ and $C$ depend only on $x^*$ and on the problem (through $\sigma$, $\bar d$, $L_J$ and $G_J$, defined in the proof), and not on $\hat x$, $z^0$ or $\eta$.
\end{lem}

\begin{proof}
Let $\sigma > 0$ be the smallest singular value of $J_\F(x^*)$,
$\bar d = \min_{v \in \F} d_v(x^*) > 0$, $L_J$ a Lipschitz constant of $J$ on $\B$, and
$G_J = \max_{x \in \B} \|J(x)\|$. Let $r = \min\{\bar d/2, \sigma/(6 L_J)\}$. If
$\|\hat x - x^*\| \leq r$, then $d_v(\hat x) \geq \bar d/2$ for $v \in \F$, and, if
$z \in \B$ and $\|z - x^*\| \leq 3r$, the smallest singular value of $J_\F(z)$ is at least
$\sigma/2$. For such $z$, since $\F \subseteq S$,
\[
J_S(z) D^2 J_S(z)^T = \sum_{v \in S} d_v^2 J_v(z) J_v(z)^T \succeq \frac{\bar d^2}{4}
J_\F(z) J_\F(z)^T \succeq \left(\frac{\bar d \sigma}{4}\right)^2 I.
\]
Hence the matrix is nonsingular and the norm of its inverse is at most
$c_1 = (4/(\bar d \sigma))^2$. The step $\Delta^j = z^{j+1} - z^j$ satisfies, for $v \in S$,
$|\Delta^j_v| = d_v^2 |J_v(z^j)^T M_j^{-1} r^j| \leq c_2 d_v^2 \|r^j\|$ and, since
$d_v \leq 1$, $\|\Delta^j\| \leq c_2 \|r^j\|$, where $c_2 = G_J c_1$. Since
$J(z^j) \Delta^j = J_S(z^j) \Delta^j_S = -r^j$, if the segment $[z^j, z^{j+1}]$ is contained
in $\B$, then
\[
\|r^{j+1}\| = \|h(z^j + \Delta^j) - h(z^j) - J(z^j) \Delta^j\|
\leq \frac{L_J}{2} \|\Delta^j\|^2 \leq c_3 \|r^j\|^2, \quad c_3 = \frac{L_J c_2^2}{2}.
\]
Let $\bar r = \min\{1, 1/(2 c_3), r/(2 c_2), 1/(2 c_2)\}$ and $C = \max\{2 c_2, c_3\}$. We prove by
induction that, while the algorithm does not stop, $z^j \in \B$, $\|z^j - x^*\| \leq 2r$,
$\|r^j\| \leq 2^{-j} \|r^0\|$ and $|z^j_v - z^0_v| \leq 2 c_2 d_v^2 (1 - 2^{-j}) \|r^0\|$ for
$v \in S$. This is true for $j = 0$. If it is true for $j$, then $M_j$ is nonsingular,
$|z^{j+1}_v - z^0_v| \leq 2 c_2 d_v^2 (1 - 2^{-j}) \|r^0\| + c_2 d_v^2 2^{-j} \|r^0\|
= 2 c_2 d_v^2 (1 - 2^{-j-1}) \|r^0\|$, and $\|z^{j+1} - z^0\| \leq 2 c_2 \|r^0\| \leq r$, so
that $\|z^{j+1} - x^*\| \leq 2r$. Since $2 c_2 \|r^0\| \leq 1$ and $d_v \leq 1$, we have
that $|z^{j+1}_v - z^0_v| \leq d_v$, and $z^0_v = \hat x_v$ is at distance at least
$d_v = d_v(\hat x)$ from its bounds, so that $z^{j+1} \in \B$ and, since $\B$ is convex, the
segment $[z^j, z^{j+1}]$ is contained in $\B$. Thus,
$\|r^{j+1}\| \leq c_3 \|r^j\|^2 \leq c_3 2^{-j} \|r^0\| \, \|r^j\| \leq \|r^j\|/2$, which
completes the induction. Since $\|r^0\| \leq \bar r \leq 1$, the test of
Step~1 holds at the latest at the first $j \geq \log_2(1/\zeta)$, so the algorithm returns a point, and the bounds of the statement follow from the induction.
\end{proof}

As stated by Lemma~\ref{lemrest}, the convergence of Algorithm~\ref{cps}.1 is quadratic. However, the constant $C$ is not known, and the computable bound $j_{\max} \geq \log_2(1/\zeta)$ on the number of iterations, which is used in Algorithm~\ref{cps}.3, comes from the linear bound $\|r^j\| \leq 2^{-j} \|r^0\|$. Lemma~\ref{lemrest} also holds, with a linear rate, for the simplified version of Algorithm~\ref{cps}.1 in which the matrix is not recomputed, that is, $M_j$ is replaced by $M = J_S(\hat x) D^2 J_S(\hat x)^T$ in Step~2 and $J_S(z^j)$ is replaced by $J_S(\hat x)$ in Step~3. (In Algorithm~\ref{cps}.2, $\hat x$ is the point $y$ at which the swap is performed, and this is, up to a small regularization, the version used in our implementation, described in Section~\ref{macmpec}.) In the proof, since $\|\hat x - x^*\| \leq r$, the matrix $M$ is nonsingular with $\|M^{-1}\| \leq c_1$, and $|\Delta^j_v| \leq c_2 d_v^2 \|r^j\|$ as before. Since now $J_S(\hat x) \Delta^j_S = -r^j$, the bound on the residual becomes
\[
\|r^{j+1}\| \leq \frac{L_J}{2} \|\Delta^j\|^2 + L_J \|z^j - \hat x\| \, \|\Delta^j\|
\leq \left( c_3 \|r^j\| + L_J c_2 \|z^j - \hat x\| \right) \|r^j\|.
\]
By the induction hypothesis, $\|r^j\| \leq \|r^0\|$ and $\|z^j - \hat x\| \leq \|z^0 - \hat x\| + 2 c_2 \|r^0\|$. Since $L_J c_2^2 = 2 c_3$, it follows that
\[
\|r^{j+1}\| \leq q \|r^j\|, \quad q = 5 c_3 \|r^0\| + L_J c_2 \|z^0 - \hat x\|.
\]
Since $\|z^0 - \hat x\| \leq 2r$, if $r$ also satisfies $r \leq 1/(12 L_J c_2)$ and $\bar r$ also satisfies $\bar r \leq 1/(15 c_3)$, then $q \leq 1/2$, and the rest of the proof is unchanged. Thus, the convergence is linear, with a ratio $q$ that tends to zero when $r^0$ tends to zero and $z^0$ approaches $\hat x$ (in Algorithm~\ref{cps}.2, $\|z^0 - y\| \leq |y_{i,\nu}| + t$). The bound $j_{\max} \geq \log_2(1/\zeta)$ and the results of Section~\ref{conv} are not affected.

The constraint-preserving swap is defined as follows. Let $y \in F_I$, $\nu \in \Nc$,
$i \in I_\nu$ (the leaving variable) and $j \notin I_\nu$ (the entering variable, so that
$y_{j,\nu} = 0$). Let $J$ be the face index $I$ with $I_\nu$ replaced by
$(I_\nu \setminus \{i\}) \cup \{j\}$, and let
\[
S = \{ (u, \omega) \mid u \in I_\omega, \mbox{ and } y_{u,\omega} \neq 0 \mbox{ or } \omega \notin \Nc \}
\setminus \{(i, \nu)\}
\]
be the set of indices of the variables that are not fixed at zero in $F_I$, other than the leaving one, excluding the null
variables of the constrained blocks (which, therefore, remain null). The entering
variable moves in a direction $\varsigma \in \{-1, +1\}$, which is \emph{admissible} if
$\varsigma = 1$ and $b_{j,\nu} > 0$, or $\varsigma = -1$ and $a_{j,\nu} < 0$. (In the MPCC
case, only $\varsigma = 1$ is admissible.)

In words, Algorithm~\ref{cps}.2 starts from the swap of the Telescopic Method, in which the leaving variable is set to zero and the entering variable is moved by a step $t$, which gives the point $z^0$. This point usually changes the values of the constraints, and Algorithm~\ref{cps}.1 corrects it by moving the other free variables, so that the values of the constraints at the current point are recovered. Then, a sufficient decrease condition is tested on the merit function $\Phi$, which, in Algorithm~\ref{cps}.3, is the augmented Lagrangian. If the correction fails or the decrease is not sufficient, the step of the entering variable is halved (backtracking), which modifies $z^0$, and the process is repeated.

\medskip
\noindent
{\bf Algorithm \ref{cps}.2 (Constraint-preserving swap, CPS)}

\noindent
Let $y$, $\nu$, $i$, $j$, $J$, $S$, and an admissible $\varsigma$ be as above, let $\Phi$ be a
merit function, and let $t_{\max} > 0$, $t_{\min} \in (0, t_{\max}]$, $\alpha > 0$,
$\zeta > 0$ and $j_{\max}$ be given. Let $\eta = h(y)$, $d_v = d_v(y)$ for $v \in S$,
and $t = \min\{t_{\max}, b_{j,\nu}\}$ if $\varsigma = 1$ or $t = \min\{t_{\max},
-a_{j,\nu}\}$ if $\varsigma = -1$.

\smallskip
\noindent
{\bf Step 1.} Let $z^0 = y - y_{i,\nu} e^{i,\nu} + \varsigma t e^{j,\nu}$. Apply
Algorithm~\ref{cps}.1 with $z^0$, $S$, $D = \mathrm{diag}(d_v)$, $\eta$, $\zeta$ and
$j_{\max}$. If it returns a point $w$ and $\Phi(w) \leq \Phi(y) - \alpha t^2$, stop and
return $w$.

\smallskip
\noindent
{\bf Step 2.} Set $t \leftarrow t/2$. If $t < t_{\min}$, stop without returning a point (the
swap is discarded). Otherwise, go to Step~1.

\bigskip

If Algorithm~\ref{cps}.2 returns a point $w$, then $w \in F_J \subseteq \Om$, because $w_{i,\nu} = 0$,
$w_{j,\nu} = \varsigma t \in [a_{j,\nu}, b_{j,\nu}]$, the variables fixed at zero in $F_I$
are not moved, and the variables of $S$ remain in the box. Moreover,
$\|h(w) - h(y)\| \leq \zeta$. For the augmented Lagrangian
$\Phi(x) = f(x) + (\rho/2) \|h(x) + \bar\lambda/\rho\|^2$, writing
$\mu = \bar\lambda + \rho h(y)$ and $e = h(w) - h(y)$, we have the identity
\[
\Phi(w) - \Phi(y) = f(w) - f(y) + \mu^T e + \frac{\rho}{2} \|e\|^2.
\]
Given $\varepsilon > 0$, let
\begin{equation} \label{rhotol}
\zeta(y) = \min\left\{ \frac{\varepsilon}{4 (1 + \|\mu\|)}, \;
\sqrt{\frac{\varepsilon}{2\rho}} \right\}.
\end{equation}
If the tolerance of Algorithm~\ref{cps}.1 satisfies $\zeta \leq \zeta(y)$, then
\begin{equation} \label{key}
\left| \Phi(w) - \Phi(y) - (f(w) - f(y)) \right| \leq \frac{\varepsilon}{2}.
\end{equation}
Thus, up to $\varepsilon/2$, the penalty
parameter and the multipliers disappear from the swap test. Note that $\zeta(y)$ is
computable at $y$. Algorithm~\ref{cps}.2 performs at most
$1 + \log_2(t_{\max}/t_{\min})$ restorations.

\subsection{The modified augmented Lagrangian method} \label{algal}

Algorithm~\ref{cps}.3 is the safeguarded augmented Lagrangian method of Algorithm~\ref{alsec}.1 applied to problem~(\ref{mpcc}), with the same updates of the multipliers and of the penalty parameter. The only difference is in the stopping criterion of the subproblem. Condition (b) of Step~1 concerns the constraint-preserving swaps of Algorithm~\ref{cps}.2 instead of the original swaps, so that the sufficient decrease is tested at points at which the constraints have, up to a tolerance, the same values as at the current point. Condition (c) is not required.

\medskip
\noindent
{\bf Algorithm \ref{cps}.3 (Augmented Lagrangian with telescopic subproblems and CPS)}

\noindent
Let $\tau \in (0,1)$, $\gamma > 1$, $\rho_1 > 0$, $\lambda_{\min} \leq \lambda_{\max}$,
$\bar\lambda^1 \in [\lambda_{\min}, \lambda_{\max}]^m$, $\delta > 0$, $t_{\max} > 0$,
$\alpha > 0$, sequences $\varepsilon_k \downarrow 0$ and $t_{\min,k} \downarrow 0$ with
$t_{\min,k} \leq t_{\max}$, and $x^0 \in \Om$ be given. Set $k \leftarrow 1$.

\smallskip
\noindent
{\bf Step 1 (Subproblem).}
Let $L_k(x) = f(x) + (\rho_k/2) \|h(x) + \bar\lambda^k / \rho_k\|^2$. Starting from $x^{k-1}$,
compute $x^k \in \Om$ and a face index $I^k$ with $x^k \in F_{I^k}$ such that
\begin{itemize}
\item[(a)] $\|g^{L_k}_{I^k}(x^k)\|_\infty \leq \varepsilon_k$, and
\item[(b)] no CPS decreases $L_k$ significantly from $x^k$: for every $\nu \in \Nc$,
  $i \in I^k_\nu$ with $|x^k_{i,\nu}| \leq \delta$, $j \notin I^k_\nu$ and admissible
  $\varsigma$, Algorithm~\ref{cps}.2 applied with $y = x^k$, $\Phi = L_k$,
  $t_{\min} = t_{\min,k}$, tolerance $\zeta = \zeta_k(x^k)$${} > 0$, where $\zeta_k(x^k)$ is~(\ref{rhotol}) with $y = x^k$,
  $\varepsilon = \varepsilon_k$, $\rho = \rho_k$ and $\bar\lambda = \bar\lambda^k$, and
  $j_{\max} = \max\{1, \lceil \log_2(1/\zeta_k(x^k)) \rceil\}$ (or any larger integer), either returns no point or returns a point $w$ with
  $L_k(w) \geq L_k(x^k) - \varepsilon_k$.
\end{itemize}

\smallskip
\noindent
{\bf Step 2 (Penalty parameter).}
If $k > 1$ and $\|h(x^k)\|_\infty > \tau \|h(x^{k-1})\|_\infty$, set $\rho_{k+1} = \gamma \rho_k$. Otherwise, set $\rho_{k+1} = \rho_k$.

\smallskip
\noindent
{\bf Step 3 (Multipliers).}
Set $\lambda^{k+1} = \bar\lambda^k + \rho_k h(x^k)$ and choose
$\bar\lambda^{k+1} \in [\lambda_{\min}, \lambda_{\max}]^m$ (for example,
$P_{[\lambda_{\min}, \lambda_{\max}]^m}(\lambda^{k+1})$, used in our implementation). Set
$k \leftarrow k + 1$ and go to Step~1.

\bigskip

All the inner processes of Step~1 are finite. Algorithm~\ref{cps}.1 performs at most $j_{\max}$ iterations, Algorithm~\ref{cps}.2 at most $1 + \log_2(t_{\max}/t_{\min,k})$ restorations, condition (b) involves finitely many swaps, and, by Lemma~\ref{lemwd}, Step~1 is completed after finitely many iterations of the Telescopic Method. Conditions (a) and (b) are the stopping criterion of the subproblem, and the points $w$ computed
by Algorithm~\ref{cps}.2 in (b) serve only to test $x^k$, and they are not used in Steps~2
and~3. If some CPS returns $w$ with $L_k(w) < L_k(x^k) - \varepsilon_k$, then $x^k$ is not
acceptable, and the method used at Step~1 must go on (in the implementation described below, from $w$ or from a better candidate). Algorithm~\ref{cps}.3 is stated for
problem~(\ref{mpcc}), in which the inequality constraints are written with slack variables,
because the constraint-preserving swaps keep the values of $h$. In our implementation, both
methods use slack variables. With $\lambda_{\min} = \lambda_{\max} = 0$, Algorithm~\ref{cps}.3 is the quadratic penalty
method. Since a CPS that returns a point $w$ satisfies $L_k(w) \leq L_k(x^k) - \alpha t^2$,
condition (b) says that every CPS at $x^k$ that returns a point has $\alpha t^2 < \varepsilon_k$. Note
that $\mu = \lambda^{k+1}$ in~(\ref{rhotol}) when $y = x^k$.

Algorithm~\ref{cps}.2 returns no point when, for every tried step $t \geq t_{\min,k}$, the
restoration (Algorithm~\ref{cps}.1) returns no point or returns a point $w$ with
$L_k(w) > L_k(x^k) - \alpha t^2$. The restoration returns no point when $M_j$ is singular, when
an iterate leaves $\B$, or when the tolerance is not reached in $j_{\max}$ iterations. In this case, the swap is discarded, as in a line
search that does not find a sufficient decrease, and $x^k$ passes the test (b) for this swap.
By Lemma~\ref{lemrest}, if $x^k$ is close to a feasible point at which PCC-LICQ holds and $|x^k_{i,\nu}| + t$ is small, the restoration returns a point, and the swap is discarded
only if the decrease is not sufficient. Since $L_k(w) - L_k(x^k) \approx f(w) - f(x^k)$
by~(\ref{key}), and $\varepsilon_k \to 0$, condition (b) says, asymptotically, that no CPS with
a step $t$ bounded away from zero decreases $f$ by an amount bounded away from zero. This is
what is used in Section~\ref{conv}. If a limit point $x^*$ were not B-stationary, there would
be a descent direction of $f$ along an entering variable, and, for $k$ large, the CPS along it
with a fixed step would return a point with a decrease larger than $\varepsilon_k$,
contradicting (b). Failures of the restoration for large steps do not affect this argument, because Algorithm~\ref{cps}.2 reduces $t$ down to $t_{\min,k}$, which tends to zero, and therefore reaches the region in which Lemma~\ref{lemrest} guarantees that the restoration returns a point.

Step~1 can be implemented by the Telescopic Method (Section~\ref{telsec}) applied to
the minimization of $L_k$ on $\Om$, with the swap search of Step~1 of that algorithm extended
by the CPS. In the notation of Section~\ref{telsec}, the swap search is performed at the point $y^j$
computed by the last iteration of Algorithm MIA (the point $x^k$ of Section~\ref{telsec}, with $j$ in the place of $k$), in the
face $F_{I^j}$. For every $\nu \in \Nc$, $i \in I^j_\nu$ with $|y^j_{i,\nu}| \leq \delta$,
$j' \notin I^j_\nu$ and admissible $\varsigma$, Algorithm~\ref{cps}.2 is applied with
$y = y^j$ and $\zeta = \zeta_k(y^j)$. If it returns a point $w$ with $L_k(w) < f_{best}$, then $w$
becomes the best candidate, and the face of $w$ becomes the candidate face. The point
$z^j$ from which the next iteration of Algorithm MIA starts (the point $z^k$ of Section~\ref{telsec}) is the best candidate, which belongs to its face and satisfies
$L_k(z^j) \leq f_{best} \leq L_k(y^j)$. The convergence analysis of the Telescopic
Method uses only these two properties of $z^j$. Since Algorithm MIA starts from $z^j$ and does
not increase $L_k$, the values $L_k(y^j)$ do not increase, and the iterations of Algorithm MIA
performed on a fixed face behave, in the limit, as if no swaps had been performed between them
(items~1 and~2 of Lemma~\ref{lemtel}). Therefore, the analysis is not affected when the
CPS are added to the swap search. (The original swaps
may be kept, although they are not needed in the analysis below.) The next lemma shows that, with this extended swap search, Step~1 is well defined.

\begin{lem} \label{lemwd}
Let Assumption~\ref{assum} hold and assume that Algorithm MIA satisfies Assumption~\ref{amia}. For every $k$, the Telescopic Method with the extended swap search, applied to
the minimization of $L_k$ on $\Om$, finds after finitely many iterations a point $y^j$ that
satisfies (a) and (b).
\end{lem}

\begin{proof}
Since $\B$ is bounded, $L_k$ satisfies the assumptions on $\phi$ of Section~\ref{telsec}, and we can
repeat the proof of Lemma~6.5 of \cite{bmtel}, with the points returned by the CPS in
addition to the points computed by the escaping procedure in the swap search of
the Telescopic Method. If, at iteration $j$, the swap search is
interrupted by the brute acceptance test, or a CPS returns a point $w$ with
$L_k(w) < L_k(y^j) - \varepsilon_k$, then $L_k(y^{j+1}) \leq L_k(z^j) \leq f_{best} \leq
L_k(y^j) - \min\{\eta_{\rm br}, \varepsilon_k\}$, because Algorithm MIA does not increase
$L_k$ (Assumption~\ref{amia}). Since $\{L_k(y^j)\}$
is non-increasing and bounded below, this happens only for finitely many $j$. If
$\|g^{L_k}_{I^j}(y^j)\|_\infty > \varepsilon_k$ for all $j$ large enough, the argument of
the proof of Lemma~6.5 of \cite{bmtel} (based on item~2 of Lemma~\ref{lemtel}) leads to
a contradiction. Therefore, for some $j$, the search is not interrupted,
$\|g^{L_k}_{I^j}(y^j)\|_\infty \leq \varepsilon_k$ and every CPS at $y^j$ that returns a point
returns $w$ with $L_k(w) \geq L_k(y^j) - \varepsilon_k$, that is, $y^j$ satisfies (a) and (b).
\end{proof}

\begin{rem} \label{remstep1}
The Telescopic Method may be used at Step~1 of Algorithms~\ref{alsec}.1 and~\ref{cps}.3
because it computes, after finitely many iterations, points that satisfy conditions (a)
and (b), and also (c) in Algorithm~\ref{alsec}.1 (Lemmas~\ref{lemalwd} and~\ref{lemwd}). The convergence results of Sections~\ref{alsec}
and~\ref{conv} use only these conditions, so any method that guarantees them could
be used instead. However, it is not enough to compute B-stationary points of the subproblems,
because condition (b) is not a stationarity condition for the minimization of $L_k$ on $\Om$, but a
test on specific moves (the swaps of the Telescopic Method in Algorithm~\ref{alsec}.1 and the
CPS in Algorithm~\ref{cps}.3), as the counterexample of Proposition~\ref{propcex} shows.
\end{rem}

\section{Convergence and necessity of the modification} \label{conv}

In this section we analyze the modified augmented Lagrangian method (Algorithm~\ref{cps}.3).
In Section~\ref{convb} we prove that its feasible limit points at which PCC-LICQ holds are
B-stationary, whatever the behavior of the penalty parameters, and we consider infeasible
limit points. In Section~\ref{nec} we show, by means of counterexamples, that the modification
is necessary, even under MPEC-LICQ, and that MPEC-LICQ cannot be replaced by MPEC-MFCQ.

\subsection{B-stationarity of the limit points} \label{convb}

In this subsection we prove the main result of the paper (Theorem~\ref{teoconv}), state its
consequence for MPCCs (Corollary~\ref{coroMPCC}), and consider infeasible limit points
(Proposition~\ref{propinf}).

\begin{teo} \label{teoconv}
Let Assumption~\ref{assum} hold, let $\{x^k\}$ be generated by Algorithm~\ref{cps}.3, and
let $x^*$ be a limit point of $\{x^k\}$ such that $h(x^*) = 0$ and PCC-LICQ holds at $x^*$.
Let $K$ be an infinite set such that $\lim_{k \in K} x^k = x^*$. Then
$\lim_{k \in K} \lambda^{k+1} = \lambda^*$, where $\lambda^*$ is the unique solution of
$\nabla_\F f(x^*) + J_\F(x^*)^T \lambda = 0$, $\F = \F(x^*)$, and $x^*$ is stationary for
the minimization of $\ell^* = f + \lambda^{*T} h$ on $\Om$. Consequently, $x^*$ is
B-stationary for~(\ref{mpcc}). No assumption is made on $\{\rho_k\}$.
\end{teo}

\begin{proof}
Since there are finitely many face indices, we may assume that $I^k = I$ for all $k \in K$, and then $x^* \in F_I$ because $F_I$ is closed. Let
$\bar d = \min_{v \in \F} d_v(x^*) > 0$.

\emph{Step 1: multipliers.} We have that
$\nabla L_k(x) = \nabla f(x) + J(x)^T (\bar\lambda^k + \rho_k h(x))$, so that
$\nabla L_k(x^k) = \nabla f(x^k) + J(x^k)^T \lambda^{k+1}$. Let $v \in \F$. Then $x_v$ is not fixed at zero in $F_I$ (either $x^*_v \neq 0$ or $v$ belongs to a block $\nu \notin \Nc$) and, for $k \in K$ large enough, $x^k_v$ is at distance at least $\bar d/2$
from its bounds. Since
$g^{L_k}_v(x^k) = P_{[a_v, b_v]}(x^k_v - \partial_v L_k(x^k)) - x^k_v$, we have that
$|g^{L_k}_v(x^k)| = |\partial_v L_k(x^k)|$ if $|\partial_v L_k(x^k)| \leq \bar d/2$ and
$|g^{L_k}_v(x^k)| \geq \bar d/2$ otherwise. By (a), for $k \in K$ large enough
(with $\varepsilon_k < \bar d/2$), $|\partial_v L_k(x^k)| \leq \varepsilon_k$. Thus,
$\nabla_\F f(x^k) + J_\F(x^k)^T \lambda^{k+1} \to 0$ for $k \in K$. If
$\{\lambda^{k+1}\}_{k \in K}$ were unbounded, dividing by $\|\lambda^{k+1}\|$ and taking a
convergent subsequence we would obtain a unit vector $u$ with $J_\F(x^*)^T u = 0$, which
contradicts PCC-LICQ. Therefore, $\{\lambda^{k+1}\}_{k \in K}$ is bounded and each of its
limit points solves $\nabla_\F f(x^*) + J_\F(x^*)^T \lambda = 0$, whose unique solution is
$\lambda^*$. Hence, $\lim_{k \in K} \lambda^{k+1} = \lambda^*$,
$\lim_{k \in K} \nabla L_k(x^k) = \nabla \ell^*(x^*)$, and, by (a) and the continuity of the
projection, $g^{\ell^*}_I(x^*) = 0$. In particular, $\partial_v \ell^*(x^*) = 0$ for every
$v$ such that $x_v$ is not fixed at zero in $F_I$ and $a_v < x^*_v < b_v$. These indices include $\F$ and the indices of the null variables of blocks of $\Nc$ that are not fixed at zero in $F_I$ and whose bounds satisfy $a_v < 0 < b_v$.

\emph{Step 2: stationarity.} As in the proof of item~3 of Lemma~\ref{lemtel},
$g^{\ell^*}_I(x^*) = 0$ implies condition~(\ref{s1}) for $\ell^*$,
and condition~(\ref{s2}) for the blocks without cardinality constraints. Let $\nu \in \Nc$ be such
that $\|\vx^*_\nu\|_0 < p_\nu$. As in the proof of item~3 of Lemma~\ref{lemtel}, there exists
$i \in I_\nu$ such that $x^*_{i,\nu} = 0$, and it remains to prove that
$g^{\ell^*}_{j,\nu}(x^*) = 0$ for all $j \notin I_\nu$. Assume, by contradiction, that
$g^{\ell^*}_{j,\nu}(x^*) \neq 0$ for some $j \notin I_\nu$. Since $x^*_{j,\nu} = 0 \in
[a_{j,\nu}, b_{j,\nu}]$, this means that $\kappa = |\partial_{j,\nu} \ell^*(x^*)| > 0$ and that
$\varsigma = -\mathrm{sign}(\partial_{j,\nu} \ell^*(x^*))$ is admissible.

\emph{Step 3: the CPS returns a point with a large decrease.} Let $S_k$ be the set $S$ defined before Algorithm~\ref{cps}.2 for the swap $(\nu, i, j)$ at $y = x^k$, that is, the indices of the variables that may be moved by the restoration, and let $S = S_k$. For $k \in K$ large enough, $\F \subseteq S_k$, because
$x^*_{i,\nu} = 0$ and the variables $x_v$, $v \in \F$, that belong to blocks of $\Nc$ are non-null at
$x^*$ and, hence, at $x^k$. Let $r$, $\bar r$ and $C$ be the constants of Lemma~\ref{lemrest}, $G$ an
upper bound of $\|\nabla \ell^*\|_\infty$ and of $\|J\|$ on $\B$, and $L_\ell$ a Lipschitz
constant of $\nabla \ell^*$ on $\B$. Let $t_0$ be the initial value of $t$ in
Algorithm~\ref{cps}.2 for $(\nu, i, j, \varsigma)$, which does not depend on $k$. For
$k \in K$ and $t \in (0, t_0]$, let
$z^0(t) = x^k - x^k_{i,\nu} e^{i,\nu} + \varsigma t e^{j,\nu}$ (recall that $x^k_{j,\nu} = 0$). Then $z^0(t)_S = x^k_S$, $\|z^0(t) - x^k\| \leq |x^k_{i,\nu}| + t$, and
$\|h(z^0(t)) - h(x^k)\| \leq G (|x^k_{i,\nu}| + t)$. Let
\[
\hat t = \min\left\{ t_0, \; \frac{r}{2}, \; \frac{\bar r}{2G}, \;
\frac{\kappa}{8 L_\ell C_1^2}, \; \frac{\kappa}{16\alpha} \right\},
\quad \mbox{where} \quad C_1 = 1 + C G.
\]
Let $k \in K$ be large enough so that $\|x^k - x^*\| + |x^k_{i,\nu}| \leq r/2$,
$G |x^k_{i,\nu}| \leq \bar r/2$, and $|x^k_{i,\nu}| \leq \delta$. Then, for all
$t \in (0, \hat t]$, Lemma~\ref{lemrest} applies with $\hat x = x^k$, $z^0 = z^0(t)$,
$\eta = h(x^k)$, $\zeta = \zeta_k(x^k)$ and $j_{\max} \geq \log_2(1/\zeta)$:
Algorithm~\ref{cps}.1 returns a point $w = w(t)$ that satisfies
$\|h(w) - h(x^k)\| \leq \zeta_k(x^k) \leq \varepsilon_k/4$,
$w_v = x^k_v$ for $v \notin S \cup \{(i,\nu), (j,\nu)\}$,
$|w_v - x^k_v| \leq C d_v(x^k)^2 G (|x^k_{i,\nu}| + t)$ for $v \in S$, and
$\|w - x^k\| \leq C_1 (|x^k_{i,\nu}| + t)$.

We have that
$f(w) - f(x^k) = \ell^*(w) - \ell^*(x^k) - \lambda^{*T}(h(w) - h(x^k))
\leq \ell^*(w) - \ell^*(x^k) + \|\lambda^*\| \varepsilon_k/4$. By the Lipschitz
continuity of $\nabla \ell^*$,
$\ell^*(w) - \ell^*(x^k) \leq \nabla \ell^*(x^k)^T (w - x^k) + (L_\ell/2) \|w - x^k\|^2$, and
\begin{align*}
\nabla \ell^*(x^k)^T (w - x^k) = {} & - x^k_{i,\nu} \, \partial_{i,\nu} \ell^*(x^k)
+ \varsigma t \, \partial_{j,\nu} \ell^*(x^k) \\
& + \sum_{v \in S_1} \partial_v \ell^*(x^k) (w_v - x^k_v)
+ \sum_{v \in S_2} \partial_v \ell^*(x^k) (w_v - x^k_v),
\end{align*}
where $S_1$ is the set of the variables $v \in S$ such that $a_v < x^*_v < b_v$ and
$S_2 = S \setminus S_1$. The first term is at most $G |x^k_{i,\nu}|$. By Step~1 and continuity,
for $k \in K$ large, $\varsigma \partial_{j,\nu} \ell^*(x^k) \leq -\kappa/2$ and
$|\partial_v \ell^*(x^k)| \leq \omega_k$ for $v \in S_1$, with $\omega_k \to 0$. If
$v \in S_2$, then $x^*_v$ is at one of its bounds, so that
$d_v(x^k) \leq |x^k_v - x^*_v| \to 0$. Therefore, there exists $\omega'_k \to 0$ such that the
last two sums are bounded by $\omega'_k (|x^k_{i,\nu}| + t)$. Using
$\|w - x^k\|^2 \leq 2 C_1^2 (|x^k_{i,\nu}|^2 + t^2)$ and~(\ref{key}), we obtain
\[
L_k(w) - L_k(x^k) \leq -\frac{\kappa}{2} t + \omega'_k t + L_\ell C_1^2 t^2 + \beta_k,
\]
where
$\beta_k = (G + \omega'_k) |x^k_{i,\nu}| + L_\ell C_1^2 |x^k_{i,\nu}|^2 + \|\lambda^*\|
\varepsilon_k/4 + \varepsilon_k/2 \to 0$. For $k \in K$ large enough,
$\omega'_k \leq \kappa/8$ and $\beta_k \leq \kappa \hat t/32$. Then, for all
$t \in [\hat t/2, \hat t]$, using $L_\ell C_1^2 t \leq \kappa/8$ and
$\beta_k \leq \kappa t/16$,
\[
L_k(w(t)) - L_k(x^k) \leq -\frac{\kappa}{2} t + \frac{\kappa}{8} t + \frac{\kappa}{8} t
+ \frac{\kappa}{16} t = -\frac{3\kappa}{16} t \leq -\alpha t^2,
\]
because $\alpha t \leq \alpha \hat t \leq \kappa/16$.

\emph{Step 4: contradiction.} For $k \in K$ large, $t_{\min,k} < \hat t/2$. Algorithm~\ref{cps}.2
applied with $y = x^k$ tries the values $t_0 2^{-q}$, $q = 0, 1, \ldots$, while
$t_0 2^{-q} \geq t_{\min,k}$, and, since $t_0 \geq \hat t$, one of them belongs to
$[\hat t/2, \hat t]$. Therefore, the algorithm returns a point with some $t \geq \hat t/2$, either at a
value $t > \hat t$ or, at the latest, at the first value in $[\hat t/2, \hat t]$, by Step~3.
The returned point satisfies $L_k(w) \leq L_k(x^k) - \alpha t^2 \leq L_k(x^k) - \alpha \hat t^2/4$.
Since $\varepsilon_k < \alpha \hat t^2/4$ for $k$ large, this contradicts (b). Hence,
$g^{\ell^*}_{j,\nu}(x^*) = 0$ for all $j \notin I_\nu$, and $x^*$ is stationary for the
minimization of $\ell^*$ on $\Om$. By Proposition~\ref{propbstat}, $x^*$ is B-stationary
for~(\ref{mpcc}).
\end{proof}

\begin{coro} \label{coroMPCC}
In the MPCC case, every limit point $x^*$ of the sequence generated by
Algorithm~\ref{cps}.3 such that $h(x^*) = 0$ and MPEC-LICQ holds at $x^*$ is B-stationary
and S-stationary, with multiplier $\lambda^* = \lim_{k \in K} \lambda^{k+1}$ for every infinite set $K$ such
that $\lim_{k \in K} x^k = x^*$.
\end{coro}

\begin{rem}
\begin{enumerate}
\item The penalty parameters appear in the proof only through~(\ref{key}), which is
  guaranteed by the computable tolerance~(\ref{rhotol}) of the restoration. The multiplier
  estimates $\lambda^{k+1}$ converge by PCC-LICQ, independently of $\rho_k$ and of
  the safeguarding box.
\item The scaling $d_v(x^k)^2$ is what makes the variables that tend to their bounds
  irrelevant in Step~3, because they are moved by $O(d_v(x^k)^2 t)$ and, therefore, their
  (possibly non-null) partial derivatives do not spoil the decrease. Without scaling, the
  restoration could move a variable that tends to a bound at which the derivative of $\ell^*$
  is positive, or take it out of the box. Variables of constrained blocks that tend to zero
  in the interior of their bounds may be moved by $O(t)$, but the derivative of $\ell^*$ with
  respect to them vanishes at $x^*$.
\item The tolerance~(\ref{rhotol}) behaves like $\sqrt{\varepsilon_k/\rho_k}$ when
  $\rho_k$ is large. In floating point arithmetic it cannot be smaller than the accuracy with
  which $h$ is evaluated, and in our tests, with $\varepsilon_k = 10^{-8}$ and $\rho_k \leq 10^{20}$,
  it is at least of the order of $10^{-14}$.
\item When $\rho_k \to \infty$, the guarantees of Algorithm~\ref{alsec}.1 (without CPS) are
  weaker. Its feasible limit points at which PCC-LICQ holds satisfy~(\ref{mprop}), that is,
  they are M-stationary in the MPCC case (Theorem~\ref{teounb} and Corollary~\ref{corounb}),
  but they need not be B-stationary, even under MPEC-LICQ and with exact solutions of the
  subproblems (Propositions~\ref{propcex2}--\ref{propcex4}).
\end{enumerate}
\end{rem}

As usual for augmented Lagrangian methods, the limit points may be infeasible. Proposition~\ref{propinf} (Section~\ref{alsecunb}) also holds for Algorithm~\ref{cps}.3, under Assumption~\ref{assum} and without the constraints $c$, with the same proof, which uses only Step~2 and condition (a). The constraint-preserving swaps, which keep $h$
constant, cannot decrease the infeasibility measure $\Psi$.

\subsection{Is the modification necessary? Can MPEC-LICQ be weakened?} \label{nec}

Theorem~\ref{teoal} (original swaps) requires bounded penalty parameters but no
constraint qualification, while Theorem~\ref{teoconv} (CPS) requires PCC-LICQ but nothing
on the penalty parameters. In this subsection we consider two questions in the MPCC case. Does
the original method fail under MPEC-LICQ? Can MPEC-LICQ be replaced by MPEC-MFCQ in
Theorem~\ref{teoconv}? The answer to the first question is yes, even with exact solutions of the subproblems and multipliers in the
interior of the safeguarding box (counterexamples of Propositions~\ref{propcex2},
\ref{propcex3} and~\ref{propcex4}). The answer to the second one is no, at least with swaps in a single block
(counterexample of Proposition~\ref{propmfcq}).

\paragraph{Example 2.}
As in Example~1 (Section~\ref{alseccex}), $n = 3$,
$\vx_1 = (x_1, x_2)$ is a complementarity pair, $x_3$ forms a block without cardinality
constraints, $x_1 \in [0, 2]$, $x_2 \in [0, 1]$, but with $x_3 \in [-2, 2]$, $h(x) = x_1 - x_3$, and $f(x) = -x_2/4 - x_1^2/2 +
x_3^2/2 + x_3^4/4$. The feasible set is
$X = \{(t, 0, t) \mid t \in [0,2]\} \cup \{(0, t, 0) \mid t \in [0,1]\}$, and the global
minimizer is $(0, 1, 0)$, with $f = -1/4$. MPEC-LICQ holds at every feasible point other than
$(2, 0, 2)$, because $\partial h/\partial x_3 = -1$ and $x_3$ is free at these points. The following
proposition shows that the quadratic penalty method with the original swaps (Algorithm~\ref{alsec}.1 with
$\lambda_{\min} = \lambda_{\max} = 0$) may converge to a point that is not
B-stationary even under MPEC-LICQ.

\begin{prop} \label{propcex2}
Consider Example~2 and Algorithm~\ref{alsec}.1 with $\bar\lambda^k = 0$ for all $k$ and
$\rho_1 \geq 3$. For each $k$, let $s_k = (\rho_k - 1)^{-1/2}$, $t_k = s_k \rho_k / (\rho_k - 1)$,
and $x^k = (t_k, 0, s_k)$. Then $x^k$ is a strict local minimizer of $L_k$ on $\Om$ and
satisfies conditions (a), (b) and (c) of Algorithm~\ref{alsec}.1 with $\varepsilon_k = 0$, also
with the look-ahead strategy of our implementation. Moreover, $\rho_k \to \infty$, $x^k \to 0$, and $0$
is not B-stationary, although MPEC-LICQ holds at $0$.
\end{prop}

\begin{proof}
Here $L_k(x) = f(x) + \rho_k (x_1 - x_3)^2/2$. Let $\rho = \rho_k$, $s = s_k$ and $t = t_k$.
Then $t - s = s/(\rho - 1)$, $\partial L_k(x^k)/\partial x_1 = -t + \rho(t - s) = 0$ and
$\partial L_k(x^k)/\partial x_3 = s + s^3 - \rho(t - s) = s + s^3 - t = s^3 - s/(\rho-1) = 0$.
Since $0 < t, s < 2$ (because $\rho \geq 3$), (a) holds and, since $\partial L_k(x^k)/\partial x_1 = 0$, so does (c). Near $x^k$, $x_1 > 0$ forces
$x_2 = 0$, and the Hessian of $L_k$ with respect to $(x_1, x_3)$ is
$\left(\begin{smallmatrix} \rho - 1 & -\rho \\ -\rho & \rho + 1 + 3s^2 \end{smallmatrix}\right)$,
with determinant $3 s^2 \rho - 1 - 3 s^2 = 2$ and positive diagonal, so $x^k$ is a strict
local minimizer of $L_k$ on $\Om$. The only swap is $x_1 \to x_2$, with
$\bar y = (0, 0, s)$, and, since $L_k(\bar y + u e^2) = L_k(\bar y) - u/4$, every point $w$ produced by the escaping procedure satisfies $L_k(w) \geq L_k(\bar y) - 1/4$. Moreover,
\[
L_k(\bar y) - L_k(x^k) = \frac{t^2}{2} + \frac{\rho}{2}\left(s^2 - (t-s)^2\right)
\geq \frac{\rho}{2(\rho - 1)} - \frac{\rho}{2(\rho - 1)^3} \geq \frac{1}{2} - \frac{3}{16},
\]
because $\rho/(\rho-1) > 1$ and $\rho/(\rho-1)^3 \leq 3/8$ for $\rho \geq 3$. Hence
$L_k(w) - L_k(x^k) \geq 1/2 - 3/16 - 1/4 = 1/16 > 0$, and (b) holds. Since no swap decreases
$L_k$, there are no candidates for the look-ahead strategy, which in our implementation only
improves swaps that decrease the objective function. Since $h(x^k) = t - s = (\rho_k - 1)^{-3/2}$
depends only on $\rho_k$, the penalty parameter is increased at least every other iteration,
so $\rho_k \to \infty$ and $x^k \to 0$. Finally, $(0, u, 0) \in X$ for $u \in [0,1]$
and $\nabla f(0)^T e^2 = -1/4 < 0$, so $0$ is not B-stationary.
\end{proof}

\begin{figure}[ht!]
\centering
\includegraphics{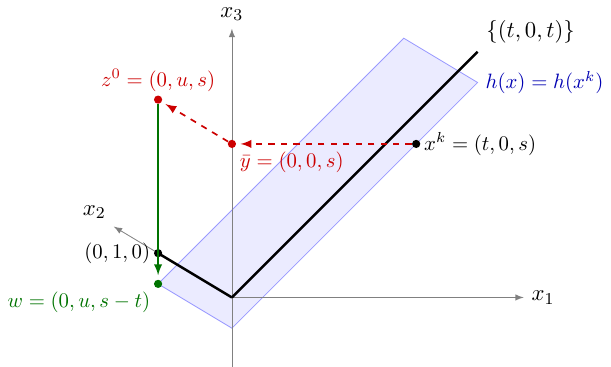}
\caption{Example~2 at an iterate $x^k = (t, 0, s)$ of Proposition~\ref{propcex2}, close to the
branch $\{(t,0,t)\}$ of the feasible set (thick lines, not to scale). The original swap
$x_1 \to x_2$ (dashed) goes to $\bar y$ and then along $x_2$, through points that are far from the
plane $h(x) = h(x^k)$, because $h(\bar y) = -s$, so that $L_k$ includes the penalty
$\rho_k s^2/2$, larger than the decrease $u/4$ of $f$ along $x_2$, and the swap is rejected.
The CPS (Algorithm~\ref{cps}.2) starts from the same point $z^0 = (0, u, s)$ and restores
the free variable $x_3$ with Algorithm~\ref{cps}.1, so that $w$ lies on the plane
$h(x) = h(x^k)$ and $L_k(w) - L_k(x^k) = f(w) - f(x^k) \approx -u/4 < 0$ when $s$ and $t$
are small. Thus, $x^k$ does not satisfy condition (b) of Algorithm~\ref{cps}.3.}
\label{figex2}
\end{figure}

Theorem~\ref{teoconv} applies to Algorithm~\ref{cps}.3 in Example~2, also with
$\lambda_{\min} = \lambda_{\max} = 0$, so that no limit point can be $0$. As Figure~\ref{figex2}
shows, near $0$ the CPS finds the descent along $x_2$ that the original swap misses. In
Proposition~\ref{propcex2}, the safeguarding box is $\{0\}$, so that $\bar\lambda^k = 0$ for all $k$, the method is the quadratic penalty method, and $\lambda^* = 0$ lies on the boundary of
this box. If the safeguarding box contains $\lambda^*$ in its interior (for example,
$\lambda_{\min} = -10^{20}$ and $\lambda_{\max} = 10^{20}$, as in Section~\ref{numer}) and the
subproblems are solved accurately, the augmented Lagrangian method with the original swaps
does not fail in this example. The multiplier estimate $\bar\lambda^2 = \lambda^2 > 0$ makes
the subproblem of the second iteration move to the branch $x_1 = 0$ (see
Section~\ref{numer}). In Proposition~\ref{propcex4} below we show, with another example, that
the method may fail also with such a safeguarding box and exact solutions of the subproblems.
The next proposition shows,
however, that the original method may fail under MPEC-LICQ with multipliers in the interior of
the safeguarding box if the subproblems are solved inexactly, with tolerances
$\varepsilon_k \to 0$ as allowed by Algorithm~\ref{alsec}.1.

\begin{prop} \label{propcex3}
Consider Example~2 and Algorithm~\ref{alsec}.1 with $\tau = 1/2$, $\gamma = 10$,
$\rho_1 \geq 16$, $\lambda_{\min} < 0$, $\lambda_{\max} \geq 1/2$, $\bar\lambda^1 = 0$,
$\bar\lambda^{k+1} = P_{[\lambda_{\min}, \lambda_{\max}]}(\lambda^{k+1})$, arbitrary $\delta$ and $\eta_{\rm br}$, $\varepsilon_1 = 1$, and
$\varepsilon_k = 8/(k+1)^3 + 1/(k(k+1)\rho_1)$ for $k \geq 2$. Let $t_0 = 0$ and, for
$k \geq 1$, $t_k = 1/(k+1)$, $h_k = (t_k - t_{k-1})/\rho_k$, and
$x^k = (t_k, 0, t_k - h_k)$. Then $x^k$ satisfies conditions (a), (b) and (c) of Algorithm~\ref{alsec}.1, $\bar\lambda^{k+1} = \lambda^{k+1} = t_k \to 0 = \lambda^*$, which belongs to
the interior of the safeguarding box, $\rho_k \to \infty$, and $x^k \to 0$, which satisfies
MPEC-LICQ and is not B-stationary.
\end{prop}

\begin{proof}
Note that $h(x^k) = x^k_1 - x^k_3 = h_k$. By induction, $\bar\lambda^k = t_{k-1}$. This is
true for $k = 1$ and, if it is true for $k$, then
$\lambda^{k+1} = \bar\lambda^k + \rho_k h_k = t_k \in (0, 1/2]$, so that
$\bar\lambda^{k+1} = t_k$. Let $t = t_k$, $s = x^k_3 = t - h_k$, and $\mu = \lambda^{k+1} = t$.

\emph{Condition (a).} We have $\partial L_k(x^k)/\partial x_1 = -t + \mu = 0$ and
$\partial L_k(x^k)/\partial x_3 = s + s^3 - \mu = s^3 - h_k$. For $k = 1$,
$s = 1/2 - 1/(2\rho_1)$ and $|s^3 - h_1| < 1 = \varepsilon_1$. For $k \geq 2$,
$h_k = -1/(k(k+1)\rho_k) < 0$ and $s \leq t + 1/(k(k+1)\rho_1) \leq 2t$, so that
$|s^3 - h_k| \leq 8/(k+1)^3 + 1/(k(k+1)\rho_1) = \varepsilon_k$. Since $x^k_1$ and $x^k_3$ are
at distance at least $1 > \varepsilon_k$ from their bounds (for $k \geq 2$) and $x_2$ is fixed in
the face, (a) holds (for $k = 1$, the projected derivatives are bounded by the derivatives). Since $\partial L_k(x^k)/\partial x_1 = 0$, condition (c) also holds.

\emph{Penalty parameters.} For $k \geq 3$, let
$r_k = |h_k|/|h_{k-1}| = (k-1)\rho_{k-1}/((k+1)\rho_k)$, and note that $|h_2|/|h_1| = 1/3$. Thus
$\rho_3 = \rho_2 = \rho_1$, $r_3 = 1/2$, so that $\rho_4 = \rho_1$, and, for $k \geq 4$, if $\rho_k = \rho_{k-1}$, then $r_k = (k-1)/(k+1) > 1/2$ and $\rho_{k+1} = 10 \rho_k$, whereas, if $\rho_k = 10\rho_{k-1}$, then $r_k < 1/2$ and $\rho_{k+1} = \rho_k$. Hence
$\rho_k = \rho_1 10^{\lfloor (k-3)/2 \rfloor}$ for $k \geq 3$, and $\rho_k \to \infty$. Clearly,
$x^k \to 0$.

\emph{Condition (b).} The only possible swap is $x_1 \to x_2$, with $\bar y = (0, 0, s)$. Since
$L_k(\bar y + u e^2) = L_k(\bar y) - u/4$, every point $w$ produced by the escaping
procedure satisfies $L_k(w) \geq L_k(\bar y) - 1/4$. Using $h(\bar y) = -s$, $s + h_k = t$ and
$s^2 - h_k^2 = t (s - h_k)$,
\[
L_k(\bar y) - L_k(x^k) = \frac{t^2}{2} - \bar\lambda^k t + \frac{\rho_k}{2} t (s - h_k).
\]
For $k = 1$, this is $1/8 + (\rho_1 - 2)/8 \geq 1/4$. For $k \geq 2$, $\bar\lambda^k = 1/k \leq 2t$
and $s - h_k = t - 2h_k \geq t$, so that
$L_k(\bar y) - L_k(x^k) \geq t^2 (\rho_k - 3)/2 \geq 1/4$ if $\rho_k \geq 3 + (k+1)^2/2$. This
holds for $k \leq 4$ because $\rho_k = \rho_1 \geq 16$, and for $k \geq 5$ because
$\rho_k \geq 16 \cdot 10^{(k-4)/2}$, whose ratio between consecutive values, $\sqrt{10}$, is
larger than the one of $3 + (k+1)^2/2$. Therefore, $L_k(w) \geq L_k(x^k)$ and (b) holds.
Finally, MPEC-LICQ holds at $0$, $\lambda^* = 0$ by Proposition~\ref{propbstat}, and $0$ is
not B-stationary, as in Proposition~\ref{propcex2}.
\end{proof}

In Proposition~\ref{propcex3} the tolerances $\varepsilon_k$ are much larger than
$|h(x^k)|$, so that the condition $\varepsilon_k \leq \eta_k \|h(x^k)\|$, with
$\eta_k \to 0$, used in \cite{abms} to prove the boundedness of the penalty parameters, does
not hold. Moreover, the branch problem $\Minimize f$ subject to $h(x) = 0$, $x_2 = 0$, does not
satisfy the second-order sufficient condition at $0$ (on the feasible branch,
$f(t, 0, t) = t^4/4$). On the other
hand, Theorem~\ref{teoconv} holds for every sequence of tolerances $\varepsilon_k \to 0$, and the
points of Proposition~\ref{propcex3} do not satisfy condition (b) of
Algorithm~\ref{cps}.3 for $k$ large, because the CPS returns $w = (0, t, -h_k)$ with
$f(w) - f(x^k) \approx -t/4$.

The following example shows that the original method may fail under MPEC-LICQ also with
multipliers in the interior of the safeguarding box and exact solutions of the subproblems
(strict local minimizers), with the parameters $\tau = 1/2$ and $\gamma = 10$ of Algencan.
As expected from the results of \cite{abms}, the second-order sufficient condition fails at
the limit point (on the feasible branch, $f$ is of order six). The iterates alternate between
two branches of local minimizers of the subproblems, on which $h$ has opposite signs.

\paragraph{Example 3.}
$n = 4$, $\vx_1 = (x_1, x_2)$ is a complementarity pair, $(x_3, x_4)$ forms a block without
cardinality constraints, $\B = [0,1]^2 \times [-1,1]^2$, $h(x) = x_1 - x_4$, and
\[
f(x) = -\frac{x_2}{4} + \frac{1}{6} \left[ \left( \frac{x_1 + x_3}{2} \right)^6
+ \left( \frac{x_1 - x_3}{2} \right)^6 \right] + x_1 x_3 (x_1 - x_4).
\]
MPEC-LICQ holds at $x^* = 0$, because $x_4$ is free and $\partial h / \partial x_4 = -1$, and
$x^*$ is not B-stationary, because $(0, u, 0, 0)$ is feasible for $u \in [0,1]$ and
$\nabla f(0)^T e^2 = -1/4$.

\begin{prop} \label{propcex4}
Consider Example~3 and Algorithm~\ref{alsec}.1 with $\tau = 1/2$, $\gamma = 10$,
$\lambda_{\min} = -1$, $\lambda_{\max} = 1$, $\rho_1 = 100$, $\bar\lambda^1 = 1/100$,
$\bar\lambda^{k+1} = P_{[\lambda_{\min}, \lambda_{\max}]}(\lambda^{k+1})$, and $\varepsilon_k = 0$
for all $k$. For each $k$, let $y_k > 0$ be the positive root of
$\rho_k y^2 - \sqrt{2}\, y - |\bar\lambda^k| = 0$, $s_k = (2 y_k^2)^{1/4}$, and
$x^k = (s_k, 0, s_k, s_k + y_k^2)$ if $\bar\lambda^k > 0$,
$x^k = (s_k, 0, -s_k, s_k - y_k^2)$ if $\bar\lambda^k < 0$.
Then $x^k$ is a strict local minimizer of $L_k$ on $\Om$ that satisfies conditions (a), (b) and (c)
of Algorithm~\ref{alsec}.1 with $\varepsilon_k = 0$, $\bar\lambda^{k+1} = \lambda^{k+1} \to 0$,
which belongs to the interior of the safeguarding box, $\rho_k \to \infty$, and $x^k \to 0$,
which satisfies MPEC-LICQ and is not B-stationary.
\end{prop}

\begin{proof}
In the face $x_2 = 0$, we use the coordinates $t_1 = (x_1 + x_3)/2$, $t_2 = (x_1 - x_3)/2$
and $\omega = h(x) = x_1 - x_4$, in which
$L_k = (t_1^6 + t_2^6)/6 + (t_1^2 - t_2^2)\omega + \bar\lambda^k \omega + \rho_k \omega^2/2$.
Let $\rho = \rho_k$, $y = y_k$, $s = s_k$, and assume that $\bar\lambda^k > 0$ (the other case is
symmetric, with $t_1$ and $t_2$ exchanged). Then $x^k$ corresponds to $t_1 = s$, $t_2 = 0$,
$\omega = -y^2$, with $s^4 = 2y^2$, and
$\partial L_k/\partial t_1 = s(s^4 + 2\omega) = 0$, $\partial L_k/\partial t_2 = 0$, and
$\partial L_k/\partial \omega = s^2 + \bar\lambda^k - \rho y^2 = \sqrt{2}\, y + \bar\lambda^k - \rho y^2 = 0$.
The Hessian is block diagonal, with the entry $5t_2^4 - 2\omega = 2y^2 > 0$ and the block
$\left(\begin{smallmatrix} 8y^2 & 2s \\ 2s & \rho \end{smallmatrix}\right)$ in
$(t_1, \omega)$, whose determinant $8\rho y^2 - 4\sqrt{2}\, y$ is positive because
$\rho y > \sqrt{2}$ (the positive root satisfies $y > \sqrt{2}/\rho$). Since $x^k_1 > 0$,
$\Om$ coincides with the face near $x^k$, so that $x^k$ is a strict local minimizer of $L_k$ on
$\Om$ and (a) holds with $\varepsilon_k = 0$ (the iterates belong to the interior of the box,
see below). Since $x^k_1$ is strictly between its bounds, $\partial L_k(x^k)/\partial x_1 = 0$, so that (c) also holds.

\emph{Multipliers and penalty parameters.} We have
$\lambda^{k+1} = \bar\lambda^k + \rho \omega = -\sqrt{2}\, y$ (and $+\sqrt{2}\, y$ if
$\bar\lambda^k < 0$). Let $a_k = |\bar\lambda^k|$ and $z_k = \rho_k |\lambda^{k+1}|$. Then
$z_k = 1 + \sqrt{1 + 2 \rho_k a_k}$ and $|h(x^k)| = y_k^2 = z_k^2/(2\rho_k^2)$. If
$|\lambda^{k+1}| \leq 1$, then $\bar\lambda^{k+1} = \lambda^{k+1}$, $a_{k+1} = z_k/\rho_k$, and
$z_{k+1} = 1 + \sqrt{1 + 2 (\rho_{k+1}/\rho_k) z_k}$ with
$\rho_{k+1}/\rho_k \in \{1, 10\}$. Since $z \mapsto 1 + \sqrt{1 + 20 z}$ is increasing and
its fixed point is $22$, and $z_1 = 1 + \sqrt{3}$, by induction $2 < z_k \leq 22$ and
$|\lambda^{k+1}| = z_k/\rho_k \leq 22/100 < 1$ for all $k$. Thus the multipliers alternate in
sign, $\bar\lambda^{k+1} = \lambda^{k+1}$, and $|\lambda^{k+1}| \to 0$ if $\rho_k \to \infty$.
If $\{\rho_k\}$ were bounded, then $\rho_k = \bar\rho$ for $k \geq k_1$, $z_{k+1} = 1 + \sqrt{1 + 2 z_k}$, $z_k \to 4$, and $|h(x^k)|/|h(x^{k-1})| = (z_k/z_{k-1})^2 \to 1 > \tau$, so that the penalty parameter would be
increased, a contradiction. Hence $\rho_k \to \infty$, $y_k \to 0$, $s_k \to 0$ and
$x^k \to 0$. Moreover, $y_k = |\lambda^{k+1}|/\sqrt{2} \leq 0.16$ and $s_k \leq 0.48$, so that
the iterates belong to the interior of $\B$.

\emph{Condition (b).} The only swap is $x_1 \to x_2$, with $\bar y = (0, 0, x^k_3, x^k_4)$.
Since $L_k(\bar y + u e^2) = L_k(\bar y) - u/4$, every point $w$ produced by the escaping
procedure satisfies $L_k(w) \geq L_k(\bar y) - 1/4$. Using $s^2 = \sqrt{2}\, y$,
$h(\bar y) - h(x^k) = -s$ and $|\bar\lambda^k| = \rho y^2 - \sqrt{2}\, y$, a direct computation
gives
\[
L_k(\bar y) - L_k(x^k) = \frac{65 \sqrt{2}}{96} y^3 + \frac{\rho y}{\sqrt{2}} \pm \sqrt{2}\, y s,
\]
with the sign $+$ if $\bar\lambda^k > 0$ and $-$ otherwise. Since $\rho y > \sqrt{2}$ and
$\sqrt{2}\, y s = 2^{3/4} y^{3/2} < 0.12$, we obtain $L_k(\bar y) - L_k(x^k) > 1 - 0.12 > 1/4$,
so that $L_k(w) > L_k(x^k)$ and (b) holds.
\end{proof}

The following example shows that Theorem~\ref{teoconv} does not hold if MPEC-LICQ is
replaced by MPEC-MFCQ. The reason is that a descent direction may require two simultaneous
swaps.

\paragraph{Example 4.}
$n = 6$, with complementarity pairs $(y_1, z_1) = (x_1, x_2)$ and $(y_2, z_2) = (x_3, x_4)$, a
block $(w_1, w_2) = (x_5, x_6)$ without cardinality constraints, $\B = [0,1]^6$, and
\[
f(x) = y_1 + y_2 - z_1 - z_2 + 3 w_1 + 3 w_2,
\qquad
h(x) = z_1 - z_2 + w_1 - w_2.
\]
The global minimizer is $(0, 1, 0, 1, 0, 0)$, with $f = -2$. The next proposition shows that the method may fail in this example.

\begin{prop} \label{propmfcq}
In Example~4, MPEC-MFCQ holds at $x^* = 0$, MPEC-LICQ does not hold, and $x^*$ is not
B-stationary. Let $I$ be the face index in which $y_1$, $y_2$, $w_1$ and $w_2$ are the variables that are not fixed at zero in $F_I$. Then, for every $\rho > 0$ and $\bar\lambda$ with
$-3 \leq \bar\lambda \leq 3$, the point $x^* = 0$ satisfies conditions (a) and (b) of
Algorithm~\ref{cps}.3 with $\varepsilon_k = 0$. More strongly, for every single swap
($y_1 \to z_1$ or $y_2 \to z_2$) and every point $w$ of the new face with $h(w) = h(0)$, we have
$f(w) \geq f(0)$. Therefore, the constant sequence $x^k = 0$ (with $\bar\lambda^1 = 0$) is
generated by Algorithm~\ref{cps}.3, and its limit is not B-stationary.
\end{prop}

\begin{proof}
At $x^* = 0$ all the variables are at their lower bounds and both pairs are biactive. The
equality constraints of the tightened problem are $h = 0$ and $y_1 = z_1 = y_2 = z_2 = 0$,
whose gradients are linearly independent, because $\nabla h$ has non-null components in
$w_1$ and $w_2$. The direction $d$ with $d_{w_1} = d_{w_2} = 1$ and the other components null
satisfies $\nabla h^T d = 0$ and is strictly feasible for the active bounds $w \geq 0$. Thus,
MPEC-MFCQ holds. MPEC-LICQ fails because $\F(0) = \emptyset$. The points
$(0, u, 0, u, 0, 0)$, $u \in [0,1]$, are feasible, with $f = -2u$, so $0$ is not
B-stationary. Since $h(0) = 0$, $\nabla L_k(0) = \nabla f(0) + \bar\lambda \nabla h(0)$, whose
components corresponding to the variables that are not fixed at zero in $F_I$ ($y_1$, $y_2$, $w_1$, $w_2$) are $1$, $1$,
$3 + \bar\lambda$ and $3 - \bar\lambda$. Since these variables are at their lower bounds,
$g^{L_k}_I(0) = 0$ and (a) holds. Also, $\lambda^{k+1} = \bar\lambda^k$, so that, with
$\bar\lambda^1 = 0$, $\bar\lambda^k = 0$ and $\rho_k = \rho_1$ for all $k$. Since $h(w) = h(0)$
implies $L_k(w) - L_k(0) = f(w) - f(0)$, it suffices to consider $f$. In the swap $y_1 \to z_1$, the points of
the new face have $z_2 = 0$, and $h(w) = 0$ with $z_1 = t \geq 0$ gives $w_2 = t + w_1$, so that
$f(w) = y_2 - t + 3 w_1 + 3(t + w_1) \geq 2t \geq 0$. The swap $y_2 \to z_2$ is analogous
($w_1 = t + w_2$). Hence no CPS returns a point, whatever the restoration procedure, and (b) holds. (In
Algorithm~\ref{cps}.2, the restoration returns no point because the moving variables $w_1, w_2$ have
$d_v = 0$.)
\end{proof}

A descent direction at $0$ requires $z_1$ and $z_2$ to enter simultaneously. Under
MPEC-LICQ, the free variables can always compensate the variation of $h$ caused by a single
entering variable, which is why swaps in one block at a time suffice. Under MPEC-MFCQ, a
complete test would have to consider simultaneous swaps in several blocks, guided by the
linearization of the problem, which leads to combinatorial subproblems such as the linear
programs with equilibrium constraints of \cite{nl26}. In Example~4, the original swaps of the
Telescopic Method, which do not preserve $h$, do escape from $0$, because, from $\bar y = 0$ along
$z_1$, $L_k(u e^{z_1}) = -u + \rho u^2/2 < 0$ for small $u$. Indeed, our
implementation, which keeps both kinds of swaps, finds the global minimizer in the first outer
iteration (Section~\ref{numer}).

\begin{rem}
In Example~4 (MPEC-MFCQ), the limit
point $x^* = 0$ of Algorithm~\ref{cps}.3 in the counterexample of Proposition~\ref{propmfcq} is A-stationary but not M-stationary, because
M-stationarity of the first pair requires $\lambda \geq 1$, and that of the second, $\lambda \leq -1$. The methods of \cite{gd22,jkmw23} and \cite{ars25}, which also keep the
complementarity constraints in the subproblems, obtain M-stationary limit points under their
constraint qualifications, whereas, with the constraint-preserving swaps, Theorem~\ref{teoconv} gives
B-stationarity under MPEC-LICQ.
\end{rem}

\section{Numerical experiments} \label{numexp}

We report experiments with the examples of Sections~\ref{alsec} and~\ref{conv}
(Section~\ref{numer}), with the MacMPEC collection (Section~\ref{macmpec}), and with problems
with partitioned cardinality constraints (Section~\ref{cardexp}).

We implemented Algorithm~\ref{cps}.3 in Fortran, as a driver around the implementation of
the Telescopic Method of \cite{bmtel}, the subroutine \texttt{pcc\_telescopic}. Algorithm MIA
is Gencan, with exact second derivatives and its default parameters. As in Step~3 of
Algorithm~5.1 of \cite{bmtel}, a single iteration of Gencan is performed on each visited face
($G = 1$, in the notation of Section~\ref{macmpec}). We used $\varepsilon = 10^{-8}$,
$\delta = \eta_{\rm br} = +\infty$, and look-ahead with 10 candidates, in whose faces Gencan is run
until convergence. As in Algencan, the safeguarded multipliers are
the projections $\bar\lambda^{k+1} = P_{[\lambda_{\min}, \lambda_{\max}]^m}(\lambda^{k+1})$. In all the experiments,
$\varepsilon_k = 10^{-8}$ for all $k$. The subroutine \texttt{pcc\_telescopic} stops when
condition (a) holds with this tolerance and the swap test of \cite{bmtel} is satisfied, and the
CPS search described below checks condition (b). For Algorithm~\ref{alsec}.1, \texttt{pcc\_telescopic} also checks condition (c) when it stops. If (c) does not hold, the variables of the face that violate it are set to zero, provided that this decreases $L_k$, and the method goes on. Step~1 is implemented as follows. First, \texttt{pcc\_telescopic} is applied to $L_k$. When it stops, a CPS search is
performed over all the blocks with cardinality constraints, leaving and entering variables,
and admissible directions, with $\delta = +\infty$, $t_{\max} = 1$, $t_{\min} = 10^{-8}$,
$\alpha = 10^{-4}$, $\varepsilon = 10^{-8}$ in~(\ref{rhotol}), and $j_{\max} = 50$. The
restoration used in the other experiments is described in Section~\ref{macmpec}. If some CPS returns $w$ with
$L_k(w) < L_k(x^k) - 10^{-8}$, the best such $w$ is taken as the new initial point and
\texttt{pcc\_telescopic} is called again, in the face of $w$. The outer iterations use
$\tau = 0.5$ and $\gamma = 10$. The codes were compiled with gfortran 15.2 (option
\texttt{-O3}), and all the experiments were run on an Apple M1 Max computer (8 performance
cores, 64 GB of memory, macOS 13.7), with eight sequential runs in parallel. The codes, the
test problems and the complete results (final points, residuals, stationarity measures and
objective values of all the methods) are available at \url{https://www.ime.usp.br/~egbirgin/}. Tables with the results by instance of all the methods are given in the supplementary material.

\subsection{Examples} \label{numer}

In the examples of Sections~\ref{alsec} and~\ref{conv}, we used
$\rho_1 = 10$, $\bar\lambda^1 = 0$, $x^0 = (1, 0, 0.9)$ in
Examples~1 and~2 and $x^0 = (0.5, 0, 0.5, 0, 0, 0)$ in Example~4 (Section~\ref{nec}), and
stop when $\|h(x^k)\|_\infty \leq 10^{-8}$. Since
$\nabla L_k(x^k) = \nabla f(x^k) + \nabla h(x^k)\lambda^{k+1}$, condition (a) says that, at the
final iterate, the projected gradient of the Lagrangian $f + (\lambda^{k+1})^T h$ on the face
is at most $10^{-8}$. We used two settings of the multipliers. In the
augmented Lagrangian method (AL), $\lambda_{\min} = -10^{20}$ and $\lambda_{\max} = 10^{20}$. In
the quadratic penalty method (QP), $\lambda_{\min} = \lambda_{\max} = 0$, so that $\bar\lambda^k = 0$ for all $k$. To show that the CPS works independently of the penalty parameter, we also ran the
method with the CPS activated only from a prescribed outer iteration $k_{cps}$ on.

\paragraph{Example 1 (Section~\ref{alseccex}).}
By Proposition~\ref{propcex}, Algorithm~\ref{alsec}.1 (without CPS) may generate iterates that converge to
$x^* = 0$, which is not B-stationary, and this is what our implementation does. MPEC-LICQ
does not hold at $x^* = 0$ (all the variables are at their bounds, so $\F(x^*) = \emptyset$)
nor at $(0, 1, 0)$, where $\nabla h = 0$, so Theorem~\ref{teoconv} does not apply. Nevertheless,
the CPS works in this example. At $x^k = (s_k + s_k^3, 0, s_k)$, the restoration moves only $x_3$, with
$d_3(x^k) = s_k > 0$ and $\partial h(0, t, s_k)/\partial x_3 = s_k(1 + 2t) \neq 0$, and
returns $w = (0, t, x_3)$ with $(1 + 2t) x_3^2 = s_k^6$, at which $f(w) \approx -t/2$.

Table~\ref{tabres} summarizes the runs. Telescopic AL and CPS AL denote
Algorithms~\ref{alsec}.1 (original swaps only) and~\ref{cps}.3 (both kinds of swaps) with the
AL setting of the multipliers, and Telescopic QP and CPS QP denote the same algorithms with the
QP setting.
In the runs of the original method that failed, the iterates behaved as predicted by
Proposition~\ref{propcex} and Proposition~\ref{propcex2} (in Example~2, because of the
tolerance of the subproblems, the iterates are only approximations of the points of
Proposition~\ref{propcex2} when $\rho_k$ is large). In Examples~1 and~2, in all the runs with
CPS, exactly one CPS was performed, with $t = 1$, after which the Telescopic Method converged in the
new face to an approximation of the global minimizer $(0,1,0)$, where the outer iterations
stopped. In Example~1 with the CPS activated at $k = 40$ ($\rho_{40} = 10^{20}$), the CPS was
performed at $x^{40} = (3.0010 \times 10^{-3}, 0, 3.0009 \times 10^{-3})$, with
$f(x^{40}) = -6.1 \times 10^{-11}$, and returned
$w = (0, 1, 4.8 \times 10^{-8})$ with $f(w) = -0.5$ and $L_{40}(w) = -0.5$. Thus, the penalty
parameter played no role. The restoration tolerance~(\ref{rhotol}) was
$7 \times 10^{-15}$. In Example~1 with $k_{cps} = 1$, the CPS was rejected at $k = 1$
(the iterate was far from feasible, with $f(x^1) = -0.68$) and accepted at $k = 2$, at
$x^2 = (1.079, 0, 0.714)$ with $f(x^2) = -0.26$, returning $w = (0, 1, 0.210)$ with
$f(w) = -0.477$. In Example~4 (Section~\ref{nec}), started in the face of $y_1$ and $y_2$, the
original swaps of \texttt{pcc\_telescopic} found the global minimizer
$(0, 1, 0, 1, 0, 0)$ in the first outer iteration, with or without CPS.

\begin{table}[ht!]
\centering
\small
\setlength{\tabcolsep}{4pt}
\begin{tabular}{llrrrlrr}
\hline
Example & Method & $k_{cps}$ & Outer it. & CPS & Final point & $f$ & $\|h\|_\infty$ \\
\hline
\multirow{4}{*}{1} & Telescopic AL & -- & 21 & -- & $(4.9\mathrm{e}{\mbox{-}2}, 0, 4.9\mathrm{e}{\mbox{-}2})$ & $-4.3\mathrm{e}{\mbox{-}6}$ & $6.8\mathrm{e}{\mbox{-}9}$ \\
 & CPS AL & 1 & 2 & 1 & $(0, 1, 7.3\mathrm{e}{\mbox{-}20})$ & $-0.5$ & $8.0\mathrm{e}{\mbox{-}39}$ \\
 & CPS AL & 21 & 21 & 1 & $(0, 1, 4.2\mathrm{e}{\mbox{-}12})$ & $-0.5$ & $2.6\mathrm{e}{\mbox{-}23}$ \\
 & CPS AL & 40 & 40 & 1 & $(0, 1, 6.9\mathrm{e}{\mbox{-}9})$ & $-0.5$ & $7.1\mathrm{e}{\mbox{-}17}$ \\
\hline
\multirow{4}{*}{2} & Telescopic QP & -- & 11 & -- & $(2.3\mathrm{e}{\mbox{-}3}, 0, 2.3\mathrm{e}{\mbox{-}3})$ & $1.5\mathrm{e}{\mbox{-}12}$ & $2.3\mathrm{e}{\mbox{-}9}$ \\
 & CPS QP & 1 & 1 & 1 & $(0, 1, -1.4\mathrm{e}{\mbox{-}16})$ & $-0.25$ & $1.4\mathrm{e}{\mbox{-}16}$ \\
 & CPS QP & 11 & 11 & 1 & $(0, 1, 0)$ & $-0.25$ & $0$ \\
 & Telescopic AL & -- & 6 & -- & $(0, 1, 3.2\mathrm{e}{\mbox{-}10})$ & $-0.25$ & $3.2\mathrm{e}{\mbox{-}10}$ \\
\hline
\multirow{2}{*}{4} & Telescopic AL & -- & 1 & -- & $(0, 1, 0, 1, 0, 0)$ & $-2$ & $0$ \\
 & CPS AL & 1 & 1 & 0 & $(0, 1, 0, 1, 0, 0)$ & $-2$ & $0$ \\
\hline
\end{tabular}
\caption{Results. ``Outer it.'' is the number of outer iterations, ``CPS'' the number of
constraint-preserving swaps performed. With $k_{cps} > 1$, the stopping test is applied only
for $k \geq k_{cps}$. The penalty parameters at $k = 21$ and $k = 40$ of Example~1 are $10^{10}$
and $10^{20}$, and at $k = 11$ of Example~2 it is $10^6$.}
\label{tabres}
\end{table}

\subsection{The MacMPEC collection} \label{macmpec}

MacMPEC \cite{macmpec} is a collection of AMPL models of MPCCs with general constraints. We
used the model and data files distributed with MacMPEC.jl \cite{macmpecjl}, which also lists
the best objective value known for each instance (most of them obtained with the filter SQP
method of \cite{fl04}). Each of the 193 model/data combinations was translated by AMPL into a
\texttt{.nl} file, read by the AMPL Solver Library (ASL) \cite{gay}. Instance \texttt{ex9.1.2}
has a binary variable, which was reformulated as in Section~\ref{probl} for all the methods. Four of them
(\texttt{pack-rig2-16}, \texttt{pack-rig2c-16}, \texttt{pack-rig2-32} and
\texttt{pack-rig2c-32}), marked as infeasible in the list of \cite{macmpecjl}, are counted as failures of all
the methods in the columns ``Feas.'' and ``Stat.'' of Table~\ref{tabmac} (see item~\ref{itpsi} below). We call \emph{small} the 116 instances with at most 100 variables
in the AMPL model and \emph{large} the other 77 (the largest has 3436 variables).

With the option \texttt{cc\_simplify} of the ASL, every complementarity condition becomes
$0 \leq x_a \perp x_b \geq 0$ for two variables (the ASL adds auxiliary variables and
equality constraints when a side of the condition is an expression, and splits mixed
complementarity conditions), and the other constraints have the form
$\ell_i \leq c_i(x) \leq u_i$. Every inequality constraint is written as $c_i(x) - s_i = 0$,
with a slack $s_i \in [\ell_i, u_i]$. In this way the problem has the form considered in
Section~\ref{probl}, in which every complementarity pair is a block with $p = 1$ and the remaining
variables (original and auxiliary variables and slacks) form a block without cardinality
constraints. Functions and first and second derivatives are evaluated by the ASL. When an
evaluation fails (for example, outside the domain of a function), the value $10^{99}$ is
returned for the objective function of the subproblem, so that the trial point is rejected.

We compared variations of the augmented Lagrangian method with telescopic subproblems with two
methods from the literature. All of them use Gencan (the bound-constrained solver of
Algencan 4.0.0) for the subproblems:
\begin{itemize}[itemsep=2pt,parsep=0pt,topsep=4pt]
\item CPS AL: Algorithm~\ref{cps}.3 (original swaps and constraint-preserving
  swaps);
\item Telescopic AL: the same method with the original swaps only
  (Algorithm~\ref{alsec}.1);
\item Algencan: Algencan 4.0.0 \cite{bm20} applied to the nonlinear programming reformulation in
  which each pair is replaced by $x_a \geq 0$, $x_b \geq 0$, $x_a x_b \leq 0$ (augmented
  Lagrangian methods applied to this reformulation were analyzed in \cite{isu12,ass18});
\item Scholtes: the regularization method of Scholtes \cite{sch01}, with $x_a x_b \leq t$ for
  $t = 1, 10^{-1}, \dots, 10^{-8}$ and finally $t = 0$, each problem being solved by
  Algencan 4.0.0 from the solution and the multipliers of the previous one.
\end{itemize}
In CPS AL and Telescopic AL, the subproblems are solved by the implementation described at the beginning of
Section~\ref{numexp}, with three differences. First, at most $G = 10$ iterations of Gencan are
performed on each visited face before the next swap search (Algorithm~5.1 of \cite{bmtel}
performs a single one, and Remark~5.1 of \cite{bmtel} allows any number). Other values of $G$ are considered below. Second, there is no look-ahead. Third, we use simultaneous swaps.
When improving swaps are found in several blocks, the best swap of each block is performed at
the same time. If the resulting point is not better than the best single swap, halves of the
list are tried and, at last, the best single swap is performed. This is allowed by Remark~5.1
of \cite{bmtel}, by which any point of the new face with a functional value not larger than
$f_{best}$ may be chosen. The values of the
constraints are updated incrementally after a swap. The CPS search is performed when the
Telescopic Method stops at a point with $\|h(x)\|_\infty \leq 10^{-5}$, since condition (b) is
required in the proof of Theorem~\ref{teoconv} only for the iterates close to a feasible
limit point. Its restoration uses simplified Newton iterations, in which the matrix of Algorithm~\ref{cps}.1 is replaced by the matrix
$J_S D^2 J_S^T + \kappa_{\rm reg} I$, computed at the starting point $y$ and kept fixed in all the iterations. Here $J_S$ is the Jacobian of $h$ with respect to the moving
variables, $D = \mathrm{diag}(d_v)$ and
$\kappa_{\rm reg} = 10^{-12} \max\{1, \max_i (J_S D^2 J_S^T)_{ii}\}$, and the matrix is factorized once per search by MA57,
through the linear solver interface of Algencan. Lemma~\ref{lemrest} also holds for this version, as shown after its proof. The column of the leaving variable is
removed by the Sherman--Morrison formula. At most 50 iterations are performed, and the
restoration also stops without returning a point if an iterate leaves $\B$ or if the residual
$\|h(z^j) - \eta\|$ becomes larger than twice its initial value. The candidates (pair, entering variable) are sorted
by the first-order prediction of the variation of $L_k$, only candidates with negative
prediction are tried, and the first CPS that returns a point is accepted ($\alpha = 10^{-4}$,
$t_{\max} = \max\{1, \|x\|_\infty\}$, $t_{\min} = 10^{-2}$, $\varepsilon = 10^{-8}$). Since $t_{\min}$ and $\varepsilon$ are fixed in the implementation, Theorem~\ref{teoconv} does not apply literally, and the final points are verified a posteriori, as described below.

The initial point of CPS AL and Telescopic AL is computed by a light Phase~I, in which Algencan, with at most 10 outer
iterations and tolerances $10^{-4}$, is applied to the reformulation used by Algencan from the initial
point of the AMPL model. If it stops at a point whose sup-norm of the constraint violation is
larger than $10^{-4}$, it is applied again from a random perturbation of the AMPL point
(fixed seed), until such a point is found or the time limit is reached. The best point found,
in which the smaller variable of each pair is set to zero and the slacks are set to the
projections of the constraint values, is the initial point of the method, and its
multipliers are the initial safeguarded multipliers $\bar\lambda^1$. The outer iterations use
the parameters of Algencan, $\tau = 0.5$, $\gamma = 10$, $\lambda_{\min} = -10^{20}$,
$\lambda_{\max} = 10^{20}$ and $\rho_1 = \max\{10^{-8}, \min\{10 \max\{1, |f(x^0)|\} /
\max\{1, \frac12 \|h(x^0)\|^2\}, 10^8\}\}$. There are at most 50 outer iterations, and the
method stops when $\|h(x^k)\|_\infty \leq 10^{-8}$ and the subproblem was solved
with the required tolerance, when $\|h(x^k)\|_\infty \leq 10^{-8}$ at three consecutive
iterations in which the subproblem solver did not reach the required tolerance (as in
Algencan), or when $\rho_k > 10^{20}$. If a subproblem is unbounded below, its solution is discarded and the
subproblem is solved again from the previous iterate with $\rho_k$ multiplied by $\gamma$.
To assess the CPS independently of Phase~I, CPS AL and Telescopic AL were also run without Phase~I, from the
initial point of the AMPL model projected onto the bounds (with the same treatment of pairs and slacks), in variants denoted by CPS AL$_0$ and Telescopic AL$_0$.

In Algencan and Scholtes, the feasibility, complementarity and optimality tolerances of Algencan were set
to $10^{-8}$, with at most 50 outer iterations and the default values of the other
parameters, from the initial point of the AMPL model projected onto the bounds. Each run had a time limit of one hour. In CPS AL and Telescopic AL,
the limit is checked by the code between subproblems (and between Phase~I attempts), and the best point found is returned. A run still active after 3800 seconds was interrupted and
counted as a failure. Algencan and Scholtes were interrupted after one hour.

The final points were checked independently of the methods. The original AMPL model (without
simplifications) was read by the ASL, and we computed at the final point the maximum violation
of the bounds and the constraints and the complementarity residual. For a condition
$\ell \leq x_j \leq u$ complementary to $c_i(x)$, the residual is
$|x_j - P_{[\ell, u]}(x_j - c_i(x))|$, which is $\min\{|x_j - \ell|, |c_i(x)|\}$ when
$u = +\infty$. A point is \emph{feasible} when both quantities are at most $10^{-6}$. At
feasible points, stationarity was checked in the simplified problem by solving, with the
bounded-variable least-squares method, $\nabla f(x) = \sum_{i \in \mathcal{A}} \mu_i
\nabla c_i(x) + \sum_{j \in \mathcal{V}} \pi_j e^j$, where $\mathcal{A}$ and $\mathcal{V}$
are the active constraints and bounds (relative tolerance $10^{-6}$), with the signs given by
the activity (free multipliers for equalities and for the null variable of a nondegenerate
pair). The point is \emph{S-stationary} when the relative residual is at most $10^{-6}$ with
nonnegative multipliers for both variables of every biactive pair. Otherwise, if there are at
most 10 biactive pairs, all the branches were checked (in each biactive pair, the multiplier
of one of the variables is free and the other is nonnegative), and the point is
\emph{B-stationary} when the residual is at most $10^{-6}$ in every branch (under MPEC-LICQ,
this is B-stationarity). A run is a \emph{success} when its final point is feasible and S- or
B-stationary. At infeasible final points, we checked the stationarity of the infeasibility
measure $\Psi$ (Proposition~\ref{propinf}), with the slacks at their best values. The point is
\emph{stationary for $\Psi$} when the projected gradient of $\Psi$ on the face of the point (for
Algencan and Scholtes, on the box, with the complementarity constraints written as
$x_a x_b \leq 0$) is at most $10^{-6} \|\nabla \Psi\|_\infty$. The relative test excludes nearly
feasible points, at which $\nabla \Psi$ is tiny. Finally, the objective value of a feasible point was compared with the best known
value $f_{best}$, and it \emph{reaches} $f_{best}$ when it is at most
$f_{best} + 10^{-3} \max\{1, |f_{best}|\}$ (for minimization).

Table~\ref{tabmac} summarizes the results.

\begin{table}[ht!]
\centering
\small
\setlength{\tabcolsep}{4pt}
\begin{tabular}{l rrrr rrr rrr}
\hline
 & \multicolumn{4}{c}{All (193)} & \multicolumn{3}{c}{Small (116)} & \multicolumn{3}{c}{Large (77)} \\
Method & Feas. & Stat. & $\Psi$ & Best & Feas. & Stat. & Best & Feas. & Stat. & Best \\
\hline
CPS AL               & 185 & 185 & 2 & 165 & 115 & 115 & 110 & 70 & 70 & 55 \\
Telescopic AL        & 185 & 185 & 2 & 165 & 115 & 115 & 110 & 70 & 70 & 55 \\
CPS AL$_0$           & 175 & 174 & 13 & 148 & 108 & 108 & 99 & 67 & 66 & 49 \\
Telescopic AL$_0$    & 175 & 174 & 13 & 145 & 108 & 108 & 97 & 67 & 66 & 48 \\
Algencan             & 179 & 171 & 8 & 156 & 110 & 106 & 102 & 69 & 65 & 54 \\
Scholtes             & 168 & 165 & 8 & 150 & 109 & 106 & 102 & 59 & 59 & 48 \\
\hline
\end{tabular}
\caption{MacMPEC: number of instances in which the final point is feasible (``Feas.''), is
feasible and S- or B-stationary (``Stat.''), and reaches the best known value (``Best''), with $G = 10$
in CPS AL, Telescopic AL, CPS AL$_0$ and Telescopic AL$_0$. ``$\Psi$'' is the number of instances in which the
final point is infeasible and stationary for $\Psi$. The four infeasible instances are counted as failures in ``Feas.'' and ``Stat.''.}
\label{tabmac}
\end{table}

We highlight the following facts.
\begin{enumerate}[itemsep=2pt,parsep=0pt,topsep=4pt]
\item Stationarity. All the feasible final points of CPS AL and Telescopic AL, and all but one of CPS AL$_0$ and of Telescopic AL$_0$ (both in \texttt{bem-milanc30-s}), were certified S- or B-stationary, against 171 of the 179 feasible points of Algencan and 165 of the 168 of Scholtes. Four points of CPS AL are B-stationary with biactive pairs (\texttt{ex9.2.2}, \texttt{qpec2}, \texttt{ralph1}, \texttt{scholtes4}). In these instances, the complementarity residuals of Algencan are between $3.7 \times 10^{-5}$ and $7.5 \times 10^{-5}$, a typical consequence of the degeneracy of the reformulation $x_a x_b \leq 0$, and, with tolerance $10^{-4}$ for this residual, Algencan and Scholtes would be feasible in 184 and 182 instances. In \texttt{bem-milanc30-s}, CPS AL, Telescopic AL and Algencan found S-stationary points with objective value $0.0944$, much smaller than the value $83.28$ of the list. In the 179 instances in which both CPS AL and Algencan found feasible points, the objective values coincide (up to $10^{-3}$) in 168, CPS AL is better in 7 and Algencan in 4.
\item \label{itpsi} Failures and infeasible points. The failures of CPS AL are the four infeasible instances and \texttt{TSC-7}, \texttt{TSC-8}, \texttt{TSC-12} and \texttt{taxmcp}, in which Phase~I did not find a point with violation at most $10^{-4}$ within the time limit. In the four infeasible instances, all the runs that stopped before the time limit found the same point, stationary for $\Psi$, with $\|h\|_\infty = 1.0 \times 10^{-4}$ (\texttt{pack-rig2-16}, \texttt{pack-rig2c-16}) and $2.3 \times 10^{-4}$ (\texttt{pack-rig2-32}, \texttt{pack-rig2c-32}), which supports their infeasibility. Without Phase~I, CPS AL$_0$ and Telescopic AL$_0$ also stopped with $\rho_k > 10^{20}$ at infeasible points stationary for $\Psi$ in nine other instances, as predicted by Proposition~\ref{propinf}.
\item Phase~I and CPS. With Phase~I, no CPS was accepted, and CPS AL and Telescopic AL returned the same points (the CPS search found candidates in 113 instances, but the Telescopic Method had already stopped at S-stationary points). Without Phase~I, the CPS improved the final point of Telescopic AL$_0$ in five instances, namely \texttt{ex9.2.5} ($f = 9 \to 5$), \texttt{kth3} ($1 \to 0.5$) and \texttt{qpec-100-2} ($-4.2517 \to -6.5907$), reaching the best known values, and \texttt{qpec-200-4} ($-3.2184 \to -4.9845$) and \texttt{qpec-100-1} ($0.2828 \to 0.2782$). Phase~I increased the number of successes from 174 to 185. In \texttt{ex9.1.4}, for example, the first subproblem without Phase~I, with a small penalty parameter, chooses branches of the pairs with no feasible point near the limit, and the CPS, which preserve the values of the constraints, cannot correct this.
\item Time. In the 167 instances in which all four methods found feasible points, the total times were 814 s (CPS AL), 730 s (Telescopic AL), 599 s (Algencan) and 1321 s (Scholtes), so that the CPS searches increased the time by about 12\%. The four runs of CPS AL that reached the time limit account for 14401 of its 18906 seconds. CPS AL$_0$ and Telescopic AL$_0$ took 4099 seconds (2569 in \texttt{bem-milanc30-s}) and 1782 seconds.
\item Iterations of Gencan per visited face. With $G = 1$, the swap search is performed after every iteration of Gencan, and its cost dominates in the instances with many pairs (1835 seconds in \texttt{incid-set1c-32}, with 961 pairs, against 52 seconds with $G = 10$). CPS AL had 183, 184, 185 and 184 successes with $G = 1$, $5$, $10$ and $50000$, and the total times in the 183 instances solved by the four versions were 4924, 1273, 1102 and 1263 seconds.
\end{enumerate}

Nurkanovi\'c and Leyffer \cite{nl26} solved 191 MacMPEC problems (translated to CasADi, with IPOPT as NLP solver) and counted as a success the computation of a B-stationary point (certified a posteriori for the methods other than MPECopt). Their highest success rate, 94.24\%, was obtained by MPECopt with regularization in Phase~I, followed by 93.2\% for a regularization method, 92.14\% for an $\ell_\infty$ penalty method, 90.0\% for IPOPT applied to the NLP reformulation and 72.52\% for a mixed-integer reformulation solved by Bonmin. With the same criterion (feasible and S- or B-stationary), the success rates of our experiments in their list of 191 problems are 97.4\% (CPS AL and Telescopic AL), 91.6\% (CPS AL$_0$ and Telescopic AL$_0$), 90.1\% (Algencan) and 85.9\% (Scholtes). Their list contains 14 repeated instances under different names (for example, \texttt{pack-rig-k} and \texttt{pack-rig2p-k}), and the rate of CPS AL in its 177 distinct instances is 97.2\%. The comparison is only indicative, because the formulations, NLP solvers, tolerances, time limits and computers are different. Infeasible final points are counted as failures, as in \cite{nl26}, although MPECopt also has a convergence result for infeasible points (Theorem~4.2 of \cite{nl26}).

In summary, as predicted by the theory, the feasible final points of CPS AL were S- or B-stationary, and CPS AL and Telescopic AL succeeded in 97.4\% of the problems of \cite{nl26}. With Phase~I, the original swaps were already sufficient in this collection, and the CPS added the guarantee of Theorem~\ref{teoconv} at a small cost. When the initial point does not lie in a good branch, the CPS improve the points at which the original swaps stop (CPS AL$_0$ versus Telescopic AL$_0$). The remaining failures are related to Phase~I, not to the CPS.

\subsection{Partitioned cardinality constraints} \label{cardexp}

In this subsection the blocks have more than two variables and, in two of the three families,
more than one non-null variable is allowed in each block. We used the implementation of
Section~\ref{macmpec}, extended to general blocks, with the same variants. CPS AL and Telescopic AL
start from the point given by Phase~I, that is, by
Algencan applied to the problem without cardinality constraints, followed by keeping the
largest entries of each block. CPS AL$_0$ and Telescopic AL$_0$ are the same methods without Phase~I. When the
variables are non-negative, we also ran Algencan and the regularization method of
Scholtes applied to the reformulation of \cite{bks}. The time limit was 600 seconds. Every final point was verified
independently for feasibility with tolerance $10^{-6}$, cardinality (entries with absolute value
at most $10^{-8}$ are considered null), and stationarity in the sense of
Definition~\ref{defstat} for the Lagrangian, with tolerance $10^{-6}$ (under PCC-LICQ, B-stationarity). We then say that the final point is \emph{verified}. For each instance, we
report the relative difference between the objective function value of each method and the
best value found by all the methods (the gap). In the three families described below, the
constraint-preserving swaps improved the original swaps, never worsening the final point, and the resulting method was better than the methods applied to the
reformulation of \cite{bks} and close to the reference values computed by CPLEX.

\subsubsection{Portfolio selection}

\paragraph{Mean-variance portfolios.}
The problem is
\[
\Minimize x^T Q x \mbox{ subject to } \bar r^T x \geq r_{\min}, \; \textstyle\sum_{j=1}^n x_j = 1, \; 0 \leq x \leq u, \; \|x\|_0 \leq s,
\]
where $Q$ is the covariance matrix of the returns of $n$ assets, $\bar r$ is the vector of mean
returns, $r_{\min}$ is the minimum required return and $u$ is the vector of upper bounds (one
block, with $p_1 = s$). We used the 90 instances of Frangioni and Gentile \cite{fg06} (available at
\url{https://commalab.di.unipi.it/datasets/MV/}), with $n = 200$, $300$ and $400$,
ignoring the minimum buy-in constraints, as in \cite{krs21,jkmw23}, and $s \in \{5, 10, 20\}$
(270 problems). The initial point
of CPS AL$_0$ and Telescopic AL$_0$ is $x = 0$. As a reference, we used CPLEX~22.2 applied to the mixed-integer
formulation, with one thread and a time limit of 600 seconds, which it reached in 266 problems, so that its values are not certified to be optimal. Table~\ref{tabmv} shows, for
each method, the number of verified final points and the mean gap. All the final points of Telescopic AL, Telescopic AL$_0$, Scholtes and CPLEX
were verified.

\begin{table}[ht!]
\centering
\begin{tabular}{lccc}
\hline
 & $s = 5$ & $s = 10$ & $s = 20$ \\
\hline
CPS AL     & 90 / 0.35\% & 90 / 0.32\% & 88 / 0.07\% \\
CPS AL$_0$    & 90 / 0.48\% & 90 / 0.24\% & 84 / 0.25\% \\
Telescopic AL     & 90 / 5.65\% & 90 / 4.51\% & 90 / 4.05\% \\
Telescopic AL$_0$    & 90 / 12.9\% & 90 / 14.7\% & 90 / 11.2\% \\
Algencan     & 81 / 10.1\% & 76 / 7.44\% & 71 / 5.99\% \\
Scholtes     & 90 / 3.66\% & 90 / 3.59\% & 90 / 3.54\% \\
CPLEX & 90 / 0.01\% & 90 / 0.00\% & 90 / 0.01\% \\
\hline
\end{tabular}
\caption{Mean-variance portfolios (90 problems for each $s$): number of verified final points
and mean gap.}
\label{tabmv}
\end{table}

Comparing CPS AL with Telescopic AL problem by problem, CPS AL found a better point in 194 problems and the same point
in 74, and never a worse one. CPS AL was better than Scholtes in 245 problems and worse in 13,
and better than Algencan in 188 problems and worse in 29. The eight failures of CPS AL and CPS AL$_0$ occurred for $n = 400$ and $s = 20$, by the
time limit. CPLEX found better points than CPS AL in 157 problems. In all of them, the two points have
different sets of non-null variables, that is, CPS AL stops at another local minimizer. The
median relative difference was $0.2\%$ in these problems and $0.03\%$ over all the problems.
The constraint-preserving swaps were accepted in most problems, and they are responsible for
the difference between CPS AL and Telescopic AL. They make the method slower. The total time of
CPS AL was 11580 seconds, against 2962 for Telescopic AL, because there are $s(n - s)$
candidate swaps.

The instances with $n = 200$ were also used in \cite{jkmw23}, with the same cardinalities. A
direct comparison is not possible, because the code of \cite{jkmw23} is not available and the
results are reported only graphically, relative to the values obtained by CPLEX with a time
limit of 30 minutes. According to \cite{jkmw23}, the values obtained by their augmented
Lagrangian method are usually not too far from those of CPLEX, and those of the boosted
version are usually very close to them. In these 90 problems, the values obtained by CPS AL were
equal to or better than those of CPLEX (with 600 seconds) in 41 problems. The median relative
difference was $0.02\%$ (mean $0.25\%$, maximum $4.8\%$).

\paragraph{Portfolios with a risk constraint and sector cardinality.}
The problem is
\[
\Maximize \bar r^T x \mbox{ subject to } x^T C x \leq \sigma_{\max}^2, \; \textstyle\sum_{j=1}^n x_j = 1, \; 0 \leq x \leq u, \; \|\vx_\nu\|_0 \leq p, \; \nu = 1, \ldots, N,
\]
where $\bar r$ and $C$ are the mean returns and the covariance matrix of the assets, and each
block $\vx_\nu$ is a sector, in which at most $p$ assets may be chosen. We used the data of the
S\&P~500 index used in \cite{bmtel} ($n = 497$ assets in $N = 11$ sectors, with $\bar r$ and
$C$ computed from 10{,}000 scenarios), with $p \in \{1, 2, 3\}$, $u \in \{0.2, 1\}$ (the same
bound for all the assets) and $\sigma_{\max} \in \{0.15, 0.20, 0.25\}$ (18 problems). We
excluded $\sigma_{\max} = 0.10$, for which no method, CPLEX included, found a feasible point.
The reference
values were computed by CPLEX applied to the mixed-integer quadratically constrained
formulation (600 seconds). All the final points of CPS AL and Telescopic AL were verified,
as were 17 of the 18 points of CPS AL$_0$, Telescopic AL$_0$ and Scholtes, and only
6 of Algencan. Because of its integrality tolerances, the points of CPLEX were checked only for feasibility and cardinality, and 17 of them passed. Table~\ref{tabsp} shows the number of verified points and the
mean gap over them (the best value is the best one among the verified points of all the
methods, CPLEX included).

\begin{table}[ht!]
\centering
\begin{tabular}{lcc}
\hline
 & Verified points & Mean gap (\%) \\
\hline
CPS AL                 & 18 & 0.15 \\
CPS AL$_0$             & 17 & 0.08 \\
Telescopic AL          & 18 & 0.50 \\
Telescopic AL$_0$      & 17 & 3.8 \\
Algencan               & 6 & 6.8 \\
Scholtes               & 17 & 0.94 \\
CPLEX                  & 17$^*$ & 0.02 \\
\hline
\end{tabular}
\caption{Portfolios with a risk constraint and sector cardinality (18 problems): number of verified final points and mean gap. $^*$Feasibility and cardinality only.}
\label{tabsp}
\end{table}

CPS AL was better than Telescopic AL in 9 problems and equal in the other 9. It was better
than Scholtes in 8 problems, equal in 8 and worse in 1. CPS AL and Telescopic AL took about 4 and 3
seconds per problem.

\subsubsection{Sparse principal component analysis}

The problem is
\[
\Minimize -x^T A x \mbox{ subject to } \|x\|^2 = 1, \; -1 \leq x \leq 1, \; \|x\|_0 \leq p,
\]
where $A \in \R^{n \times n}$ is a covariance matrix, so that $x^T A x$ is the variance
explained by the direction $x$, and at most $p$ entries of $x$ may be non-null (one block). We
used the sample covariance matrices $A$ of $n$ observations of a normal vector in $\R^n$ with
covariance $I + \theta v v^T$, where $v$ has 10 non-null entries $\pm 1/\sqrt{10}$, for
$n \in \{50, 100, 200\}$, $\theta \in \{1, 3\}$ and five random seeds, and $p \in \{5, 10, 20\}$
(90 problems). Since $x = 0$ is stationary for the augmented Lagrangian, CPS AL$_0$ and Telescopic AL$_0$ start from
the canonical vector of the largest diagonal entry of $A$. We compared with the truncated
power method (TP) of Yuan and Zhang \cite{yz13}, started from the leading eigenvector
truncated to its $p$ largest entries. All the final points of all the methods were verified. Table~\ref{tabspca} shows the mean gap of each method and the number of
problems in which it found the best value.

\begin{table}[ht!]
\centering
\begin{tabular}{lcc}
\hline
 & Mean gap (\%) & Best value \\
\hline
CPS AL                 & 0.30 & 70 \\
CPS AL$_0$             & 1.5 & 55 \\
Telescopic AL          & 0.52 & 60 \\
Telescopic AL$_0$      & 2.0 & 45 \\
TP                     & 3.3 & 21 \\
\hline
\end{tabular}
\caption{Sparse principal component analysis (90 problems): mean gap and number of problems in which each method found the best value.}
\label{tabspca}
\end{table}

CPS AL was better than Telescopic AL in 13 problems and
equal in the other 77, and better than TP in 65 and worse in 2. The constraint-preserving
swaps were accepted 22 times. All the runs took less than one second.

\section{Conclusions} \label{concl}

We considered the augmented Lagrangian method of Algencan in which the bounds and the
partitioned cardinality constraints (in particular, complementarity constraints) are kept in
the subproblems, which are solved by the Telescopic Method. The B-stationarity that the
Telescopic Method guarantees for the subproblems is inherited by the limit points when the
penalty parameters remain bounded, but not in general. The reason is that the variation of the penalty term
along a swap grows with the penalty parameter, so that the swaps of the subproblems may stop
detecting descent directions of the original problem. Constraint-preserving swaps remove the dependence on the penalty parameter. With them, every feasible limit point at which PCC-LICQ (MPEC-LICQ for MPCCs)
holds is B-stationary (S-stationary for MPCCs), whatever the behavior of the penalty
parameters and of the safeguarded multipliers, which includes the quadratic penalty method.
Without them, if the stopping criterion of the subproblems also includes a derivative test outside the current face, the feasible limit points are still M-stationary for MPCCs under MPEC-LICQ, as for other augmented Lagrangian methods that keep the complementarity constraints in the subproblems. We also showed that the modification is necessary under MPEC-LICQ, even with exact solutions
of the subproblems and multipliers in the interior of the safeguarding box, and that this assumption
cannot be replaced by MPEC-MFCQ, because descent may require simultaneous swaps in several
blocks.

In the MacMPEC collection, with a light Phase I, the original swaps already stopped at S- or
B-stationary points, and no constraint-preserving swap was accepted. Without Phase I, the
constraint-preserving swaps improved the final points of the original swaps in five instances,
three of which reached the best known value. The method
found certified S- or B-stationary points in 185 of the 193 instances, more than the other
methods tested, and its failures are related to Phase I, not to the swaps.
With cardinality constraints (Section~\ref{cardexp}), on the other hand, the
constraint-preserving swaps were often accepted and improved the final points of the original
swaps, at the price of a larger number of candidate swaps.

\section*{Declaration of generative AI use}
During the preparation of this work, the authors used Claude (Anthropic) to assist in drafting the
manuscript, developing and checking the proofs, implementing the codes and running the numerical
experiments. The authors reviewed and edited all the content and take full responsibility for the
content of the publication.


\begin{thebibliography}{99}

\bibitem{abms} R. Andreani, E. G. Birgin, J. M. Mart\'{\i}nez, and M. L. Schuverdt, On augmented
  Lagrangian methods with general lower-level constraints, {\it SIAM Journal on Optimization}
  18, pp.~1286--1309, 2008.

\bibitem{ars25} R. Andreani, M. da Rosa, and L. D. Secchin, On the boundedness of
  multipliers in augmented Lagrangian methods for mathematical programs with complementarity
  constraints, Optimization Online, 2025.

\bibitem{ahm11} R. Andreani, G. Haeser, and J. M. Mart\'{\i}nez, On sequential optimality
  conditions for smooth constrained optimization, {\it Optimization} 60, pp.~627--641, 2011.

\bibitem{ahss19} R. Andreani, G. Haeser, L. D. Secchin, and P. J. S. Silva, New sequential
  optimality conditions for mathematical programs with complementarity constraints and
  algorithmic consequences, {\it SIAM Journal on Optimization} 29, pp.~3201--3230, 2019.

\bibitem{ass18} R. Andreani, L. D. Secchin, and P. J. S. Silva, Convergence properties of a second
  order augmented Lagrangian method for mathematical programs with complementarity constraints,
  {\it SIAM Journal on Optimization} 28, pp.~2574--2600, 2018.

\bibitem{be13} A. Beck and Y. C. Eldar, Sparsity constrained nonlinear optimization: optimality
  conditions and algorithms, {\it SIAM Journal on Optimization} 23, pp.~1480--1509, 2013.

\bibitem{bi66} A. Ben-Israel, A Newton--Raphson method for the solution of systems of
  equations, {\it Journal of Mathematical Analysis and Applications} 15, pp.~243--252, 1966.

\bibitem{bmtel} E. G. Birgin and J. M. Mart\'{\i}nez, A Telescopic Method with guaranteed
  B-stationarity for partitioned cardinality and complementarity constraints, arXiv:2611.00897
  [math.OC], 2026.

\bibitem{bm20} E. G. Birgin and J. M. Mart\'{\i}nez, Complexity and performance of an augmented
  Lagrangian algorithm, {\it Optimization Methods and Software} 35, pp.~885--920, 2020.

\bibitem{gencan} E. G. Birgin and J. M. Mart\'{\i}nez, Large-scale active-set box-constrained
  optimization method with spectral projected gradients, {\it Computational Optimization and
  Applications} 23, pp.~101--125, 2002.

\bibitem{bmbook} E. G. Birgin and J. M. Mart\'{\i}nez, {\it Practical Augmented Lagrangian Methods
  for Constrained Optimization}, SIAM, Philadelphia, 2014.

\bibitem{bks} O. P. Burdakov, C. Kanzow, and A. Schwartz, Mathematical programs with cardinality
  constraints: reformulation by complementarity-type conditions and a regularization method,
  {\it SIAM Journal on Optimization} 26, pp.~397--425, 2016.

\bibitem{ds96} J. E. Dennis, Jr. and R. B. Schnabel, {\it Numerical Methods for Unconstrained
  Optimization and Nonlinear Equations}, Classics in Applied Mathematics 16, SIAM,
  Philadelphia, 1996.

\bibitem{fl04} R. Fletcher and S. Leyffer, Solving mathematical programs with complementarity
  constraints as nonlinear programs, {\it Optimization Methods and Software} 19, pp.~15--40,
  2004.

\bibitem{flrs06} R. Fletcher, S. Leyffer, D. Ralph, and S. Scholtes, Local convergence of SQP
  methods for mathematical programs with equilibrium constraints, {\it SIAM Journal on
  Optimization} 17, pp.~259--286, 2006.

\bibitem{fg06} A. Frangioni and C. Gentile, SDP diagonalizations and perspective cuts for a class
  of nonseparable MIQP, {\it Operations Research Letters} 35, pp.~181--185, 2007.

\bibitem{gay} D. M. Gay, Hooking your solver to AMPL, Technical Report 97-4-06, Computing Sciences
  Research Center, Bell Laboratories, 1997.

\bibitem{gd22} L. Guo and Z. Deng, A new augmented Lagrangian method for
  MPCCs---theoretical and numerical comparison with existing augmented Lagrangian methods,
  {\it Mathematics of Operations Research} 47, pp.~1229--1246, 2022.

\bibitem{hks13} T. Hoheisel, C. Kanzow, and A. Schwartz, Theoretical and numerical comparison of
  relaxation methods for mathematical programs with complementarity constraints,
  {\it Mathematical Programming} 137, pp.~257--288, 2013.

\bibitem{isu12} A. F. Izmailov, M. V. Solodov, and E. I. Uskov, Global convergence of augmented
  Lagrangian methods applied to optimization problems with degenerate constraints, including
  problems with complementarity constraints, {\it SIAM Journal on Optimization} 22,
  pp.~1579--1606, 2012.

\bibitem{jkmw23} X. Jia, C. Kanzow, P. Mehlitz, and G. Wachsmuth, An augmented Lagrangian method for
  optimization problems with structured geometric constraints, {\it Mathematical Programming}
  199, pp.~1365--1415, 2023.

\bibitem{ks15} C. Kanzow and A. Schwartz, The price of inexactness: convergence properties of
  relaxation methods for mathematical programs with complementarity constraints revisited,
  {\it Mathematics of Operations Research} 40, pp.~253--275, 2015.

\bibitem{krs21} C. Kanzow, A. B. Raharja, and A. Schwartz, An augmented Lagrangian method for
  cardinality-constrained optimization problems, {\it Journal of Optimization Theory and
  Applications} 189, pp.~793--813, 2021.

\bibitem{krs21b} C. Kanzow, A. B. Raharja, and A. Schwartz, Sequential optimality conditions for
  cardinality-constrained optimization problems with applications, {\it Computational
  Optimization and Applications} 80, pp.~185--211, 2021.

\bibitem{kllm} C. Kirches, J. Larson, S. Leyffer, and P. Manns, Sequential linearization method for
  bound-constrained mathematical programs with complementarity constraints, {\it SIAM Journal on
  Optimization} 32, pp.~75--99, 2022.

\bibitem{macmpecjl} B. Legat, MacMPEC.jl, \url{https://github.com/blegat/MacMPEC.jl}.

\bibitem{llcn06} S. Leyffer, G. L\'opez-Calva, and J. Nocedal, Interior methods for mathematical
  programs with complementarity constraints, {\it SIAM Journal on Optimization} 17,
  pp.~52--77, 2006.

\bibitem{macmpec} S. Leyffer, MacMPEC: AMPL collection of MPECs, Argonne National Laboratory,
  2000, \url{https://wiki.mcs.anl.gov/leyffer/index.php/MacMPEC}.

\bibitem{lf06} G.-H. Lin and M. Fukushima, Hybrid approach with active set identification for
  mathematical programs with complementarity constraints, {\it Journal of Optimization Theory
  and Applications} 128, pp.~1--28, 2006.

\bibitem{nl26} A. Nurkanovi\'c and S. Leyffer, A globally convergent method for computing
  B-stationary points of mathematical programs with equilibrium constraints, {\it Computational
  Optimization and Applications}, 2026, DOI 10.1007/s10589-026-00832-x.

\bibitem{or70} J. M. Ortega and W. C. Rheinboldt, {\it Iterative Solution of Nonlinear
  Equations in Several Variables}, Academic Press, New York, 1970.

\bibitem{ss00} H. Scheel and S. Scholtes, Mathematical programs with complementarity constraints:
  stationarity, optimality, and sensitivity, {\it Mathematics of Operations Research} 25,
  pp.~1--22, 2000.

\bibitem{sch01} S. Scholtes, Convergence properties of a regularization scheme for mathematical
  programs with complementarity constraints, {\it SIAM Journal on Optimization} 11,
  pp.~918--936, 2001.

\bibitem{yz13} X.-T. Yuan and T. Zhang, Truncated power method for sparse eigenvalue problems,
  {\it Journal of Machine Learning Research} 14, pp.~899--925, 2013.

\end{thebibliography}
\end{document}